\PassOptionsToPackage{unicode}{hyperref}
\PassOptionsToPackage{hyphens}{url}

\documentclass[preprint,review,11pt]{elsarticle}

\journal{}

\usepackage{iftex}
\ifPDFTeX
  \usepackage[T1]{fontenc}
  \usepackage[utf8]{inputenc}
  \usepackage{textcomp}
  \usepackage{mathptmx} 
\else
  \usepackage{unicode-math}
  \defaultfontfeatures{Scale=MatchLowercase}
  \defaultfontfeatures[\rmfamily]{Ligatures=TeX,Scale=1}
\fi

\usepackage{fullpage,times}
\usepackage{amsmath,amssymb}
\usepackage{bm}
\usepackage{mathtools}
\usepackage{graphicx}
\usepackage{longtable}
\usepackage{booktabs}
\usepackage{array}
\usepackage{calc}
\usepackage{multirow}
\usepackage{multicol}
\usepackage{makecell}
\usepackage{enumitem}
\usepackage{placeins}
\usepackage{siunitx}
\usepackage{xspace}

\usepackage{setspace}
\usepackage{etoolbox}

\usepackage{algorithm}
\usepackage[noend]{algorithmic}
\AtBeginEnvironment{algorithmic}{\setstretch{1.1}}

\usepackage{subcaption}
\usepackage{caption}
\usepackage{tabularx}
\newcolumntype{L}[1]{>{\raggedright\arraybackslash}p{#1}}
\newcolumntype{Y}{>{\raggedright\arraybackslash}X}

\usepackage[table,xcdraw,x11names]{xcolor}
\usepackage{xurl}
\usepackage{hyperref}
\hypersetup{
  colorlinks=true,
  linkcolor=Blue3,
  filecolor=Blue3,
  citecolor=Blue3,
  urlcolor=Blue3,
  pdfstartview={FitH},
  pdfauthor={Jilin Song, Amer Shalaby, Merve Bodur},
  pdftitle={Approximate Dynamic Programming for Real-time Assignment of Extraboard Transit Operators},
  pdfsubject={Public transit, crew scheduling, Markov decision process, dynamic programming},
  plainpages=false
}

\usepackage[numbers]{natbib}

\makeatletter
\patchcmd\longtable{\par}{\if@noskipsec\mbox{}\fi\par}{}{}
\makeatother

\IfFileExists{footnotehyper.sty}{\usepackage{footnotehyper}}{\usepackage{footnote}}
\makesavenoteenv{longtable}

\makeatletter
\def\maxwidth{%
  \ifdim\Gin@nat@width>\linewidth
    \linewidth
  \else
    \Gin@nat@width
  \fi
}
\def\maxheight{%
  \ifdim\Gin@nat@height>\textheight
    \textheight
  \else
    \Gin@nat@height
  \fi
}
\makeatother

\setkeys{Gin}{width=\maxwidth,height=\maxheight,keepaspectratio}

\makeatletter
\def\fps@figure{htbp}
\makeatother

\ifLuaTeX
  \usepackage{selnolig}
\fi

\makeatletter
\def\ps@pprintTitle{%
  \let\@oddhead\@empty
  \let\@evenhead\@empty
  \def\@oddfoot{\reset@font\hfil\thepage\hfil}
  \let\@evenfoot\@oddfoot
}
\makeatother

\newcommand{\delmin}{\ensuremath{\delta_{\texttt{min}}}}
\newcommand{\kfirst}{\ensuremath{k^\texttt{one}_h}}
\newcommand{\jtn}{\ensuremath{\mathcal{J}_t^\texttt{now}}}
\newcommand{\itn}{\ensuremath{\mathcal{I}_t^\texttt{now}}}
\newcommand{\dtn}{\ensuremath{\mathcal{D}_t^\texttt{now}}}
\newcommand{\tfrj}{\ensuremath{\tau_t(j)}}
\newcommand{\phiti}{\ensuremath{\phi_t(i)}}

\newcommand{\decspace}{\ensuremath{\mathcal{X}^t(S_t)}}
\newcommand{\xtjk}{\ensuremath{x^t_{jk}}}
\newcommand{\utik}{\ensuremath{u^t_{ik}}}
\newcommand{\xvec}{\ensuremath{\bm{x}^t}}

\newcommand{\stpo}{\ensuremath{S_t^\texttt{post}}}
\newcommand{\tautpo}{\ensuremath{\tau_t^\texttt{post}}}
\newcommand{\phitpo}{\ensuremath{\phi_t^\texttt{post}}}
\newcommand{\dtpo}{\ensuremath{\mathcal{D}_t^\texttt{post}}}
\newcommand{\vtst}{\ensuremath{V_t(S_t)}}
\newcommand{\vtstnext}{\ensuremath{V_{t+1}(S_{t+1})}}
\newcommand{\wset}{\ensuremath{\mathcal{W}}}
\newcommand{\vtpo}{\ensuremath{V_t^\texttt{post}}}
\newcommand{\vpoappr}{\ensuremath{\bar{V}_t^\texttt{post}}}
\newcommand{\vjtt}{\ensuremath{\bar{v}_j(\tau,t)}}
\newcommand{\vjbar}{\ensuremath{\bar{v}_j}}
\newcommand{\viphi}{\ensuremath{\bar{v}_i(\phi)}}
\newcommand{\vibar}{\ensuremath{\bar{v}_i}}
\newcommand{\mtj}{\ensuremath{\mathcal{M}_t(j)}}
\newcommand{\ytjm}{\ensuremath{y^t_{jm}}}
\newcommand{\ztk}{\ensuremath{z^t_k}}
\newcommand{\hjtau}{\ensuremath{\mathcal{H}_j(\tau)}}
\newcommand{\hjtaut}{\ensuremath{\mathcal{H}_j(\tau,t)}}
\newcommand{\vhold}{\ensuremath{v^\texttt{base}}}
\newcommand{\mhset}{\ensuremath{\mathcal{M}_h}}
\newcommand{\hranked}{\ensuremath{\mathcal{H}^\texttt{ranked}}}
\newcommand{\kranked}{\ensuremath{\mathcal{K}^\texttt{ranked}}}
\newcommand{\delot}{\ensuremath{\delta_{\texttt{ot}}}}
\newcommand{\costot}{\ensuremath{c_\texttt{ot}}}
\newcommand{\kiphi}{\ensuremath{\mathcal{K}_i(\phi)}}
\newcommand{\dasng}{\ensuremath{\mathcal{D}_t^\texttt{asng}}}
\newcommand{\noxb}{\ensuremath{N_\texttt{XB}}}
\newcommand{\vbarhld}{\ensuremath{\bar{v}^\texttt{hld}}}
\newcommand{\vbarhldot}{\ensuremath{\bar{v}^\texttt{hld-ot}}}
\newcommand{\lbartj}{\ensuremath{\bar{l}_{tj}}}
\newcommand{\lbarotti}{\ensuremath{\bar{l}^\texttt{ot}_{ti}}}

\newcommand{\nholdt}{\ensuremath{N^\texttt{hld}_t}}
\newcommand{\vhattj}{\ensuremath{\hat{v}_{tj}}}
\newcommand{\vhatotti}{\ensuremath{\hat{v}^\texttt{ot}_{ti}}}
\newcommand{\navg}{\ensuremath{n^\texttt{avg}}}
\newcommand{\wtj}{\ensuremath{w_{tj}}}
\newcommand{\wotti}{\ensuremath{w^\texttt{ot}_{ti}}}
\newcommand{\qtj}{\ensuremath{q_{tj}}}
\newcommand{\qotti}{\ensuremath{q^{\texttt{ot}}_{ti}}}
\newcommand{\nmax}{\ensuremath{N^\texttt{max}}}
\newcommand{\hjelig}{\ensuremath{\mathcal{H}_j}}
\newcommand{\hjaftert}{\ensuremath{\mathcal{H}_j^>(t)}}

\newcommand{\hopenxit}{\ensuremath{\mathcal{H}_{\xi t}^\texttt{open}}}

\newcommand{\vtilde}{\ensuremath{\tilde{v}}}
\newcommand{\ntop}{\ensuremath{n^\texttt{upper}}}
\newcommand{\rthr}{\ensuremath{r^2_\texttt{thr}}}

\newcommand{\setItemSep}{\setlength\itemsep}

\allowdisplaybreaks

\newcommand{\tcb}{\textcolor{black}}

\begin{document}\sloppy
\setlength{\parindent}{2em}

\begin{frontmatter}

\title{Approximate Dynamic Programming for Real-time Assignment of Extraboard Transit Operators}

\author[1]{Jilin Song}
\ead{jilin.song@mail.utoronto.ca}

\author[1]{Amer Shalaby}
\ead{amer.shalaby@utoronto.ca}

\author[2]{Merve Bodur\corref{cor1}}
\ead{merve.bodur@ed.ac.uk}

\cortext[cor1]{Corresponding author}

\address[1]{Department of Civil and Mineral Engineering, University of Toronto, Toronto, Canada}
\address[2]{School of Mathematics and Maxwell Institute for Mathematical Sciences, University of Edinburgh, Edinburgh, UK}

\begin{abstract}
\tcb{We investigate} real-time assignment decisions for extraboard transit operators, who are responsible for covering open work due to unexpected events such as driver absenteeism. Efficient usage of extraboard operators is critical as open work negatively affects service reliability. \tcb{We formulate} the problem as a Markov decision process to capture its stochastic and sequential nature. \tcb{We develop an approximate dynamic programming approach that selects assignments through} an integer program \tcb{using rewards from known open work and estimates of future rewards}. \tcb{We compute future-value estimates for individual operators offline using customized backward recursions}. \tcb{We incorporate a simulation-based, assignment-dependent adjustment for competition among operators for the same future work}. \tcb{We conduct} case studies \tcb{using published service schedules and constructed operator duties for one transit agency, and a larger synthetic instance}. \tcb{We examine} key performance metrics including uncovered open work and extraboard utilization rates for varying \tcb{work-source opening probabilities} and extraboard roster sizes. The approximate policy \tcb{achieves higher mean rewards than the tested myopic, practice-oriented, and scenario-based rolling-horizon benchmarks}. \tcb{We further examine the} inclusion of overtime drivers in the reserve operator roster \tcb{and passenger-delay-weighted coverage rewards}. Observations from these computational analyses provide \tcb{case-study insights into extraboard utilization, overtime use, and the prioritization of open work}.
\end{abstract}


\end{frontmatter}



\vspace{-8pt}

\section{Introduction}
In public transportation, service planning and scheduling are often done according to predetermined schedules known to transit users. During the daily operation of a public bus system, however, required operation tasks are frequently unassigned, becoming \textit{open work}. Unscheduled absence of regular-duty operators is the predominant source of open work \citep{diab2014extraboard}. In fact, urban bus operators are known to have much higher absenteeism rates than average blue-collar workers \citep{evans1994working}. Open work may also take the form of extra trips dispatched as a response to equipment breakdowns, weather, surge in passenger volume, etc. \citep{gupta2016reserve}.

Many public transit agencies have a roster of \textit{extraboard operators} at the agency's operating garages. These operators do not have fixed sets of trips to operate, unlike regular-duty operators. Instead, they are designated to cover open work. Extraboard operators are usually full-time employees who work 8-hour shifts daily, and they can make up a significant portion of the operator workforce (20\% on average at large systems) according to a survey by \citet{deannuntis2008transit}. It is common practice at many transit agencies to place newly trained operators on the extraboard roster \citep{boyle2009controlling}. Some extraboard operators cover known-in-advance open work such as planned leaves or vacations. The schedule for these extraboards can be made rather easily: they simply inherit the original work of the absent operator. 

The focus of this study is the assignment decisions for operators who cover unexpected open work in a real-time fashion. The operator workforce in question includes both extraboard operators and overtime regular-duty drivers. As labor compensation accounts for over half of all transit operating costs \citep{vuchic2017urban}, smart assignment decisions for available extraboard and overtime drivers \tcb{are} key to operational efficiency. Real-time assignment of extraboard and overtime operators is a challenging problem due to its stochasticity. Throughout the course of a workday, uncertain outcomes in the amount of open work unfold in a sequential manner. Despite the importance of these real-time assignment decisions, \tcb{practice reports describe reliance} on ad‑hoc rules and supervisor experience \tcb{\citep{boyle2009controlling,deannuntis2008transit}}. The academic literature also offers limited guidance on \tcb{real-time assignment policies that use work-source opening probabilities to anticipate future opportunities while accounting for both extraboard and overtime availability}. 

In this paper, we aim to address these research gaps. We model the assignment problem as a Markov decision process (MDP), designed to handle the problem's stochastic and sequential \tcb{nature}. \tcb{The MDP balances immediate assignment rewards against future opportunities at each decision-making epoch.} \tcb{We develop an approximate dynamic programming (ADP) policy that combines offline estimation of future rewards with online assignment optimization.} The proposed method is evaluated using \tcb{a schedule-based MiWay case study and a larger synthetic instance}. Its performance is compared against \tcb{myopic and practice-oriented rules, and a scenario-based rolling-horizon policy}.

The main contributions of this paper are as follows:
\begin{itemize}\setItemSep{-0.3em}
    \item We present, to the best of our knowledge, the first MDP framework \tcb{for} real-time extraboard operator assignment \tcb{that incorporates overtime drivers and the sequential revelation of open work}.
    \item We develop a tractable real-time assignment policy that produces high-quality decisions by \tcb{combining operator-specific backward value calculations with an assignment-dependent competition adjustment in an online integer program}.
    \item \tcb{Through schedule-based and synthetic case studies}, we \tcb{report higher mean rewards than the tested myopic, practice-oriented, and scenario-based rolling-horizon benchmarks}. \tcb{Matched component comparisons assess the effects of individual future values and the competition adjustment.}
    \item We offer \tcb{case-study} insights on the usage of overtime drivers, the effects of extraboard sizing on system performance\tcb{, and the prioritization of open work}.
\end{itemize}

The remainder of this paper is organized as follows. Section~\ref{sec:litReview} reviews relevant literature and highlights the key gaps motivating this work. This is followed by the problem definition and its mathematical formulation in Section~\ref{sec:problem}. We present the proposed solution methodology in Section~\ref{sec:methods}. Numerical results are then presented in Section~\ref{sec:results} through two case studies to evaluate the performance of the proposed method. The paper concludes in Section~\ref{sec:conclusion} with a summary of the main findings and suggestions for future research.

\section{Literature Review}
\label{sec:litReview}

\subsection{Planning and scheduling of extraboard operators}

Crew scheduling, or scheduling of transit operators, is a key task in the service planning process at public transit agencies. The Transit Cooperative Research Program (TCRP) Report 135 published by the Transportation Research Board \citep{boyle2009controlling} details operator scheduling methods often used in practice. In the academic literature, transit operator scheduling is a major topic of interest \citep{ibarra2015planning}. The problem is often formulated as a set covering problem which minimizes system costs while ensuring all required services are covered and satisfying constraints regarding work rules.

On the other hand, planning and scheduling of extraboard operators have received relatively little attention compared to the decision-making for regular operators. In the \tcb{practice reviewed by \citet{boyle2009controlling}}, decisions for extraboard operators \tcb{often drew on} transit supervisors' past experience. \tcb{The earlier survey by \citet{deannuntis2008transit} reported limited software support for these decisions despite the use of software for regular-operator schedules.} Extraboard operator scheduling has not been extensively studied in the literature.

\begin{figure}[!ht]
    \centering
    \includegraphics[width=0.805\textwidth]{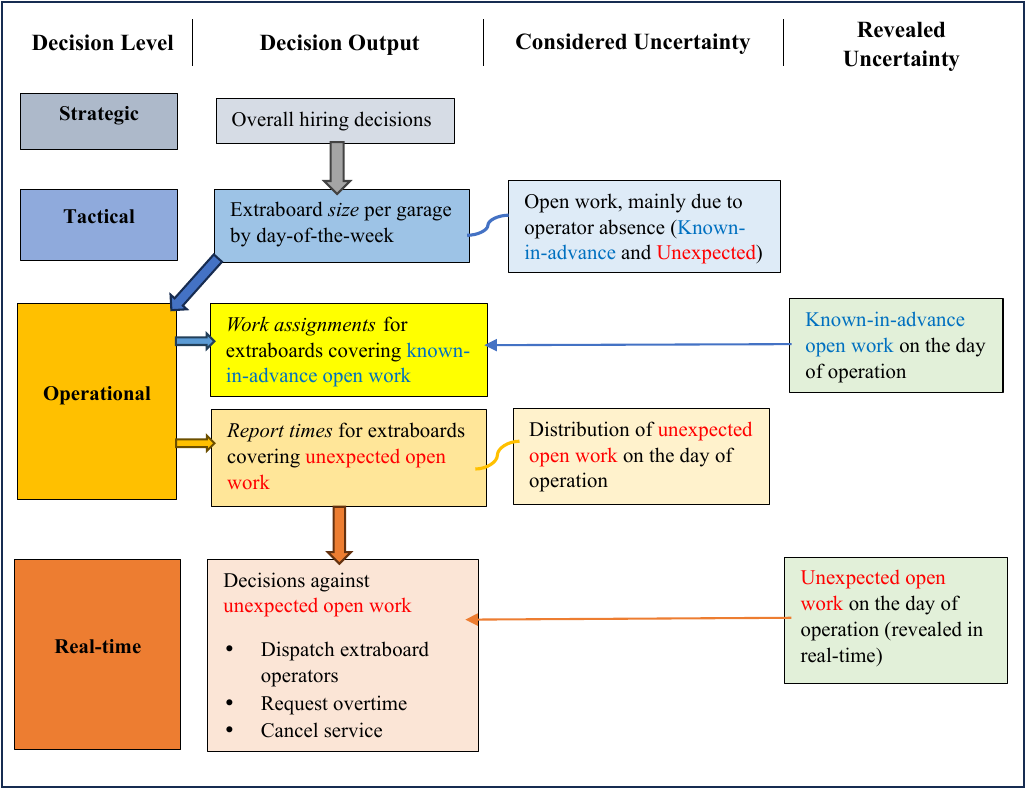}
    \caption{Levels of Extraboard Planning and Scheduling}
    \label{fig:flowchart}
\end{figure}

An early study on this topic by \citet{koutsopoulos1990scheduling} discussed three decision levels in extraboard planning and scheduling: strategic, tactical, and operational levels. Strategic-level decisions involve overall workforce sizing and hiring at the transit agency. Tactical-level decisions are made at the start of each planning period (typically every 3--4 months) regarding the size of the extraboard roster by day-of-week at each garage. Operational-level decisions are made daily. The decisions are the ``report times'' of extraboard operators for the next workday, the time at which they begin their work shift at the garage and await work assignments. Below the three aforementioned decision levels, transit supervisors also make daily \textit{real-time} decisions for usage of extraboard operators against open work throughout the workday. The interrelations among the strategic, tactical, operational, and real-time decisions were summarized by \citet{song2025extraboard}, as presented in Figure \ref{fig:flowchart}. This flowchart highlights two important aspects about extraboard decision-making: (1) outputs from upper decision levels serve as input to lower levels; (2) decisions have to be made under uncertainty that will be revealed after decisions are finalized.

In the state of practice, extraboard sizing at the tactical decision level is often made with simple, established methods. For example, many agencies use a pre-determined absence rate and multiply it by the number of regular-duty operators to determine the extraboard size \citep{deannuntis2008transit}. In the literature, various optimization-based methods have been proposed which attempt to tackle the stochasticity of the problem. \citet{koutsopoulos1990scheduling} developed \tcb{models for the daily number of extraboard operators and their report times}. \citet{ozbay2014understanding} proposed a model to prescribe the optimal number of extraboard driver-hours per day using probabilistic constraints.

Determination of report times at the operational level is another interesting decision-making problem involving uncertainty. The TCRP Report by \citet{boyle2009controlling} did not mention any common method used by transit supervisors. \citet{kaysi1990scheduling} proposed an algorithmic approach that gives extraboard report times one after the other while considering varying absenteeism rates of regular drivers at different times of the day. \citet{song2025extraboard} modeled the report time scheduling problem as a two-stage stochastic program. The two-stage structure allows for modeling \tcb{the assignment of} extraboard operators to open work tasks, which are encapsulated in the second-stage decisions. \tcb{Relatedly, \citet{li2014extraboard} jointly modeled extraboard shift-start times and work assignments under work-time and spread-time constraints, and developed heuristic approaches for that operational scheduling problem.}

\subsection{Real-time assignment decisions}

Real-time assignment decisions for extraboard operators are the focus of this study. Existing research on this topic is limited, as is the case for tactical and operational decision levels. \citet{gupta2016reserve} \tcb{studied sequential assignment of reserve drivers to jobs arriving without knowledge of future requests}. \tcb{They proposed a randomized online algorithm with a competitive-ratio guarantee relative to an offline solution, balancing acceptance of feasible jobs against retaining capacity for later requests.} Their proposed algorithm is \tcb{motivated by} fixed interval scheduling \citep{kovalyov2007fixed}, where required jobs with defined start and end times are allocated to one or more parallel machines. \tcb{Our setting additionally distinguishes extraboard and overtime availability, allows work-piece rewards beyond duration-based coverage, and uses a stochastic model of future source openings to evaluate current assignments.}

Aside from extraboard operators, transit agencies may also have regular drivers work overtime to cover open work. Even though a premium is often associated with overtime work (e.g., 1.5 times the regular pay), overtime drivers can often be more economical \tcb{than extraboard operators} due to their flexibility \citep{boyle2009controlling}. However, over-reliance on overtime has its risks. Service reliability can be compromised in the situation where few drivers are willing to work overtime when needed \citep{shiftan1994absence}. An abundance of overtime opportunities may also induce absenteeism for some drivers recognizing that overtime premium can easily make up for the lost regular pay \citep{perry1984extraboard}. At large transit agencies, where extraboard operators make up a significant portion of the driver workforce, overtime is used as a supplemental source to cover open work. \citet{deannuntis2012best} surveyed the workforce management practices at Miami-Dade Transit. The authors detailed the daily procedure for taking names of drivers willing to work overtime, and this serves as the reference for the problem definition discussed later in this paper. On the other hand, the academic literature scantly addresses scheduling and operational decisions for overtime operators. \citet{kaysi1990scheduling} mentioned that open work not covered by extraboards is either given to an overtime driver or canceled. To consider the effects of overtime drivers when deciding extraboard report times, the authors proposed to place higher weights on the morning hours in terms of service loss to reflect that there is usually less overtime availability in the morning. The study did not discuss how overtime drivers are assigned to open work in real-time operations.

\tcb{Table \ref{tab:compare_xb_studies} summarizes reviewed literature on various levels of extraboard planning and scheduling.}

\begin{table}[ht!]
\centering
\small
    \caption{\tcb{Comparison of studies in extraboard planning and scheduling}}
    \label{tab:compare_xb_studies}
    {
    \renewcommand{\arraystretch}{1.2}
    
    \begin{tabularx}{\textwidth}{L{2.1cm} L{2.2cm} L{2.2cm} Y Y}
    \toprule
    Study & Decision level & Operator type  & Treatment of uncertainty & Model / solution method \\
    \midrule
    \citet{deannuntis2008transit} & Tactical & Extraboard & Pre-defined absence rates of regular operators & Rule-of-thumb sizing based on absence rates found in industry practice \\
    \citet{koutsopoulos1990scheduling}  & Tactical / operational  & Extraboard  & Uncertain daily work requirements & Models for daily extraboard staffing and report times \\
    \citet{ozbay2014understanding}  & Tactical & Extraboard & Empirical distribution of daily open work & Prescription of extraboard driver-hours satisfying probabilistic constraints \\
    \citet{boyle2009controlling}  & Operational & Extraboard & N/A; industry practice review  & No specific method for determining report times identified; example schedules shown  \\
    \citet{kaysi1990scheduling}  & Operational & Extraboard, with overtime as consideration & Binomial distribution for absent drivers per time period & Specialized algorithm for determining extraboard report times one after another \\
    \citet{song2025extraboard}  & Operational  & Extraboard & Open work arising from absenteeism of individual regular operators & Two-stage stochastic program with customized solution algorithm \\
    \citet{li2014extraboard} & Operational & Extraboard & Open work from unplanned events represented as fixed jobs & Joint shift-start and work-assignment model with heuristic algorithms \\
    \citet{gupta2016reserve}  & Real-time       & Extraboard & Jobs arrive without knowledge of future requests  & Randomized online interval-scheduling algorithm with competitive guarantee \\
    \citet{deannuntis2012best}  & Operational / real-time discussion & Overtime & N/A; industry case study & Description of a daily procedure for managing drivers working overtime  \\
    This study & Real-time & Extraboard and overtime & Probabilistic openings of work sources over the day & MDP with individual-value approximation, competition adjustment, and online integer program \\
    \bottomrule
    \end{tabularx}
    }
\end{table}

\subsection{Sequential decision-making and the Markov decision process}

A common theme of the decision-making involving extraboards is uncertainty. However, a key difference exists between real-time assignment decisions and the upper level decisions. Both the extraboard sizing and report time problems are one-time decisions made at their respective frequencies. Real-time decisions, on the other hand, are made sequentially: assignments are made throughout the planning horizon while anticipating future uncertainties.

The Markov decision process (MDP) is a mathematical model that describes sequential decision-making problems under uncertainty. The MDP has been used to model many transportation-related problems. Example applications include dynamic fleet management \citep{godfrey2002adaptive}, dynamic crowd-shipping for online order delivery \citep{mousavi2024approximate}, \tcb{and} scheduling of autonomous buses~\citep{wu2026dynamic}. \tcb{Related stochastic decision models have also been used to derive} bus holding \tcb{policies that minimize passenger waiting time} \citep{berrebi2015real}. Ride sourcing, including ride hailing and ride sharing, is a popular topic that shares similarities with extraboard scheduling. In daily operations, the ride sourcing platforms are tasked with dynamically matching the supply side (e.g., drivers) with the demand side (e.g., passengers, delivery orders) \citep{wang2019ridesourcing}. The objective is usually solving for a policy that maps any given state to its optimal decisions. For example, \citet{yu2019markov} uses MDP to model the vacant taxi routing problem and finds an optimal routing policy that maximizes long-term profit over the working period. 

For simple MDP models, an optimal policy can be obtained with dynamic programming. In real-world contexts, MDP models may suffer from one or more ``curses of dimensionality'': having large state, outcome, and/or action spaces \citep{powell2007approximate}. As such, exact solutions are often difficult to obtain. To tackle this, the approximate dynamic programming (ADP) method has been used. \citet{yu2019integrated} \tcb{formulated an on-demand ride-pooling} problem as a \tcb{multistage stochastic program} and \tcb{combined spatial and temporal decomposition with} ADP. \citet{jiang2023adaptive} used ADP to solve a large-scale multi-driver order dispatching problem. \tcb{Another approach for tackling large-scale MDPs is to use rolling-horizon look ahead policies. For example, \citet{alnaggar2025heatmap} employed a stochastic look ahead policy for driver matching in an MDP representation of crowdsourced delivery with probabilistic driver repositioning.}

Another popular method for tackling sequential decision-making problems is reinforcement learning (RL). The environment for \tcb{an} RL model often has an underlying MDP structure. \tcb{RL includes both model-free and model-based approaches; it is therefore not distinguished from ADP solely by whether a transition model is known} \citep{li2023reinforcement}. In the transportation context, RL has been applied to a variety of problems, such as traffic signal control \citep{yau2017survey}, dynamic holding control to prevent bus bunching \citep{wang2020dynamic}, and the aforementioned dynamic ride sourcing problem \citet{al2019deeppool}.

\tcb{For real-time extraboard assignment, work-source outcomes and feasible assignments can be specified explicitly. We use this structure to estimate future values offline and optimize feasible current assignments online.} In the next sections, the proposed MDP formulation and \tcb{approximate policy} are discussed.

\section{Problem Definition and Formulation}
\label{sec:problem}

This section details the Markov decision process (MDP) formulation for the real-time extraboard work assignment problem. \tcb{We} first \tcb{introduce} the \tcb{model} parameters and notation. Then, we discuss each of the important elements in the proposed MDP model: state variables, decision variables, reward function, exogenous information, transition function, and policies.

In the proposed real-time extraboard work assignment problem, transit supervisors at a bus garage collaborate and effectively act as a single decision maker. For simplicity, we refer to this collective decision maker as the ``bus garage'' or ``garage''. The planning horizon is defined as one day of operation for the bus garage. The day is discretized into $T$ periods of equal length. We denote the set of periods in the planning horizon by $\mathcal{T}:=\{0,1,2,\hdots,T-1\}$. Durations of operator shifts and work pieces are measured in terms of number of periods. Assignment decisions are made at the beginning of each period\tcb{, after newly revealed open work is observed}. At the \tcb{beginning} of \tcb{the next} period, the state \tcb{reflects} these decisions \tcb{and any new open work revealed since}.

The bus garage is tasked with making real-time assignment decisions for \tcb{its} reserve operator workforce to cover unexpected open work. This can be interpreted as matching the ``supply'' of reserve drivers (i.e., extraboard and overtime operators) with its ``demand'' (i.e., open work). The next paragraphs discuss the supply and demand sides in detail.

The bus garage has a roster of $J$ extraboard operators (referred to as ``XBs'' hereafter) working on the day. The set of these XBs is denoted by $\mathcal{J}:=\{0,1,2,\hdots,J-1\}$. From upper-level decisions, the report time (i.e., the first period) and the last period of each XB's shift are known. For a particular XB $j\in \mathcal{J}$, these time periods are denoted by $s_j, e_j\in\mathcal{T}$, respectively. From this information, we can obtain the set of XBs \tcb{on duty} at every period in the planning horizon: $\mathcal{J}_t:=\{j\in\mathcal{J}:s_j\leq t\leq e_j\}, \forall t\in \mathcal{T}$.

In addition to XBs, the bus garage may also use drivers working overtime to cover open work. In our MDP, \tcb{the overtime-driver setup is motivated by the volunteer lists and minimum-pay guarantee described for} Miami-Dade Transit \tcb{\citep{deannuntis2012best}}. \tcb{We assume that the list of drivers willing to work overtime} (referred to as ``OTs'' hereafter) \tcb{and their report times are known} to the garage at the start of the planning horizon. We denote the set of OTs by $\mathcal{I}:=\{0,1,2,\hdots,I-1\}$. The report time of OT $i\in \mathcal{I}$ is denoted by $s_i\in \mathcal{T}$. OTs have a minimum pay duration, the length of which is denoted by $\delot$. For each period in the planning horizon, we define the set of OTs \tcb{whose assignment window includes period $t$}: $\mathcal{I}_t:=\{i\in\mathcal{I}:s_i\leq t\leq s_i+\delot\},\forall t\in \mathcal{T}$. The detailed decision-making timeline for OT $i\in \mathcal{I}$ is as follows: starting at $s_i$, this OT is said to be on ``standby'' awaiting assignment to open work. If the OT is assigned to a piece of open work, they are paid for the duration from $s_i$ to the end of that work assignment. If the OT \tcb{remains unassigned after the assignment decision at $s_i+\delot$}, they are sent home and paid for that minimum pay duration. In the event that the assigned work finishes before the end of the OT's minimum pay duration, they still receive that minimum pay. For example, consider an OT who reports at 5 PM with a minimum pay duration of one hour. They may be assigned work at any \tcb{decision epoch from} 5 to 6 PM, inclusive. If this OT gets assigned a work piece that starts at 5:30 and ends at 7, they get paid for two hours (from 5 to 7). Otherwise, after being on standby for the whole hour and never assigned any work, the OT is sent home and get paid for one hour (from 5 to 6). The per-period overtime \tcb{cost used in the assignment objective} is denoted by $\costot$.

We now discuss the demand side of the problem: open work to be assigned to XBs and OTs. As briefly mentioned in previous paragraphs, the concept of ``work piece'' is formally defined. A work piece is defined as \tcb{an indivisible} unit of work that must be carried out by a single bus operator. An example work piece can be a round trip on a particular bus route. This is the case for bus routes where the only feasible location for operator changeover (also known as ``relief location'') is at one of the terminals. The set of all required work pieces during the day is denoted by $\mathcal{K}:=\{0,1,2,\hdots,K-1\}$. Each work piece $k\in \mathcal{K}$ has a defined starting period ($a_k\in \mathcal{T}$) and duration ($\delta_k$); hence, its last period is given by $b_k=a_k+\delta_k-1$. The duration of the shortest possible work piece during the day is denoted by $\delmin=\min\{\delta_k\}_{k\in\mathcal{K}}$. An XB is eligible for any work piece that falls entirely within the XB's shift. The set of eligible work pieces for XB $j\in\mathcal{J}$ is denoted by $\mathcal{K}_j=\{k\in\mathcal{K}:(a_k\geq s_j) \ \texttt{and} \ (b_k\leq e_j) \}$. An XB may be assigned to an eligible work piece if they are \tcb{available for assignment at the current decision epoch}. Therefore, an XB may be assigned to multiple non-overlapping work pieces throughout their shift. On the other hand, OTs are eligible for any work piece with a start time within their \tcb{assignment window}. The eligible work piece set for OT $i\in\mathcal{I}$ is denoted by $\mathcal{K}_i=\{k\in\mathcal{K}:s_i\leq a_k\leq s_i+\delot\}$. An OT may be assigned to an eligible work piece if they are currently on standby (i.e., not assigned to work and \tcb{still within their assignment window}). To keep the OT shifts simple, we specify that each OT can only be assigned to at most one open work piece. An OT is sent home after finishing their assigned work piece. A reward ($c_k$) is associated with assigning piece $k\in\mathcal{K}$ to an XB or OT. This reward may also be understood as the cost saving from \tcb{avoiding} leaving the piece uncovered.

Over the planning horizon, outcomes of work pieces unfold, in which they are revealed to be either open or not open. Only \textit{open} work pieces may be covered by XBs and OTs. To characterize how work pieces become open, we define a set of ``work sources'' $\mathcal{H}:=\{0,1,2,\hdots,H-1\}$. Each work source $h\in \mathcal{H}$ consists of one or more work pieces in $\mathcal{K}$, while each work piece belongs to only one parent work source, denoted by $h_k$. The first work piece of the work source is denoted by $\kfirst$. An example work source is a scheduled run of a regular-duty operator. The operator run consists of multiple consecutive round trips, each being a work piece that belongs to the work source. The ``openness'' outcome of a work piece follows that of its parent source. All trips of an absent driver's run become open work pieces. Historical data may be used to estimate the probability that work sources become \tcb{open}. The start time of a work source is defined as the start time of its first component work piece: $a_h=a_{\kfirst}$.

\subsection{State variables} 
The state of the system at period $t\in \mathcal{T}$, before decision-making, is defined by the states of XBs, OTs, and revealed open pieces, denoted by $S_t = (\tau_t, \phi_t, \mathcal{D}_t)$. 

Vector $\tau_t = (\tfrj)_{j\in \mathcal{J}}$ captures the states of all XBs at period $t$. $\tfrj$ describes the next period when XB $j\in \mathcal{J}$ is available. Note that $\tau_0(j)$ is the report time for XB $j$. The set of working XBs currently available at $t$ is denoted by $\jtn:=\{j\in\mathcal{J}_t:\tfrj=t\}$. \tcb{An} XB \tcb{with an outstanding work assignment} at $t$ has their $\tfrj$ value set to the period immediately after the final period in their assigned work. XBs who already got off work will not appear in the set $\mathcal{J}_t$ (and thus $\jtn$) again as time progresses.

Vector $\phi_t = (\phiti)_{i\in \mathcal{I}}$ \tcb{tracks} each OT\tcb{'s standby status} as of $t\in \mathcal{T}$. For an OT $i$ who has started their overtime shift, a $\phiti$ value between 0 and $\delot$ \tcb{counts their standby periods and} indicates that they are on standby and available for assignment at $t\in\mathcal{T}$. We denote the set of these OTs by $\itn:=\{i\in\mathcal{I}_t: 0\leq\phiti\leq \delot\}$. OTs who are \tcb{no longer available for assignment (either because they were assigned work or because they remained unassigned after their final standby decision)} have their $\phi$ values fixed at $\delot+1$ for the rest of the planning horizon so that these OTs do not appear in $\itn$ again.

$\mathcal{D}_t\subseteq \mathcal{K}$ is the set of all revealed open work pieces known at $t\in\mathcal{T}$ \tcb{that have not yet been assigned and have a} start time no earlier than $t$, such that these pieces may be assigned to XBs or OTs. The set of open pieces that start exactly at $t$, a subset of $\mathcal{D}_t$, is denoted by $\dtn:=\{k\in\mathcal{D}_t:a_k=t\}$.

\subsection{Decision variables}\label{subsec:decision_variables}
At state $S_t$, assignment decisions are made to match currently available XBs ($\jtn$) and OTs ($\itn$) to open work pieces. 
For XBs, we define binary decision variables $\xtjk$ for all eligible $(j,k)$ pairs:
\[
\xtjk = 
\begin{cases} 
    1 & \text{if XB $j\in\jtn$ is assigned to piece $k\in\tcb{\mathcal{D}_t}\cap\mathcal{K}_j$} \\
    0 & \text{otherwise} \\
\end{cases}
\]
\tcb{For} OTs, binary decision variables $\utik$ are defined for all eligible $(i,k)$ pairs:
\[
\utik = 
\begin{cases} 
    1 & \text{if OT $i\in\itn$ is assigned to piece $k\in\dtn\cap\mathcal{K}_i$} \\
    0 & \text{otherwise} \\
\end{cases}
\]
\tcb{For an XB, a work sequence consists of one or more non-overlapping eligible open pieces. At $t$, an XB can be assigned only to a sequence whose first piece starts at $t$; all pieces in that sequence, including those starting later, are \emph{binding assignments}. The XB remains unavailable for other assignments through the last assigned piece, including any gaps between selected pieces.}
If a currently available XB is not assigned to any open work piece, they will hold for one period. \tcb{An OT not assigned a work piece starting at $t$ remains on standby for one period, unless $t=s_i+\delot$, when the OT is released after the assignment decision.}

Let $\xvec$ \tcb{and} $\decspace$ denote the vector of all decision variables and the decision space, respectively. $\xvec\in\decspace$ is defined by the following constraints:
\begin{align}
    &\sum_{k\in\dtn\cap \mathcal{K}_j}\xtjk\leq 1 &&\forall j\in\jtn \label{eqn:cons11}\\
    &\sum_{k\in\dtn\cap \mathcal{K}_i}\utik\leq 1 &&\forall i\in\itn \label{eqn:cons12} \\
    &\sum_{j\in\jtn:k\in\mathcal{K}_j}\xtjk + \sum_{i\in\itn:k\in\mathcal{K}_i}\utik\leq 1 &&\forall k\in\dtn \label{eqn:cons13} \\
    &\tcb{\sum_{j\in\jtn:k\in\mathcal{K}_j}\xtjk\leq 1} &&\tcb{\forall k\in\mathcal{D}_t\setminus\dtn} \label{eqn:cons14}\\
    &\tcb{\xtjk\leq\sum_{k'\in\dtn\cap\mathcal{K}_j}x^t_{jk'}} &&\tcb{\forall j\in\jtn,\ k\in(\mathcal{D}_t\setminus\dtn)\cap\mathcal{K}_j} \label{eqn:cons15}\\
    &\tcb{\sum_{\substack{k\in\mathcal{D}_t\cap\mathcal{K}_j:\\a_k\leq\ell\leq b_k}}\xtjk\leq 1} &&\tcb{\forall j\in\jtn,\ \ell=t+1,\ldots,e_j} \label{eqn:cons16}\\
    &\xtjk\in \{0,1\} &&\forall j\in\jtn, k\in\tcb{\mathcal{D}_t} \cap \mathcal{K}_j \label{eqn:cons17} \\
    &\utik\in \{0,1\} &&\forall i\in\itn, k\in\dtn \cap \mathcal{K}_i \label{eqn:cons18}
\end{align}
Constraints \eqref{eqn:cons11} and \eqref{eqn:cons12} specify that an XB or an OT can be assigned to at most one eligible open work piece \tcb{starting at $t$}, respectively. Constraint \eqref{eqn:cons13} specifies that a work piece \tcb{starting at $t$} can be covered by at most one eligible driver, whether that is an XB or OT. \tcb{Constraint \eqref{eqn:cons14} prevents multiple XBs from being assigned to the same later-starting piece. Constraint \eqref{eqn:cons15} requires a piece starting at $t$ in any XB sequence that includes later pieces. Constraint \eqref{eqn:cons16} prevents the pieces assigned to an XB from overlapping.} Constraints \tcb{\eqref{eqn:cons17}} and \tcb{\eqref{eqn:cons18}} define the variables as binary.

\subsection{Reward function}
When an open work piece $k\in\mathcal{K}$ is assigned to an XB or OT, a reward $c_k$ is gained. \tcb{For an XB sequence beginning at $t$, rewards for all its pieces, including those starting later, are credited when the sequence is committed.} If an OT is assigned a work piece that extends beyond their minimum pay duration, they are paid for the additional work hours at $\costot$ per period. Note that regular pay of XBs and minimum pay of OTs are not included in the reward function. This is because these are committed expenditures that must be paid out regardless of decisions. The \tcb{net} reward from making decisions $\xvec$ at $t\in\mathcal{T}$ is given by:
\begin{align}
    g_t(\xvec)= & \sum_{j\in \jtn}\sum_{k\in \tcb{\mathcal{D}_t}\cap\mathcal{K}_j} c_k \xtjk + \sum_{i\in \itn}\sum_{k\in \dtn\cap\mathcal{K}_i} c_k \utik 
 -\costot \sum_{i\in \itn}\sum_{k\in \dtn\cap\mathcal{K}_i} \bigg[ a_k+\delta_k-(s_i+\delot)\bigg]^+ \utik \label{eqn:immediate_reward}
\end{align}
where the notation $[\cdot]^+$ denotes $\max(0,\cdot)$. \tcb{The positive-part term is zero if the assigned piece ends within the OT's minimum-pay interval; otherwise, it counts only the periods beyond that interval.}

\subsection{Exogenous information}
\label{sec:exog}
We use $\wset_t$ to denote the newly revealed set of open work pieces at $t\in\mathcal{T}$. $\wset_t$ can alternatively be interpreted as the aggregated set of open pieces that are made known to the garage between $t-1$ and $t$ (for $t\geq1$). $\wset_0$ denotes the set of open work pieces that are known at the very start of the day of operation. 

The openness outcome of a work piece follows that of the source it belongs to. The outcome of source $h\in\mathcal{H}$ is assumed to follow a Bernoulli trial with known probability $p_h$, which is also referred to as the ``openness probability''. \tcb{We assume that the outcomes of distinct work sources are independent.} \tcb{When source $h$ is open, all its component pieces, including those starting later, are revealed together and added to the system at $a_h$.}

\subsection{State transition functions} \label{subsec:transition}
The system starts from an initial state $S_0=(\tau_0,\phi_0,\mathcal{D}_0)$, which captures the information known at the start of the planning horizon. $\tau_0=(\tau_0(j))_{j\in\mathcal{J}}$ is the vector of report times for the XB roster. No OT has spent time on standby in the initial state: $\phi_0(i)=0, \forall i\in\mathcal{I}$. The initial set of open work pieces $\mathcal{D}_0$ is the same as $\wset_0$. The transition from state $S_t$ to $S_{t+1}$ depends on both decisions $\xvec$ and exogenous information $\wset_{t+1}$. We introduce the post-decision state, $\stpo=(\tautpo,\phitpo,\dtpo)$, to represent the state reached immediately after decision-making and before the arrival of new information. We divide the transition into two parts: (i) from $S_t$ to $\stpo$; (ii) from $\stpo$ to $S_{t+1}$. We discuss transitions (i) and (ii) for $\tau_t$, $\phi_t$, and $\mathcal{D}_t$ separately.

\textit{Transition (i) of $\tau_t$}. XBs that receive a work assignment have their $\tau_t$ value updated to the period immediately after the \tcb{last piece of the assigned sequence}. Available XBs without assignment will hold and be available again in the next period. XBs not currently available do not have their $\tau_t$ values updated. This transition function is presented in mathematical form as follows: 
\begin{equation}
    \tautpo(j)=
    \begin{cases}
        \displaystyle\mathop{\tcb{\max}}\limits_{k\in\tcb{\mathcal{D}_t}\cap\mathcal{K}_j\tcb{:\,\xtjk=1}}(b_k+1) &\forall j\in\jtn:\displaystyle\sum_{k\in\tcb{\mathcal{D}_t}\cap\mathcal{K}_j}\xtjk\tcb{\geq}1\\
        \tfrj+1 &\forall j\in\jtn: \displaystyle\sum_{k\in\tcb{\mathcal{D}_t}\cap\mathcal{K}_j}\xtjk=0 \\
        \tfrj &\forall j\in\mathcal{J}\backslash\jtn
    \end{cases}
\end{equation}
\textit{Transition (ii) of $\tau_t$}. Because there is no exogenous information associated with XBs, the state of XBs at $t+1$ will remain the same as their post-decision state from the last period: $\tau_{t+1}=\tautpo$.

\textit{Transition (i) of $\phi_t$}. Currently available OTs either get assigned a work piece or \tcb{remain unassigned at the current decision epoch}. OTs being assigned have their $\phi_t$ updated to $\delot+1$ so that they do not become available again. OTs who are \tcb{unassigned} have $\phi_t$ increase by one. \tcb{At their final standby decision, this update gives a post-decision value of $\delot+1$ and releases them.} $\phi_t$ is not updated for unavailable OTs. Presented in mathematical form:
\begin{equation}
    \phitpo(i)=
    \begin{cases}
        \delot+1 &\forall i\in\itn: \displaystyle\sum_{k\in\dtn \cap \mathcal{K}_i}\utik = 1\\
        \phiti+1 &\forall i\in\itn: \displaystyle\sum_{k\in\dtn \cap \mathcal{K}_i}\utik = 0 \\
        \phiti &\forall i\in\mathcal{I}\backslash\itn
    \end{cases}
\end{equation}
\textit{Transition (ii) of $\phi_t$}. Similar to XBs, there is no exogenous information associated with OTs. Therefore, $\phi_{t+1}=\phitpo$.

\textit{Transition (i) of $\mathcal{D}_t$}. The post-decision state of open pieces is obtained by removing all pieces with start time $t$ \tcb{and any later-starting pieces committed to XBs}: $\dtpo=\mathcal{D}_t\backslash\tcb{(\dtn\cup\dasng)}$, \tcb{where $\dasng=\{k\in\mathcal{D}_t:\sum_{j\in\jtn:k\in\mathcal{K}_j}\xtjk=1\}$}. Recall that \tcb{a piece still awaiting assignment when its start time arrives must then be assigned or dropped}.

\textit{Transition (ii) of $\mathcal{D}_t$}. Newly revealed open work pieces, $\mathcal{W}_{t+1}$, are added to the post-decision set of open pieces: $\mathcal{D}_{t+1}=\dtpo\tcb{\cup}\wset_{t+1}$.

\subsection{Optimal and approximate policies}
At state $S_t$, \tcb{an optimal policy makes} driver assignment decisions to maximize the immediate reward plus the expected reward across the remaining planning horizon. Optimal decisions can be obtained by solving the Bellman equation:
\begin{equation} \label{eqn:bell1}
    \vtst=\max_{\xvec\in\decspace} g_t(\xvec) + \mathbb{E}_{\wset_{t+1}} [\vtstnext|S_t,\xvec]
\end{equation}
where state $S_{t+1}$ is transitioned from $S_t$ with decisions $\xvec$ and exogenous information $\wset_{t+1}$. $\vtst$ is the \textit{value function} \tcb{evaluated at the current state}. \tcb{We set $V_T(S_T)=0$ for every terminal state $S_T$.} An optimal policy could theoretically be found by obtaining the exact value function for every possible state using backward dynamic programming. However, due to the system's very large state space, that approach is computationally intractable. In the next section, we detail an approximate policy that is tractable to solve. 

We propose a tractable solution method that finds an approximation for the value function. The approach results in an approximate policy that produces high-quality XB and OT assignment decisions. Its details are presented in the next section.

\section{Solution Methodology}
\label{sec:methods}

\tcb{This section presents a customized approximate dynamic programming (ADP) approach for the real-time extraboard assignment problem formulated as an MDP. Optimizing assignments solely against already revealed work would disregard the value of retaining operators for work revealed later. We therefore approximate the post-decision value function by summing estimated future values for individual XBs and OTs, with adjustments for competition among operators for the same future work. Separate backward dynamic programming algorithms calculate the individual values for a specified range of report times and operator states, while a simulation-based algorithm estimates the adjustment parameters. These values and parameters are computed off-line and used in an integer program for real-time assignment decisions, reducing the computational effort required during daily operations. The individual values can also be reused across daily combinations of report times within the precomputed range.}

\tcb{The following subsections first define the post-decision value function and its individual-operator approximation. We then introduce an integer programming formulation of the approximate policy and develop the backward dynamic programs used to estimate individual operator values. Finally, we present adjustments intended to mitigate overestimation from summing these values and describe the simulation-based estimation of the adjustment parameters.}

\subsection{Post-decision value function}
We define the post-decision value function, denoted by
\begin{equation}
    \vtpo(\stpo)=\mathbb{E}_{\wset_{t+1}}[\vtstnext|\stpo] = \mathbb{E}_{\wset_{t+1}} [\vtstnext|S_t,\xvec]
\end{equation} 
\tcb{For a decision $\xvec$ made in state $S_t$, $\stpo$ is the resulting post-decision state.} Accordingly, the Bellman equation can be rewritten as follows:
\begin{equation}
    \vtst=\max_{\xvec\in\decspace} \tcb{\bigl\{} g_t(\xvec) + \vtpo(\stpo) \tcb{\bigr\}}
\end{equation}
The post-decision state space \tcb{remains prohibitively large for exact computation}. We \tcb{therefore} approximate the post-decision value function \tcb{as explained next}.

\subsection{Value function approximation}\label{subsec:value_approximation}
The post-decision state $\stpo=(\tautpo,\phitpo,\dtpo)$ describes the states of both the driver supply (in XBs and OTs) and demand (in open work pieces). \tcb{We approximate its value as a function of the driver supply vectors:}
\begin{equation}
    \vtpo(\stpo) \approx \vpoappr(\tautpo,\phitpo)
\end{equation}
\tcb{The pieces in $\dtpo$ are known at $t$ and remain available for later decisions. The integer program described in the next subsection evaluates candidate plans over the known pieces in $\mathcal{D}_t$, while the individual values introduced below estimate opportunities from work that may be revealed after $t$.}
\tcb{With this separation in place, we define a tractable supply-based approximation as a sum of individual XB and OT values:}
\begin{equation} \label{eqn:post_simple}
    \vpoappr(\tautpo,\phitpo) = \sum_{j\in\tcb{\mathcal{J}_t}}\vjbar(\tcb{\tautpo(j)},t) + \sum_{i\in\tcb{\mathcal{I}_t}}\vibar(\tcb{\phitpo(i)})
\end{equation}
\tcb{The sums include only operators within their shift or assignment window at $t$, matching the domains of the individual value functions defined below.}
The individual value for XB $j\in\mathcal{J}$ is denoted by $\vjtt$. Parameter $\tau$ \tcb{denotes a possible next-availability period for} XB $j$ \tcb{after the decision; in~\eqref{eqn:post_simple}, $\tau$ equals} $\tautpo(j)$. $\vjtt$ can be interpreted as \tcb{an estimate of} the expected \tcb{future} reward from covering open pieces that \textit{will be revealed} after $t$ by the XB $j$ working in their remaining shift interval $[\tau,e_j]$. Note that the parameter $t$, time of decision-making, is a necessary input to the value function. This is because the openness outcomes of more work sources are revealed as time progresses. For example, consider two scenarios: (1) at 6 AM, an XB is assigned to a 2-hour long open work piece; (2) at 7 AM, an XB is assigned to a 1-hour long open work piece. \tcb{The XB's next-availability time is 8 AM in both scenarios. Nevertheless, the individual future-value estimates need not coincide, because the set of eligible work sources whose outcomes remain unknown may differ between 6 and 7 AM.}

The individual value function for OT $i\in\mathcal{I}$ is denoted by $\viphi$, where $\phi$ \tcb{denotes a possible post-decision OT state; in~\eqref{eqn:post_simple}, $\phi$ equals} $\phitpo(i)$. Similar to that of an XB, this value can be understood as \tcb{an estimate of} the expected \tcb{net} reward \tcb{from a possible assignment starting during the OT's remaining standby window} $[s_i+\phi,s_i+\delot]$. \tcb{For an OT still on standby, $s_i+\phi$ is the next decision epoch, so a separate $t$ argument is unnecessary. Under the single-piece assignment rule specified in Section~\ref{sec:problem}, an OT assigned to a work piece or released has no subsequent assignment opportunities. The reward and any overtime cost of an assigned piece enter the immediate reward $g_t$ in~\eqref{eqn:immediate_reward}, while the OT's remaining future value is zero.}

\tcb{A practical advantage of the proposed decomposition is that individual value estimates can be reused across daily rosters with different combinations of XB/OT report times. By contrast, a joint value function approximation may need to be recalculated for each new combination of report times, which can be computationally expensive. When the work-piece and work-source data, the OT cost coefficient, and the minimum-pay duration remain fixed, the garage can instead retrieve and combine individual values precomputed for the day’s report times and shifts, without recalculating those values. This is useful when report-time combinations vary from day to day within the precomputed set.}

\subsection{Integer programming model for the approximate policy}
\tcb{For the online optimization, we incorporate the individual future-value estimates into an integer program (IP) for XB and OT assignment decisions. The IP selects a plan over known open work, which may include XB sequences and OT pieces starting after $t$. These future-starting selections are provisional: they inform current decisions without committing the operator and may be revised later.}

\tcb{To represent possible XB assignments in the plan,} we \tcb{construct a set of candidate} work sequences, $\mathcal{M}_t$, \tcb{each consisting of one or more non-overlapping} open work pieces in $\mathcal{D}_t$. We further define the set of eligible work sequences for each XB $j\in\jtn$, denoted by $\mtj$, which is a subset of $\mathcal{M}_t$. Each sequence in this set only consists of eligible work pieces for XB $j$. That is: $\mtj=\{m\in\mathcal{M}_t:k\in\mathcal{K}_j,\forall k\in m\}$. Let $(a_m, b_m)$ denote the first and last periods covered by sequence $m$, defined as \tcb{the earliest start and latest end among its component pieces}, respectively.

\tcb{We represent XB selections using sequence-indexed variables $\ytjm$ in place of the piece-indexed variables $\xtjk$ from Section~\ref{subsec:decision_variables}. We also introduce $\ztk$ to indicate planned coverage by either an XB or an OT:}
\begin{align*}
    \ytjm &= 
\begin{cases} 
    1 & \text{if XB $j\in\jtn$ \tcb{is selected to cover} sequence $m\in\mtj$ \tcb{in the plan}} \\
    0 & \text{otherwise} \\
\end{cases} \\
\ztk &= 
\begin{cases} 
    1 & \text{if open piece $k\in\mathcal{D}_t$ \tcb{is selected for coverage in the plan} (by either an XB or OT)} \\
    0 & \text{otherwise} \\
\end{cases}
\end{align*}
\tcb{For OTs, we continue to use piece-indexed variables, defined here over all eligible known pieces in the plan:}
\[
\utik = 
\begin{cases} 
    1 & \text{if OT $i\in\itn$ \tcb{is selected to cover} piece $k\in\tcb{\mathcal{D}_t}\cap\mathcal{K}_i$ \tcb{in the plan}} \\
    0 & \text{otherwise} \\
\end{cases}
\]
\tcb{Using the additive approximation in \eqref{eqn:post_simple}, we obtain a provisional plan} at a given state $S_t$ \tcb{by solving the following IP}:
\begin{align}
    \max \ &\sum_{k\in\mathcal{D}_t} c_k\ztk -\costot \sum_{i\in \itn}\sum_{k\in \mathcal{D}_t\cap\mathcal{K}_i} \bigg[ a_k+\delta_k-(s_i+\delot)\bigg]^+ \utik\notag\\
    + &\sum_{j\in\jtn} \bigg( \sum_{m\in\mtj}\ytjm\vjbar(\tau=b_m+1,t) + 
    \Big(1-\sum_{m\in\mtj}\ytjm\Big)\vjbar(\tau=t+1,t)\bigg)\notag\\ 
    + &\sum_{i\in\itn}\Big(1-\sum_{k\in\mathcal{D}_t\cap\mathcal{K}_i}\utik\Big)\vibar(\phi=\phiti+1) \label{eqn:obj} \\
    \textrm{s.t.} \ &\sum_{m\in\mtj}\ytjm\leq 1 &&\-\hspace{-3.2cm}\forall j\in\jtn \label{eqn:cons21} \\
    &\sum_{k\in\mathcal{D}_t\cap \mathcal{K}_i}\utik\leq 1 &&\-\hspace{-3.2cm}\forall i\in\itn \label{eqn:cons22} \\
    &\sum_{j\in\jtn}\sum_{\substack{m\in\mtj:\\k\in m}} \ytjm + \sum_{\substack{i\in\itn:\\k\in\mathcal{K}_i}} \utik = \ztk &&\-\hspace{-3.2cm}\forall k\in\mathcal{D}_t \label{eqn:cons23} \\
    &\ytjm\in \{0,1\} &&\-\hspace{-3.2cm}\forall j\in\jtn, m\in\mtj \label{eqn:cons24} \\
    &\ztk\in \{0,1\} &&\-\hspace{-3.2cm}\forall k\in\mathcal{D}_t \label{eqn:cons25} \\
    &\utik\in \{0,1\} &&\-\hspace{-3.2cm}\forall i\in\itn, k\in\mathcal{D}_t \cap \mathcal{K}_i \label{eqn:cons26}
\end{align}
The objective function \eqref{eqn:obj} \tcb{scores the provisional plan using} three terms: \tcb{planned reward net of OT costs}, \tcb{estimated} future values of available XBs, and \tcb{estimated} future values of available OTs. The first line of the objective \tcb{includes the reward for each selected piece, including pieces starting after $t$, and the additional cost of pieces selected for OTs. Realized reward is calculated from the assignments committed at $t$}.
The second line denotes the sum of individual future values for available XBs ($\forall j\in\jtn$): the first term inside the big parentheses is the \tcb{value estimate applied when an XB is selected in the plan for} work sequence $m\in\mtj$, \tcb{using $\tau=b_m+1$ as its planned next availability}. The second term represents the value if the XB is \tcb{not selected in the plan}. The third line of the objective function is the sum of future values for available OTs ($\forall i\in\itn$). \tcb{This plan-level term is included only when an OT has no selected piece; the physical standby update depends on the executed decision described below.}
Constraint \eqref{eqn:cons21} ensures that each XB is \tcb{selected for} at most one eligible work sequence. Constraint \eqref{eqn:cons22} ensures that each OT is \tcb{selected for} at most one eligible work piece. \eqref{eqn:cons23} \tcb{links $z_{tk}$ to the planned coverage of} each open work piece in $\mathcal{D}_t$ \tcb{by an XB or OT and, because $z_{tk}$ is binary, prevents more than one operator from being selected for the same piece}. Constraints \eqref{eqn:cons24}, \eqref{eqn:cons25}, \eqref{eqn:cons26} define all decision variables as binary.

From the garage's perspective, \tcb{assignment decisions are made} at the decision frequency dictated by the length of time periods in $\mathcal{T}$ (e.g., every 15 minutes). \tcb{The IP is solved when at least one operator is available and an open piece starts at $t$. After a solve, each selected XB sequence with $a_m=t$ is committed in full, including its later pieces, and each selected OT piece with $a_k=t$ is assigned. All other available XBs hold for one period, and all other available OTs remain on standby for one period or are released after their final standby decision; this includes operators selected only for future-starting work. When the IP is not solved, we take all selection variables as zero and apply the same hold, standby, or release rules. In either case, the post-decision state is updated as follows.}

\textit{Transition (i) of $\tau_t$}. XB $j\in\tcb{\jtn}$ who is assigned to work sequence $m\in\mtj$ \tcb{with $a_m=t$} will have the $\tautpo(j)$ value updated to $b_m+1$. The transition function remains the same for XBs who are holding \tcb{(including those with only a future-starting sequence selected in the plan)} or are unavailable.
\begin{equation}
    \tautpo(j)=
    \begin{cases}
        \displaystyle\sum_{m\in\mtj\tcb{:\,a_m=t}}(b_m+1)\ytjm &\forall j\in\jtn:\displaystyle\sum_{m\in\mtj\tcb{:\,a_m=t}}\ytjm=1\\
        \tfrj+1 &\forall j\in\jtn: \displaystyle\sum_{m\in\mtj\tcb{:\,a_m=t}}\ytjm=0 \\
        \tfrj &\forall j\in\mathcal{J}\backslash\jtn
    \end{cases}
\end{equation}

\textit{Transition (i) of $\mathcal{D}_t$}. \tcb{As specified in Section~\ref{subsec:transition}, all pieces in $\dtn$ leave the open-work set after the decision, whether or not they are covered. Here, $\dasng$ consists of the pieces in selected XB sequences with $a_m=t$, including any that start later. Thus $\dtpo=\mathcal{D}_t\backslash(\dtn\cup\dasng)$; pieces selected only provisionally remain in $\dtpo$.}

\tcb{The OT transition in Section~\ref{subsec:transition} remains applicable because it counts only $\utik$ for pieces with $a_k=t$; the exogenous-information transitions also remain unchanged.}

\subsection{Calculation of value functions}\label{subsec:value_calculation}
We employ customized backward dynamic programming \tcb{algorithms} to \tcb{compute the individual future-value estimates for XBs and OTs in \eqref{eqn:post_simple}}.
We first \tcb{describe the calculation for} a particular XB $j\in\mathcal{J}$\tcb{, whose shift spans periods $s_j$ through $e_j$}. \tcb{The algorithm computes} $\vjtt$ \tcb{for all} $\tau\in[s_j+1,e_j+1]$ \tcb{and} $t\in[s_j,\tau-1]$. \tcb{Here,} $e_j+1$ \tcb{marks the end of the XB's shift, and} $t<\tau$ \tcb{reflects that a decision at} $t$ \tcb{leaves the XB next available in a later period.}

\tcb{For each XB, we store $\vjtt$ in a lookup table rather than impose the kind of prescribed functional relationship often used in parametric ADP approximations. This representation accommodates variation arising from the heterogeneity and temporal distribution of potentially open work pieces and from differences in which work-source outcomes remain unknown at $t$. Indexing the table only by $(\tau,t)$ also keeps the backward calculation manageable. If $m$ multi-piece work sources reveal their outcomes after $t$ but before $\tau$, and each retains a distinct piece eligible for the XB at or after $\tau$, their binary outcomes can generate up to $2^m$ sets of known future work before assignment decisions are considered. A direct backward calculation conditioned on those histories would face these possibilities, whereas our table has $L_j(L_j+1)/2$ entries for an XB shift of $L_j=e_j-s_j+1$ periods (e.g., 465 for the 30-period shifts in our case study). We compute the table offline for use in the online IP. The backward recursion uses source-opening probabilities and candidate sequences of eligible pieces within each source, as detailed in Algorithm~\ref{alg:bdp}.}

\begin{algorithm}
\caption{Calculation of \tcb{individual future-value estimates} for XB $j\in\mathcal{J}$}
\label{alg:bdp}
\begin{algorithmic}
    \STATE \textit{Step 0:} Set $\vjtt=0,\forall\tau\in[e_j-\delmin+2,e_j+1],t\in[s_j,\tau-1]$.
    \FOR{$\tau=(e_j-\delmin+1),(e_j-\delmin),\hdots,(s_j+2),(s_j+1)$}
        \STATE \textit{Step 1:} Obtain $\hjtau=\{h\in\mathcal{H}:\exists k\in h \ \texttt{s.t.} \ (a_k=\tau) \ \texttt{and} \ (k\in \mathcal{K}_j) \}$.
        \FOR{$t=(\tau-1),(\tau-2),\hdots,(s_j+1),(s_j)$}
            \STATE \textit{Step 2:} Obtain $\hjtaut=\{h\in\hjtau:a_h>t\} $.
            \STATE \textit{Step 3:} Denote $\vhold=\vjbar(\tau+1,t)$. Then, $\forall h\in\hjtaut$: obtain $\mhset$ and compute $v_h$ using \eqref{eqn:vh}.
            \STATE \textit{Step 4:} Rank $\{h\in\hjtaut:v_h>\vhold\}$ by descending $v_h$; denote the ranked set by $\hranked$.
            \STATE \textit{Step 5:} \tcb{If $\hranked=\emptyset$, set $\vjtt=\vhold$; otherwise,} compute $\vjtt$ according to \eqref{eqn:compv}.
        \ENDFOR
    \ENDFOR
    \RETURN $\vjtt, \forall\tau\in [s_j+1,e_j+1], \forall t\in [s_j,\tau-1]$.
\end{algorithmic}
\end{algorithm}

Recall that the minimum duration of all work pieces \tcb{in $\mathcal{K}$} is denoted by $\delmin$. If the XB becomes available at or after $e_j-\delmin+2$, their remaining shift interval is too short to cover another piece. \textit{Step 0} sets the \tcb{corresponding value estimates to} zero for all \tcb{$t\in[s_j,\tau-1]$}. We \tcb{then} start the \tcb{backward recursion} at $\tau=e_j-\delmin+1$ and decrement $\tau$ by one at each iteration until $\tau=s_j+1$. \tcb{In} \textit{Step 1}\tcb{, for} each $\tau$, we \tcb{identify the} work sources \tcb{containing} an eligible piece that \tcb{starts} at $\tau$. We call this piece the ``\tcb{piece of interest}'' \tcb{and denote} the set of these sources by $\hjtau$. \tcb{In} \textit{Step 2}\tcb{, we iterate over $t=\tau-1,\ldots,s_j$ and retain only sources satisfying $a_h>t$, obtaining} $\hjtaut$. \tcb{Since each source's outcome is revealed at its start time $a_h$, this subset contains the candidate} sources whose \tcb{outcomes remain unknown at} $t$.

In \textit{Step 3}, we first denote the \tcb{estimated} value of holding the XB \tcb{for one period by} $\vhold\tcb{=\vjbar(\tau+1,t)}$. Then, \tcb{for each} $h\in\hjtaut$, we calculate the \tcb{estimated} value of assigning the XB to \tcb{its highest-valued candidate} sequence, denoted by $v_h$. \tcb{Let} $\{k_1,k_2,\hdots,k_n\}$ \tcb{denote the pieces of $h$ that are eligible for XB $j$ and start at or after $\tau$, listed in chronological order. Here,} $k_1$ is the ``\tcb{piece of interest}'' \tcb{starting at $\tau$}, and $k_n$ is the last \tcb{eligible} piece. The candidate set of work sequences is defined as $\mhset=\big\{\{k_1\},\{k_1,k_2\},\{k_1,k_2,k_3\},\hdots,\{k_1,k_2,k_3,\hdots,k_n\}\big\}$.
\tcb{For the work sources considered here, pieces within a source do not overlap, so these candidate sequences are temporally feasible for the XB.}
The \tcb{estimated} value associated with $h\in\hjtaut$ is obtained by finding the highest-valued work sequence from $\mhset$:
\begin{equation} \label{eqn:vh}
    v_h=\max_{m\in \mhset} \tcb{\biggl\{}\sum_{k\in m}c_k+\vjbar(\tau=b_m+1,t)\tcb{\biggr\}}
\end{equation}
\tcb{Each candidate's value combines the reward from covering its pieces with the estimated continuation value when the XB becomes available again at $b_m+1$.}

In \textit{Step 4}, we rank work sources from highest to lowest $v_h$, discarding any that produces a value no greater than $\vhold$\tcb{, since holding has at least as large an estimated value}. The ranked work sources are denoted by $\hranked=\{h_1,h_2,\hdots,h_N\}$. \tcb{In} \textit{Step 5}, we compute $\vjtt$ \tcb{by combining these estimates with the source-opening probabilities}. \tcb{Under the assumed independence of source outcomes, $h_r$ is the highest-ranked open source when it is open and all higher-ranked sources are closed, which has probability} $p_{h_{\tcb{r}}}\prod_{\tcb{q}=1}^{\tcb{r}-1}(1-p_{h_{\tcb{q}}})$. \tcb{If} none of \tcb{the ranked sources is} open, \tcb{the holding estimate} $\vhold$ \tcb{contributes with probability $\prod_{r=1}^{N}(1-p_{h_r})$}. \tcb{If $N=0$, we set $\vjtt=\vhold$; otherwise, the estimate} is calculated as follows:
\begin{equation} \label{eqn:compv}
    \vjtt = v_{h_1}p_{h_1} + \sum_{r=2}^{N}v_{h_r}p_{h_r}\prod_{q=1}^{r-1}(1-p_{h_q}) + \vhold\prod_{q=1}^{N}(1-p_{h_q})
\end{equation}

\tcb{To calculate} \tcb{individual future-value estimates for} OTs\tcb{, we introduce an auxiliary value $F_i(\phi,B)$ to account for the option to wait for a later piece of a newly revealed source. Here, $B$ contains eligible pieces revealed open within the calculation that start at or after $s_i+\phi$. The value $F_i(\phi,B)$ is the optimal expected net reward for OT $i$ just before source outcomes at $s_i+\phi$ are observed. Each lookup entry starts with $B=\varnothing$, excluding work already known at the actual decision epoch, as in Section~\ref{subsec:value_approximation}.}

\tcb{The} value of assigning each piece $k\in\tcb{\mathcal{K}_i}$ to the OT \tcb{is given by} the following equation:
\begin{equation} \label{eqn:vk}
    v_k=c_k-\costot\Big[ a_k+\delta_k-(s_i+\delot)\Big]^+
\end{equation}

\tcb{For $\phi\in[1,\delot]$, let $B^+$ be the random set containing $B$ and all eligible pieces of sources revealed to be open at $s_i+\phi$. Its contents depend on the source outcomes observed at that period. After observing these outcomes, the OT may take a piece starting at that period or hold and retain the later pieces:}
\begin{equation} \label{eqn:ot_general}
    \tcb{\begin{aligned}
    F_i(\phi,B)=\mathbb{E}\Big[\max\Big\{&
    \max\bigl(\{0\}\cup\{v_k:k\in B^+,\ a_k=s_i+\phi\}\bigr), F_i\bigl(\phi+1,\{k\in B^+:a_k>s_i+\phi\}\bigr)
    \Big\}\Big].
    \end{aligned}}
\end{equation}
\tcb{The expectation in \eqref{eqn:ot_general} averages the maximized value over the possible combinations of source outcomes at $s_i+\phi$. The weights in this average follow from each source $h$ starting at that period opening independently with probability $p_h$. With terminal value $F_i(\delot+1,\varnothing)=0$, backward recursion gives $\viphi=F_i(\phi,\varnothing)$ for $\phi\in[1,\delot+1]$; computed auxiliary values are stored for reuse.}

\tcb{If computing these values becomes computationally expensive, a cheaper approximation can consider only the first pieces of newly revealed sources, omitting their later-piece opportunities. Algorithm~\ref{alg:bdpot} summarizes this scalar recursion.}

\begin{algorithm}
\caption{\tcb{Scalar recursion for} \tcb{individual future-value estimates} for OT $i\in\mathcal{I}$}
\label{alg:bdpot}
\begin{algorithmic}
    \STATE \textit{Step 0:} Set $\vibar(\phi=\delot+1)=0$.
    \FOR{$\phi=(\delot),(\delot-1),\hdots,2,1$}
        \STATE \textit{Step 1:} Obtain $\kiphi=\{k\in\mathcal{K}_i: \tcb{k=k^\texttt{one}_{h_k},}\ a_k=s_i+\phi\}$.
        \STATE \textit{Step 2:} Denote $\vhold=\vibar(\phi+1)$. Then, $\forall k\in\kiphi$: compute $v_k$ using \eqref{eqn:vk}.
        \STATE \textit{Step 3:} Rank $\{k\in\kiphi:v_k>\vhold\}$ by descending $v_k$; denote the ranked set by $\kranked$.
        \STATE \textit{Step 4:} \tcb{If $\kranked=\varnothing$, set $\viphi=\vhold$; otherwise,} compute $\viphi$ according to \eqref{eqn:compvot}.
    \ENDFOR
    \RETURN $\viphi, \forall\phi\in [1,\delot+1]$.
\end{algorithmic}
\end{algorithm}

\tcb{For each $\phi$, we obtain the holding value $\vhold$ and} rank \tcb{the eligible first} work pieces by decreasing $v_k$ value. Similar to Algorithm \ref{alg:bdp} for XBs, pieces with $v_k\leq\vhold$ are removed from the ranked set because it is \tcb{at least as valuable} to hold the OT \tcb{as} to assign them such pieces \tcb{under this approximation}. The \tcb{remaining} ranked set is denoted $\kranked=\{k_1,k_2,\hdots,k_N\}$\tcb{, with} respective parent sources $h_1,h_2,\hdots,h_N$. \tcb{If $N=0$, we set $\viphi=\vhold$. Otherwise,} we calculate $\viphi$ \tcb{by combining the assignment and holding values with the source-opening probabilities,} following the same logic as Step 5 of Algorithm \ref{alg:bdp}:
\begin{equation} \label{eqn:compvot}
    \viphi = v_{k_1}p_{h_1} + \sum_{r=2}^{N}v_{k_r}p_{h_r}\prod_{q=1}^{r-1}(1-p_{h_q}) + \vhold\prod_{q=1}^{N}(1-p_{h_q})
\end{equation}
\tcb{The scalar recursion is exact for the individual-OT calculation when pieces within each source do not overlap and $\delot\leq\delmin$. Any later piece of a source revealed at $a_h\geq s_i+1$ starts at or after $a_h+\delmin$, which exceeds $s_i+\delot$. Consequently, no later piece needs to be retained: \eqref{eqn:ot_general} is evaluated only at $B=\varnothing$ and reduces to \eqref{eqn:compvot}. These conditions hold in our case study (Section~\ref{sec:ot_experiments}), where the scalar form therefore evaluates the general recursion exactly.}

The \tcb{individual} value \tcb{estimates} for XBs and OTs serve as \tcb{inputs} to the \tcb{IP}. These \tcb{estimates} can be \tcb{computed offline} given the daily bus schedule, open work probabilities, and \tcb{specified XB shifts and OT} report times. \tcb{For example, if} report times are finalized around noon \tcb{on the preceding} day\tcb{, the corresponding estimates} can \tcb{be computed afterward for use on} the day of operation.

Alternatively, \tcb{we} can \tcb{precompute estimates for candidate XB shifts and OT report times} at the start of \tcb{a} planning period. \tcb{These estimates can be reused across workdays with the same work-piece and work-source data, OT cost coefficient, and minimum-pay duration.} For example, \tcb{weekdays within a three-month scheduling period may share these inputs. Once the roster for} a particular weekday \tcb{is finalized}, \tcb{the} corresponding \tcb{estimates} can be retrieved \tcb{for use in the IP}. \tcb{Separate} sets of \tcb{estimates} can be \tcb{computed} for other groups of workdays\tcb{, such as} Saturdays \tcb{and} Sundays.

\subsection{Algorithmic adjustments to account for overestimation of value functions}
\tcb{The future-value terms in \eqref{eqn:obj} use the additive approximation in \eqref{eqn:post_simple}, which combines values estimated separately for each XB and OT. Each individual calculation in Section~\ref{subsec:value_calculation} considers the operator's opportunities without competition from other operators. Consequently, several operators' estimates may reflect the same future work, although each piece can be covered only once. Their summed values can therefore overestimate the future reward attainable by the roster.} We call this the ``oversupply effect''\tcb{. Such overestimation can make retaining operators for future work appear more attractive than warranted. We therefore introduce an adjustment to the holding values of operators with no selected work in the provisional plan.}

We first describe the \tcb{adjustment} for XBs. At decision epoch $t\in\mathcal{T}$, \tcb{for an available} XB $j\in\tcb{\jtn}$, \tcb{let} $\vbarhld_{tj}\tcb{=}\vjbar(t+1,t)$ \tcb{denote its individual future-value estimate after holding for one period. To represent competition among XBs for future work, we introduce a hypothetical} ``position'' $n$\tcb{, defined as the} number of other XBs \tcb{assumed to have} priority over XB $j$. \tcb{Position is a modeling device for estimating competition and does not prescribe the order of actual assignments. A higher position represents greater competition for future work, motivating a larger reduction in the XB's nominal holding value. We denote the adjusted future-value estimate by} $\vhattj$. \tcb{Motivated by computational tractability, we use a linear approximation for the position effect, with $\lbartj$ giving the change in the affine holding-value estimate per unit increase in position. This choice allows the adjustment to be incorporated into the online IP through auxiliary variables and linear constraints.}

\begin{figure}[!ht]
    \centering
    \includegraphics[width=0.65\textwidth]{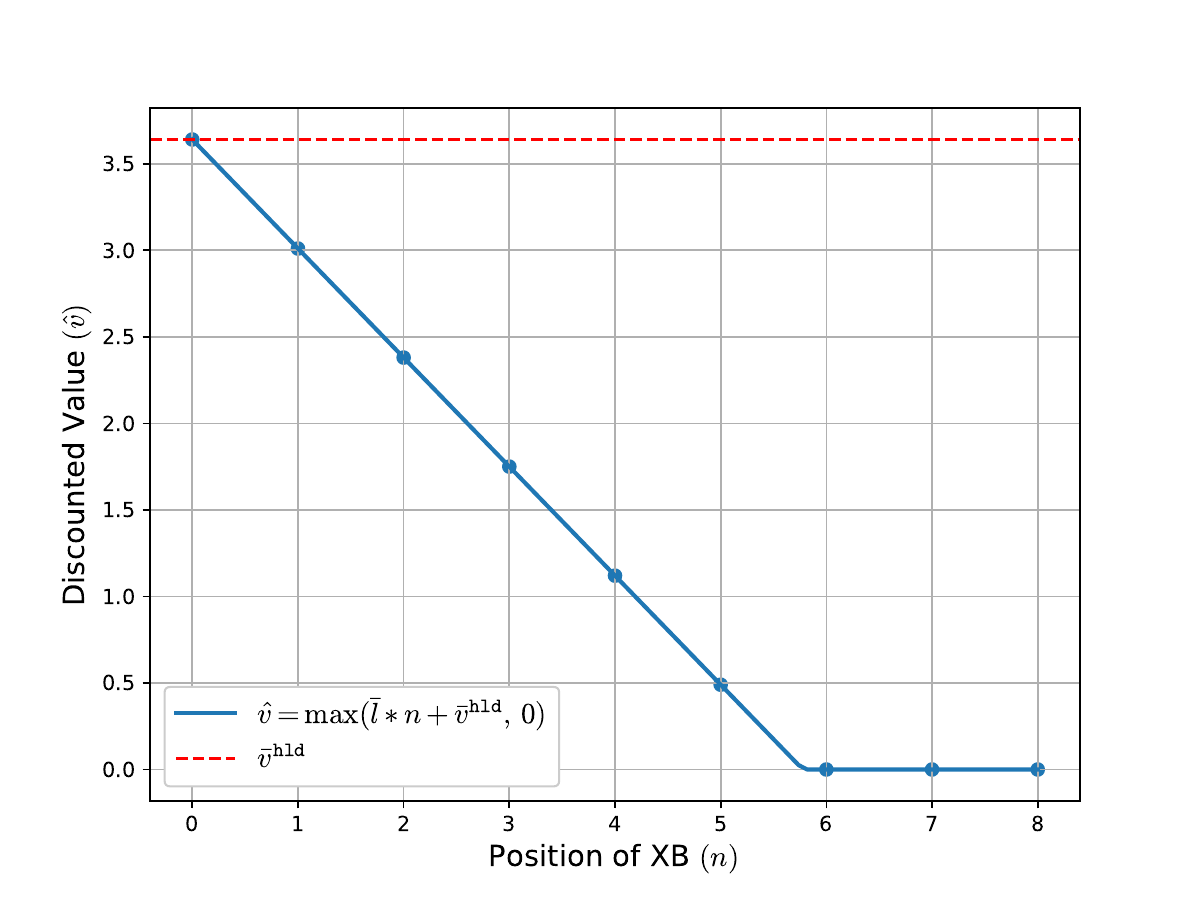}
    \caption{Discounted Future Value for Holding due to Oversupply of Operators}
    \label{fig:disc_val}
\end{figure}

\tcb{Figure~\ref{fig:disc_val} illustrates the adjusted future-value estimate, with the indices} $t$ and $j$ \tcb{suppressed for clarity}. At position \tcb{$n=0$}, \tcb{meaning that the XB is assumed to have} top priority \tcb{for future work, the adjusted estimate equals the nominal holding value} $\vbarhld$. \tcb{We use} \tcb{$\bar{l}\leq0$,} \tcb{so each additional position reduces the affine estimate by} $\tcb{-}\bar{l}$. \tcb{When $\vbarhld+\bar{l} \cdot n\leq0$, we set the adjusted estimate to zero. This lower bound reflects that an XB can remain idle for the rest of its shift and earn zero reward.} Thus, the \tcb{adjusted future-value estimate for a} holding XB \tcb{at} position $n$ is expressed as $\hat{v} = \max(\bar{l} \cdot n + \bar{v}^\mathtt{hld},\, 0)$. \tcb{The estimation of} $\bar{l}$ \tcb{is described in Section~\ref{subsec:loss_estimation}}.

\tcb{The position-based adjustment must account for the different combinations of XBs that may be left without selected sequences in the provisional plan. Within each such group, the actual order of future assignments depends on the work revealed and subsequent decisions. Jointly choosing work assignments and hypothetical priority positions would introduce an additional ordering problem into the online optimization. To avoid this additional complexity in repeated online solves, we evaluate each XB's adjusted holding value at the \textit{average} position within the group. This approximation represents competition through the size of the group, which is determined by the plan's assignment decisions.} 
\tcb{Let} $\nholdt$ \tcb{denote the} number of \tcb{currently available XBs with at least $\delmin$ periods remaining in their shifts and no selected sequence in the provisional plan:}
\begin{equation}
    \nholdt=\sum_{\substack{j\in\jtn\tcb{:}\\\tcb{t\leq e_j-\delmin+1}}}(1-\sum_{m\in\mtj}\ytjm)
\end{equation}
\tcb{XBs with fewer than $\delmin$ periods remaining cannot cover another piece and have zero nominal holding value. An XB selected only for a future-starting sequence holds physically at $t$ but is not counted in $\nholdt$, since this quantity is used to score the provisional plan. For $\nholdt\geq1$, the hypothetical positions within the group range from $0$ to $\nholdt-1$. We therefore use the} average position
\begin{equation}
    \navg = \frac{\nholdt-1}{2}, \quad \nholdt\geq 1
\end{equation}
\tcb{for each XB in the group.} Note that $\navg$ is not necessarily an integer under this definition. For instance, when there are two XBs \tcb{in the group}, the \tcb{hypothetical positions are $0$ and $1$, and we evaluate both XBs' adjusted holding values at $\navg=0.5$.}

We \tcb{incorporate this average-position adjustment into the online IP} (\ref{eqn:obj})-(\ref{eqn:cons26}) \tcb{by introducing} auxiliary variables, $\wtj\geq0, \forall j\in \jtn$, \tcb{for the adjusted holding-value contributions. The future-value terms associated with selected sequences remain unchanged.} The second line of the objective function (\ref{eqn:obj}) now reads
\[\sum_{j\in\jtn} \bigg( \sum_{m\in\mtj}\ytjm\vjbar(\tau=b_m+1,t) + \wtj
\bigg)\]
\tcb{The variable $\wtj$ represents the adjusted holding value $\vhattj$ when XB $j$ has no selected sequence in the provisional plan, and is zero otherwise.} \tcb{Since $\wtj$ enters the maximization objective with a positive coefficient,} the following \tcb{upper bounds enforce this relation at optimality}:
\begin{align}
    \wtj&\leq\vhattj &&\forall j\in\jtn \label{eqn:w_and_vhat} \\
    \wtj&\leq\vbarhld_{tj} \cdot (1-\sum_{m\in\mtj}\ytjm) &&\forall j\in\jtn
\end{align}
To \tcb{incorporate} the expression $\vhattj = \max(\lbartj \cdot \navg + \vbarhld_{tj},\, 0)$ \tcb{into the IP, we introduce} a binary variable $\qtj$. The following inequalities \tcb{activate the affine upper bound when $\qtj=1$ and force $\vhattj=0$ when $\qtj=0$}:
\begin{align}
    \vhattj&\leq\vbarhld_{tj} \cdot \qtj &&\forall j\in\jtn \label{eqn:vhat_1} \\
    \vhattj&\leq\vbarhld_{tj}+\lbartj \cdot \navg-\lbartj \cdot \frac{\nmax-1}{2}(1-\qtj) &&\forall j\in\jtn \label{eqn:vhat_2} \\
    \vhattj&\geq 0 &&\forall j\in\jtn \label{eqn:vhat_3}
\end{align}
\tcb{Here,} $\nmax$ denotes \tcb{an upper bound on $\nholdt$}. \tcb{It} can be chosen as the maximum number of XBs \tcb{simultaneously on duty, calculated from their} report times and \tcb{shift} durations. \tcb{This ensures $\navg\leq(\nmax-1)/2$ whenever $\nholdt\geq1$. Thus, when $\qtj=0$, the right-hand side of \eqref{eqn:vhat_2} is nonnegative, allowing $\vhattj=0$. Combining $\wtj\leq\vhattj$ with these upper bounds allows us to eliminate $\vhattj$ by imposing the bounds directly on} $\wtj$ in the final formulation.

\tcb{We apply the same form of position-based adjustment to} OTs. At decision epoch $t\in\mathcal{T}$, \tcb{let} $\vbarhldot_{ti} = \vibar(t+1-s_i)$ \tcb{denote the} nominal \tcb{holding-value estimate} for an available OT $i\in\itn$. \tcb{As} for XBs, $\lbarotti\tcb{\leq0}$ \tcb{gives} the \tcb{change in the affine holding-value estimate per unit increase in position}, and $\vhatotti$ \tcb{denotes the adjusted holding-value estimate. The position used for an OT represents the number of XBs assumed to have priority for future work.}
From numerical experiments with a myopic decision rule, we find that \tcb{giving XBs priority over OTs yields higher average net rewards than the reverse ordering (Section~\ref{sec:ot_experiments}). Motivated by this comparison, we place OTs not selected in the plan behind all XBs counted in $\nholdt$ for the holding-value calculation}, such that $n=\nholdt$. The \tcb{adjusted holding-value estimate} for the OT is therefore expressed as $\vhatotti = \max(\lbarotti \cdot \nholdt + \vbarhldot_{ti},\, 0)$.

\tcb{This adjustment omits competition among OTs, a simplification motivated by their shorter standby window relative to XB shifts in our case study.} We use the variables $\wotti\tcb{\geq0}, \forall i\in\itn$ \tcb{for the adjusted OT holding-value contributions, with binary variables} $\qotti$ \tcb{introduced as for XBs. At optimality, $\wotti$ equals $\vhatotti$ when OT $i$ has no selected piece in the plan and is zero otherwise.} \tcb{An OT selected only for a later-starting piece remains on standby physically but does not receive this plan-level holding term.}

\tcb{Combining these adjustments and substituting the expressions for} $\navg$ and $\nholdt$ \tcb{yields the final IP, which is solved online to implement} the approximate policy:
\begin{align}
    \max \ &\sum_{k\in\mathcal{D}_t} c_k\ztk -\costot \sum_{i\in \itn}\sum_{k\in \mathcal{D}_t\cap\mathcal{K}_i} \bigg[ a_k+\delta_k-(s_i+\delot)\bigg]^+ \utik\notag\\
    + &\sum_{j\in\jtn} \bigg( \sum_{m\in\mtj}\ytjm\vjbar(\tau=b_m+1,t) + 
    \wtj \bigg)\notag\\
    + &\sum_{i\in\itn}\wotti \label{eqn:obj3} \\
    \textrm{s.t.} \ &\sum_{m\in\mtj}\ytjm\leq 1 &&\forall j\in\jtn \label{eqn:cons31} \\
    &\sum_{k\in\mathcal{D}_t\cap \mathcal{K}_i}\utik\leq 1 &&\forall i\in\itn \label{eqn:cons32} \\
    &\sum_{j\in\jtn}\sum_{\substack{m\in\mtj:\\k\in m}} \ytjm + \sum_{\substack{i\in\itn:\\k\in\mathcal{K}_i}} \utik = \ztk &&\forall k\in\mathcal{D}_t \label{eqn:cons33} \\
    &\wtj\leq\vbarhld_{tj} \cdot (1-\sum_{m\in\mtj}\ytjm) &&\forall j\in\jtn \label{eqn:cons34} \\
    &\wtj\leq\vbarhld_{tj} \cdot \qtj &&\forall j\in\jtn \label{eqn:cons35} \\
    &\wtj\leq\vbarhld_{tj}+\lbartj \cdot \frac{\displaystyle\sum_{\substack{\tcb{j'}\in\jtn\tcb{:}\\\tcb{t\leq e_{j'}-\delmin+1}}}(1-\sum_{m\in\mathcal{M}_t(\tcb{j'})}y^t_{\tcb{j'}m})-1}{2}-\lbartj \cdot \frac{\nmax-1}{2}(1-\qtj) &&\-\hspace{+0.05em}\forall j\in\jtn \label{eqn:cons36} \\
    &\wotti\leq\vbarhldot_{ti} \cdot (1-\sum_{k\in\mathcal{D}_t\cap\mathcal{K}_i}\utik) &&\-\hspace{-0em}\forall i\in\itn \label{eqn:cons37} \\
    &\wotti\leq\vbarhldot_{ti} \cdot \qotti &&\-\hspace{-0em}\forall i\in\itn \label{eqn:cons38} \\
    &\wotti\leq\vbarhldot_{ti}+\lbarotti \cdot \sum_{\substack{j\in\jtn\tcb{:}\\\tcb{t\leq e_j-\delmin+1}}}(1-\sum_{m\in\mtj}\ytjm)-\lbarotti \cdot\nmax(1-\qotti) &&\-\hspace{-0em}\forall i\in\itn \label{eqn:cons39} \\
    &\ytjm\in \{0,1\} &&\-\hspace{-3.0em}\forall j\in\jtn, m\in\mtj \label{eqn:cons310} \\
    &\ztk\in \{0,1\} &&\-\hspace{-3.0em}\forall k\in\mathcal{D}_t \label{eqn:cons311} \\
    &\utik\in \{0,1\} &&\-\hspace{-3.0em}\forall i\in\itn, k\in\mathcal{D}_t \cap \mathcal{K}_i \label{eqn:cons312} \\
    &\wtj\geq 0 && \-\hspace{-3.0em}\forall j\in\jtn \label{eqn:cons313} \\
    &\wotti\geq 0 &&\-\hspace{-3.0em}\forall i\in\itn \label{eqn:cons314} \\
    &\qtj\in \{0,1\} &&\-\hspace{-3.0em}\forall j\in\jtn \label{eqn:cons315} \\
    &\qotti\in \{0,1\} &&\-\hspace{-3.0em}\forall i\in\itn \label{eqn:cons316}
\end{align}

In the objective function (\ref{eqn:obj3}), \tcb{adjusted holding-value contributions for XBs in $\jtn$ and OTs in $\itn$ with no selected work in the provisional plan} are represented by $w$ and $w^\texttt{ot}$ variables, respectively. \tcb{As in \eqref{eqn:obj}, the first line scores known work selected in the provisional plan; realized reward depends only on executed assignments.} \tcb{Constraints} (\ref{eqn:cons34})-(\ref{eqn:cons36}) incorporate \tcb{the holding-value adjustments} for XBs, and (\ref{eqn:cons37})-(\ref{eqn:cons39}) are analogous constraints for OTs. \tcb{When $\nholdt=0$, all available XBs either have selected sequences or have zero nominal holding value, so the constraints enforce $\wtj=0$ without requiring an interpretation of an average position.} Lastly, (\ref{eqn:cons313})-(\ref{eqn:cons316}) \tcb{specify the domains of} the auxiliary variables $w$, $w^\texttt{ot}$, $q$, and $q^\texttt{ot}$.

\subsection{Estimation of loss parameters}\label{subsec:loss_estimation}
We \tcb{estimate the} loss parameters $\bar{l}$ and $\bar{l}^\texttt{ot}$ \tcb{for} individual XBs and OTs, respectively, \tcb{using simulated work-source opening outcomes}. As with the XB and OT \tcb{value estimates}, these loss parameters \tcb{can be computed offline, either daily or at the start of a planning period, and reused while the relevant inputs remain unchanged}.

\tcb{To estimate $\lbartj$ for an XB $j$ holding at $t$, we simulate opening outcomes for eligible work sources revealed after $t$. In each scenario, we approximate the effect of other XBs having priority as follows. At position $n=0$, the XB has top priority, which we represent by an assignment to the first eligible source to open. At position $n=1$, one XB has priority ahead of it, so we consider an assignment to the second such source, as if the first opportunity had been taken. More generally, position $n$ is represented by an assignment to the $(n+1)$-th eligible source to open. We value each assignment by its reward plus the XB's estimated future value after completion. We average these values across scenarios at each hypothetical position. We then fit a linear regression of these averages against position $n$, using $n=0,1,\ldots,N^{\texttt{fit}}$, where $N^{\texttt{fit}}$ denotes the maximum position initially considered for fitting. The fitted coefficient of position is used to estimate $\lbartj$, as detailed in Algorithm~\ref{alg:loss}.}

\begin{algorithm}[!ht]
\caption{Calculation of the loss parameter for XB $j\in\mathcal{J}$}
\label{alg:loss}
\begin{algorithmic}
    \STATE \tcb{\textbf{Fitting parameters:} integer $N^{\texttt{fit}}\geq1$ and threshold $0\leq\rthr\leq1$.}
    \STATE \textit{Step 0:} Set $\lbartj=0,\forall t\in[\tcb{s_j},e_j]$. \tcb{Retain zero if no fit is accepted.}
    \STATE \textit{Step 1:} Obtain $\hjelig=\{h\in\mathcal{H}: \kfirst\in\mathcal{K}_j \}$.
    \STATE \textit{Step 2:} $\forall h\in\hjelig$: obtain $\tcb{m_j(h)}$ and calculate $\vjbar(h)$ according to (\ref{eqn:vhloss}).
    \FOR{$t=s_j,s_j+1,s_j+2,\hdots,e_j-\delmin$}
        \STATE \textit{Step 3:} Obtain $\hjaftert=\{h\in\hjelig: a_h>t\}$.
        \STATE \textit{Step 4:} Simulate random scenarios (denoted by $\Xi$) of openness outcomes for work sources in $\hjaftert$, \tcb{sampling each source once per scenario with opening probability $p_h$}. The set of open work sources in each scenario $\xi\in\Xi$ is denoted by $\hopenxit$.
        \FOR{$\xi\in\Xi$}
            \STATE \textit{Step 5:} Rank the work sources in $\hopenxit$ by ascending start time ($a_h$), with ties broken randomly. Denote the ranked work sources as $\{h_{\xi,0}, h_{\xi,1},\hdots,h_{\xi,n},\hdots,h_{\xi,(|\tcb{\hopenxit}|-1)} \}$.
            \STATE \textit{Step 6:} Get $\vjbar(\xi,n)=\vjbar(h_{\xi,n}), \forall n\in \tcb{\{0,1,\ldots,N^{\texttt{fit}}\}}$. For values of $n$ where $h_{\xi,n}$ does not exist, $\vjbar(\xi,n)=0$.
        \ENDFOR
        \STATE \textit{Step 7}: Calculate $\vtilde_j(n)=\sum_{\xi\in\Xi}\vjbar(\xi,n)/|\Xi|, \forall n\in \tcb{\{0,1,\ldots,N^{\texttt{fit}}\}}$
        \FOR{$\ntop = \tcb{N^{\texttt{fit}}},\tcb{N^{\texttt{fit}}}-1,\tcb{N^{\texttt{fit}}}-2,\hdots,2,1$}
            \STATE \tcb{If $\vtilde_j(0)=\cdots=\vtilde_j(\ntop)$, \textbf{break}.}
            \STATE \textit{Step 8}: Perform linear regression \tcb{with a fitted intercept} on $\vtilde_j(n)$ against $n$, for $n\in \tcb{\{}0,1,\hdots,\ntop\tcb{\}}$. Record its $r^2$ and slope ($m$).
            \IF{$r^2\geq \rthr$}
                \STATE $\lbartj=\tcb{\min\{m,0\}}$
                \STATE \textbf{break}
            \ENDIF
        \ENDFOR
    \ENDFOR
    \RETURN $\lbartj, \forall t\in [s_j,e_j]$.
\end{algorithmic}
\end{algorithm}

Algorithm \ref{alg:loss} \tcb{estimates $\lbartj$} for XB $j\in\mathcal{J}$ \tcb{at each period} $t\in[s_j,e_j]$. \textit{Step 0} \tcb{initializes all these coefficients to zero, which is retained whenever no fit is accepted.} \tcb{In particular, for} $t \geq e_j-\delmin+1$, \tcb{the remaining shift after holding} is too short to cover any work piece, \tcb{so} $\vbarhld_{tj}=0$ \tcb{and no holding-value correction is needed}. \textit{Steps 1 \tcb{and} 2} \tcb{calculate an assignment value for each eligible work source}. We define the set $\hjelig$, which consists of all work sources whose first component work piece ($\kfirst$) is eligible for the XB. \tcb{For each such source $h$, we form a sequence $m_j(h)$ containing all its pieces eligible for the XB, ordered by start time.} The associated \tcb{source value combines the assignment reward with the estimated continuation value after completing this sequence}:
\begin{equation} \label{eqn:vhloss}
    \vjbar(h)=\sum_{k\in \tcb{m_j(h)}}c_k+\vjbar(\tau=b_{\tcb{m_j(h)}}+1,t=a_h)
\end{equation}
\tcb{We then estimate} $\lbartj$ \tcb{for} $t=s_j\tcb{,\ldots,e_j-\delmin}$. In \textit{Step 3}, we \tcb{form} $\hjaftert$, a subset of $\hjelig$ \tcb{containing sources revealed strictly after} $t$. In \textit{Step 4}, \tcb{we simulate opening outcomes for these sources, sampling each source once per scenario according to a Bernoulli trial with opening probability} $p_h$. The set of random scenarios is denoted by $\Xi$, indexed by $\xi$. \tcb{In \textit{Step 5}, for each scenario $\xi\in\Xi$, the revealed open work sources are ranked in ascending order of their start times, with ties broken randomly. In \textit{Step 6}, position $n=0,1,\ldots,N^{\texttt{fit}}$ receives the source value $\vjbar(\xi,n)=\vjbar(h_{\xi,n})$, or zero if no corresponding source exists.} In \textit{Step 7}, \tcb{we average these values across scenarios to obtain} $\vtilde_j(n)$\tcb{, suppressing the fixed decision epoch $t$ in this notation. Zero outcomes remain in these averages: they represent scenarios with too few open sources to provide work at the corresponding position.}

\begin{figure}[htbp]
  \centering
  \begin{subfigure}{0.5\textwidth}
    \centering
    \includegraphics[width=\linewidth]{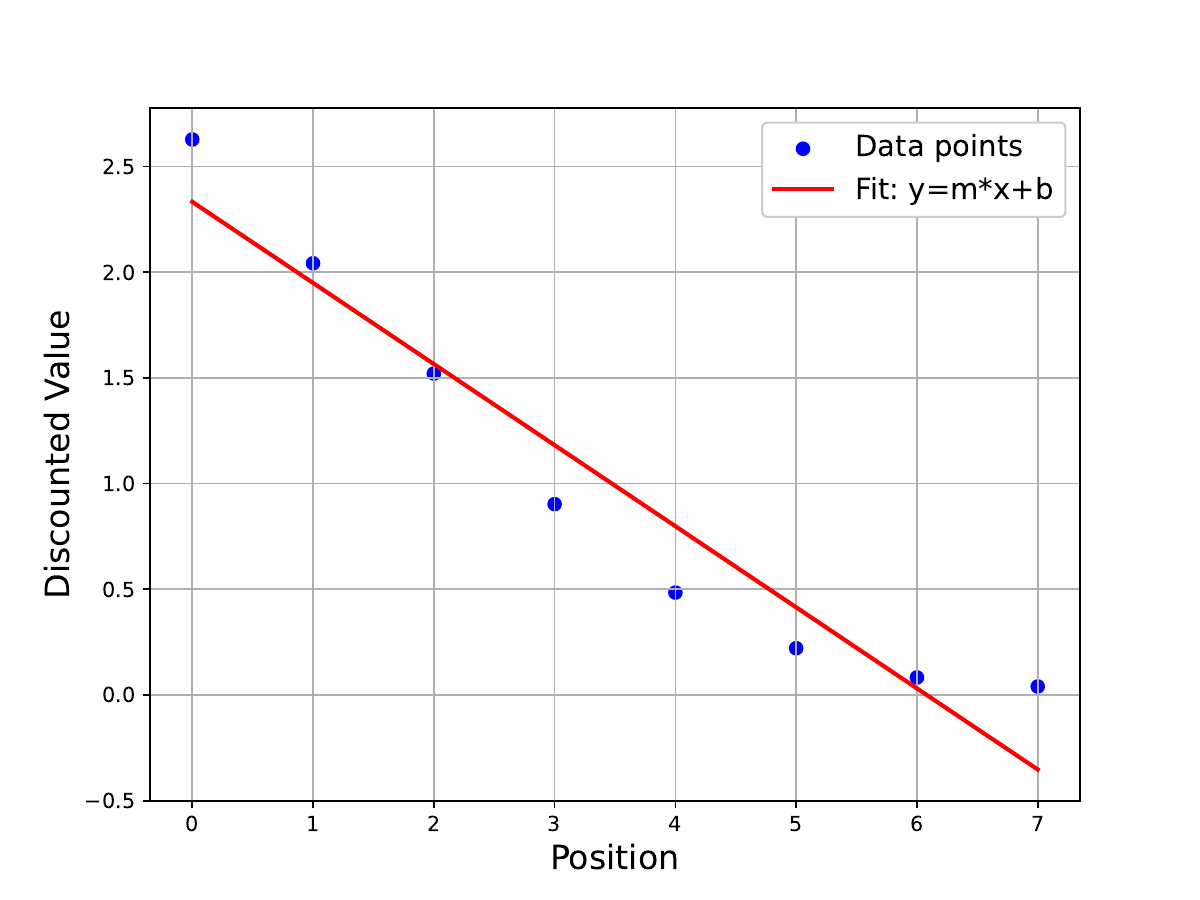}
    \caption{All data points; $r^2=0.930$}
  \end{subfigure}\hfill
  \begin{subfigure}{0.5\textwidth}
    \centering
    \includegraphics[width=\linewidth]{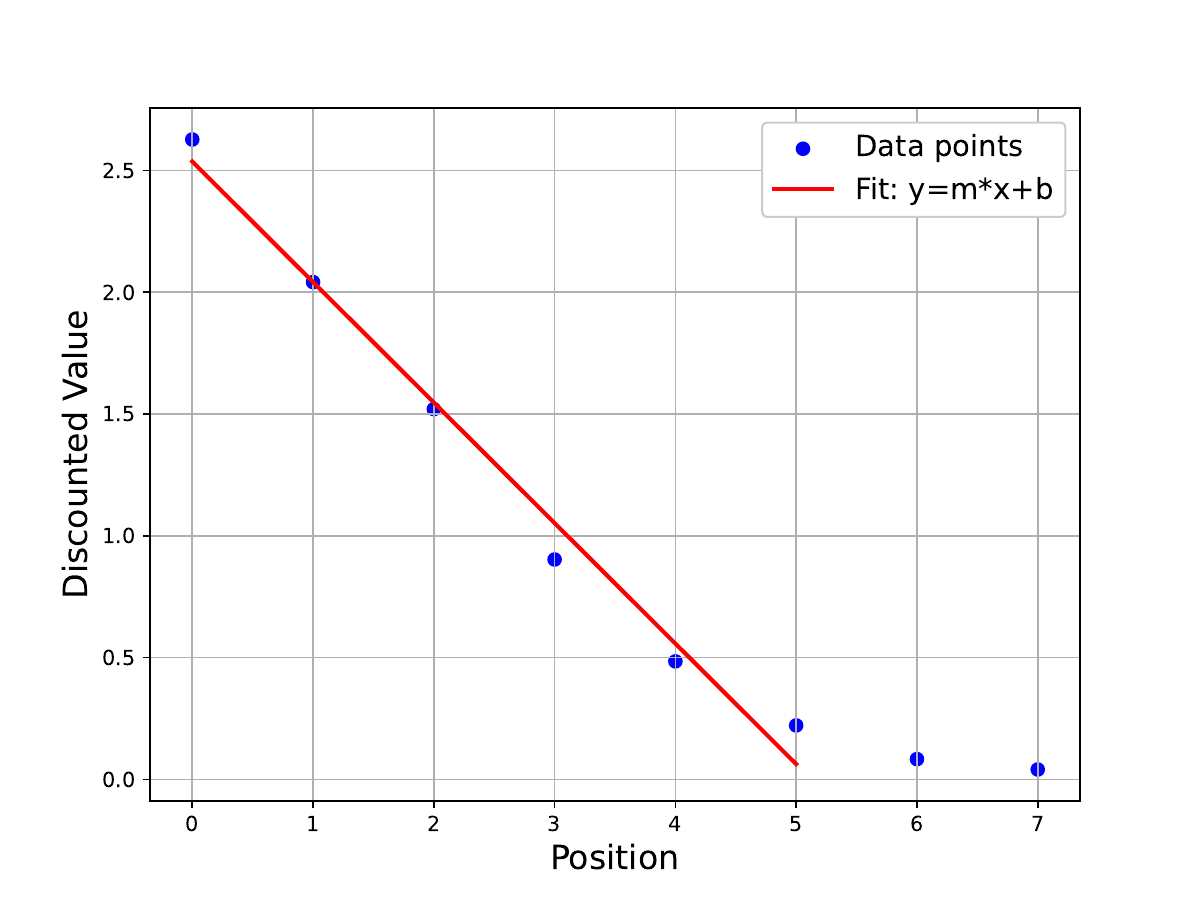}
    \caption{$n\in \tcb{\{}0,1,2,\hdots,5\tcb{\}}$\tcb{;} $r^2=0.986$}
  \end{subfigure}
  \caption{Example Linear Regression for Estimating the Loss Parameter}
  \label{fig:ols_comp}
\end{figure}

\tcb{In \textit{Step 8}, we fit a linear regression of} $\vtilde_j(n)$ \tcb{against position $n$, estimating both an intercept and a coefficient of position.} Initially, all $\vtilde_j(n)$ values from $n=0$ to $n=\tcb{N^{\texttt{fit}}}$ are included. \tcb{At later positions, the averages may flatten near zero because few scenarios contain enough open sources, making a single linear relationship less suitable over the full range. We therefore decrease the maximum included position $\ntop$ one at a time until the fitted line satisfies $r^2\geq\rthr$. This selects the longest initial range of positions meeting the fit criterion, with at least two positions included. For an accepted fit with coefficient $m$, we set $\lbartj=\min\{m,0\}$ to preserve the nonpositive correction convention. If the values in the remaining range are constant or no fit is accepted, the initialized value of zero is retained. Only the coefficient of position is used in the holding-value adjustment; its value at position zero remains the nominal holding value $\vbarhld_{tj}$, rather than the fitted intercept.} 
\tcb{Figure~\ref{fig:ols_comp} illustrates the fitting-range selection for a threshold $\rthr=0.97$. In panel (a), fitting all positions gives $r^2=0.930$, below the threshold. In panel (b), fitting positions $n=0,1,\ldots,5$ gives $r^2=0.986$, which satisfies the criterion. The coefficient of the latter fitted line is used to estimate $\lbartj$ for the illustrated $(j,t)$ pair.}

Algorithm \ref{alg:loss_ot} \tcb{estimates $\lbarotti$} for OT $i\in\mathcal{I}$\tcb{, following the same scenario-generation and regression procedure as for XBs.} \tcb{Each eligible work source is scored by the highest net reward among its pieces eligible for the OT, since the OT can receive at most one assignment. Here, position $n$ represents the number of XBs assumed to have priority over the OT. The fitting-range selection, nonpositive coefficient convention, and zero fallback follow Algorithm~\ref{alg:loss}. Only the fitted coefficient of position enters the adjustment; the value at position zero remains the nominal OT holding value $\vbarhldot_{ti}$. In the MiWay OT experiments, the conditions established after Algorithm~\ref{alg:bdpot} ensure that only the first piece of a source revealed after $t$ can be eligible.}

\begin{algorithm}[!ht]
\caption{Calculation of the loss parameter for OT $i\in\mathcal{I}$}
\label{alg:loss_ot}
\begin{algorithmic}
    \STATE \tcb{\textbf{Fitting parameters:} integer $N^{\texttt{fit}}\geq1$ and threshold $0\leq\rthr\leq1$.}
    \STATE \textit{Step 0:} Set $\lbarotti=0, \ \tcb{\forall t\in[s_i,s_i+\delot]}$. \tcb{Retain zero if no fit is accepted.}
    \STATE \textit{Step 1:} Obtain $\tcb{\mathcal{H}_i=\{h\in\mathcal{H}:\exists k\in\mathcal{K}_i,\ h_k=h\}}$.
    \STATE \textit{Step 2:} $\forall \tcb{h\in\mathcal{H}_i}$: compute $\vibar(\tcb{h})=\tcb{\max_{k\in\mathcal{K}_i:\,h_k=h}v_k}$ \tcb{using \eqref{eqn:vk}}.
    \FOR{$t=s_i,s_i+1,s_i+2,\hdots,s_i+\delot-1$}
        \STATE \textit{Step 3:} Obtain $\tcb{\mathcal{H}_i^>(t)=\{h\in\mathcal{H}_i:a_h>t\}}$.
        \STATE \textit{Step 4:} Simulate random scenarios (denoted by $\Xi$) of openness outcomes for \tcb{work sources in $\mathcal{H}_i^>(t)$, sampling each source once per scenario with opening probability $p_h$}. The set of open work \tcb{sources} in each scenario $\xi\in\Xi$ is denoted by $\tcb{\hopenxit}$.
        \FOR{$\xi\in\Xi$}
            \STATE \textit{Step 5:} Rank the work \tcb{sources} in $\tcb{\hopenxit}$ by ascending start time ($\tcb{a_h}$), with ties broken randomly. Denote the ranked work sources as $\tcb{\{h_{\xi,0},h_{\xi,1},\hdots,h_{\xi,(|\hopenxit|-1)}\}}$.
            \STATE \textit{Step 6:} Get $\vibar(\xi,n)=\vibar(\tcb{h_{\xi,n}}), \forall n\in \tcb{\{0,1,\ldots,N^{\texttt{fit}}\}}$. For values of $n$ where $\tcb{h_{\xi,n}}$ does not exist, $\vibar(\xi,n)=0$.
        \ENDFOR
        \STATE \textit{Step 7}: Calculate $\vtilde_i(n)=\sum_{\xi\in\Xi}\vibar(\xi,n)/|\Xi|, \forall n\in \tcb{\{0,1,\ldots,N^{\texttt{fit}}\}}$
        \FOR{$\ntop = \tcb{N^{\texttt{fit}}},\tcb{N^{\texttt{fit}}}-1,\tcb{N^{\texttt{fit}}}-2,\hdots,2,1$}
            \STATE \tcb{If $\vtilde_i(0)=\cdots=\vtilde_i(\ntop)$, \textbf{break}.}
            \STATE \textit{Step 8}: Perform linear regression \tcb{with a fitted intercept} on $\vtilde_i(n)$ against $n$, for $n\in \tcb{\{}0,1,\hdots,\ntop\tcb{\}}$. Record its $r^2$ and slope ($m$).
            \IF{$r^2\geq \rthr$}
                \STATE $\lbarotti=\tcb{\min\{m,0\}}$
                \STATE \textbf{break}
            \ENDIF
        \ENDFOR
    \ENDFOR
    \RETURN $\lbarotti, \forall t\in [s_i,s_i+\delot]$.
\end{algorithmic}
\end{algorithm}

\section{Numerical Results}
\label{sec:results}
The performance of the proposed extraboard assignment method is tested \tcb{in} two \tcb{case studies}. The first is based on published bus schedules of a Canadian public transit agency. The second \tcb{uses} a \tcb{synthetic} schedule mirroring important characteristics of the first, designed to \tcb{evaluate performance} at a larger scale. \tcb{The} following subsections \tcb{describe the} case \tcb{studies}, \tcb{experimental} settings, comparison benchmarks\tcb{, and results}.

\subsection{MiWay case study}
\begin{figure}[!ht]
    \centering
    \includegraphics[width=0.75\textwidth]{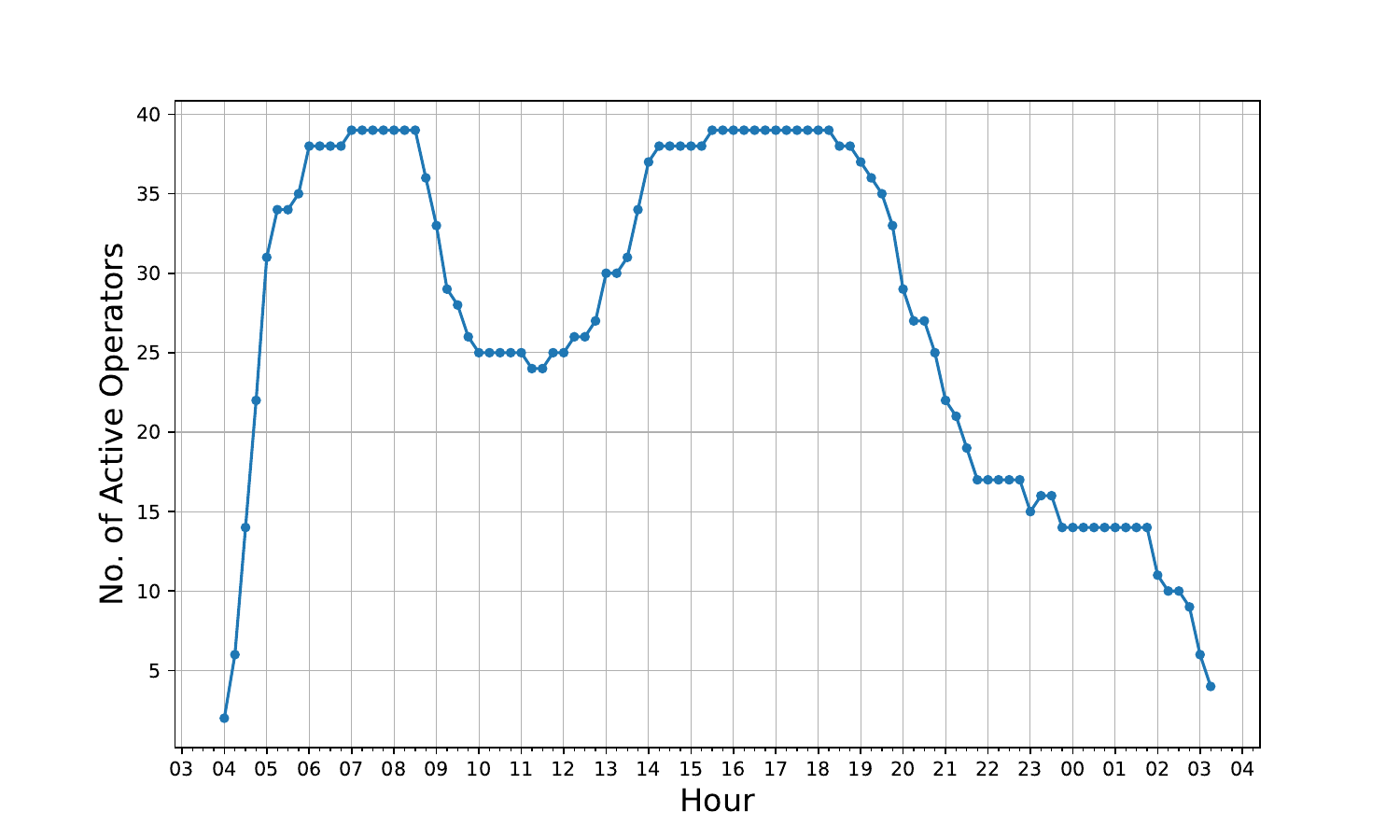}
    \caption{Daily Temporal Profile of Active Operators \tcb{in the} MiWay \tcb{Case Study}}
    \label{fig:temp_profile}
\end{figure}

The first case study \tcb{considers} a bus garage \tcb{serving} seven routes of MiWay, the transit agency serving the city of Mississauga, Ontario, Canada, located just west of Toronto. The spring 2024 weekday schedules of these routes are obtained from MiWay's \tcb{publicly available} General Transit Feed Specification (GTFS) \tcb{data} \citep{mississauga_developer_2020}. \tcb{Since these data do not identify operator duties, we construct a regular-driver schedule to approximate the daily service profiles of the selected routes.} \tcb{The modeled} day of operation starts at 4:00 AM and ends at 3:30 AM \tcb{on} the next day, \tcb{defining a} decision-making horizon \tcb{of} 94 fifteen-minute periods. \tcb{The constructed} schedule \tcb{includes} 87 regular drivers, each \tcb{operating} on one of the seven routes. \tcb{Of these,} 68 have straight shifts \tcb{comprising 7.5 productive hours, while} the other 19 have split shifts \tcb{comprising} two disjoint work sessions \tcb{totaling 7 productive hours}. Split shifts \tcb{accommodate} increased \tcb{service during} AM and PM peak hours. Figure~\ref{fig:temp_profile} displays the daily temporal profile of active operators \tcb{in the constructed schedule, with clearly visible} AM and PM peaks. In total, \tcb{the schedule contains} 643 \tcb{productive} driver-hours. Each driver run, or each \tcb{session of a} split shift, is divided into one or more consecutive work pieces. \tcb{Each} work piece \tcb{represents a} round trip on \tcb{its route, with} driver changeover \tcb{assumed to occur at} the terminal closest to the garage. \tcb{Work-piece durations range} from 1 to 2.75 hours, depending on \tcb{the} route and \tcb{time of day}.

In this \tcb{case study, we consider} two types of work sources: (1) scheduled \tcb{work} of \tcb{regular-duty drivers and} (2) \tcb{extra trips required in response to} unexpected events, such as bus \tcb{breakdowns} or \tcb{surges} in passenger demand. \tcb{For type (1), each straight run and each session of a split shift constitutes a separate work source, giving 106 sources in total. A source becomes open when its scheduled driver is unavailable for that run or session, making all its component pieces open. We set} the base \tcb{openness} probability \tcb{to} 7\% \tcb{for each type (1) source. This case-study value is comparable in magnitude to the Toronto Transit Commission's (TTC's) reported system-wide employee absenteeism rate of 7.14\% in June 2017} \citep{ttc2017ceo}. Type (2) \tcb{sources} are defined \tcb{for routes 5, 7, 11, and 51, with potential extra trips starting at fifteen-minute intervals within specified time windows}. Each \tcb{of these 334 sources} consists of a single \tcb{round-trip} work piece, \tcb{whose duration is specified by route and start time}. \tcb{Their} base openness \tcb{probabilities are case-study inputs that vary by route and time of day, ranging from 1\% to 12\%, to represent variation in the likelihood of additional work. Consistent with the source-based uncertainty model, outcomes are sampled independently across work sources, including separate sessions of split shifts. For both source types, we multiply the base probabilities by 1.5 and 2 to obtain the medium and high probability settings, respectively.}

\tcb{Table~\ref{tab:case_study_data_input} summarizes the sources and modeling basis of the key inputs used in the MiWay case study. MiWay's publicly available GTFS schedules provide the service information used to construct the instance. Agency-specific operator absence, open-work, trip-level ridership, and dispatch records were not available for this study. The corresponding inputs therefore combine literature-informed choices with values specified for the case study, as detailed in the table. The numerical experiments examine different openness probability levels, XB workforce sizes, OT availability and costs, and work-piece reward definitions. In practical implementation, a transit agency can use its own operational records to estimate these inputs.}

\begin{table}[ht!]
\centering
    \caption{\tcb{Sources and modeling basis for key inputs of the MiWay case study}}
    \label{tab:case_study_data_input}
    {
    \renewcommand{\arraystretch}{1.2}
    \scalebox{0.9}{%
    \begin{tabular}{L{3.1cm} L{6.4cm} L{6.7cm}}
    \toprule
    Input element & Case-study specification & Source or rationale \\
    \midrule
    Bus schedules & Seven MiWay routes; spring 2024 weekday schedules & Published MiWay GTFS data \citep{mississauga_developer_2020} \\
    Regular-driver schedule & 87 constructed duties (68 straight, 19 split) & Constructed to approximate daily service profiles informed by GTFS \\
    Work piece construction & Regular runs/sessions divided into route round trips; each extra trip is one piece & Round-trip pieces informed by route patterns; changeover modeled at the terminal closest to the garage \\
    $p_h$ for regular runs/sessions & Base 7\%; medium: $1.5\times$ base; high: $2\times$ base & Comparable in magnitude to the TTC-wide employee absence rate of 7.14\% in June 2017 \citep{ttc2017ceo} \\
    $p_h$ for extra trips & Routes 5, 7, 11, 51; base 1\%--12\% by route/time; medium: $1.5\times$ base; high: $2\times$ base & Route/time scenarios represent varying likelihoods of extra work from disruptions or demand surges \\
    Source outcomes & Independent Bernoulli outcomes by source (split sessions separate); all pieces revealed together at source start & Source-based uncertainty model in Section~\ref{sec:exog} \\
    XB report times & Fixed rosters from a sequential coverage heuristic using a four-route regular-work profile & Report-time scheduling approach motivated by \citet{kaysi1990scheduling} \\
    OT roster and report times & Fixed roster; report times staggered across AM/PM peaks; one-hour minimum-pay duration & Peak-coverage design; volunteer-list and minimum-pay structure motivated by \citep{deannuntis2012best} \\
    Passenger-weight inputs & Route/time ridership and headway (min); passenger-hour weight = ridership $\times$ headway$/60$ (see Section~\ref{sec:passenger_weights}) & Constructed passenger-exposure proxy; historical 2008 daily-boardings rankings provide context \citep{mississauga2009ridership} \\
    Practice-oriented piece division & Straight runs: two approximately equal pieces; split sessions: one piece each & Larger assignment blocks are motivated by transit scheduling practice \citep{boyle2009controlling} \\
    ``Overrule'' probability & 20\% in the practice-oriented benchmark & Operator refusal or supervisor discretion documented by \citep{deannuntis2012best} \\
    \bottomrule
    \end{tabular}%
    }
    }
\end{table}

\subsubsection{Experiment setup}
We conduct numerical experiments to test the performance of the proposed approximate policy under various settings for the following important parameters:
\begin{itemize}\setItemSep{-0.3em}
    \item Level of openness probabilities ($p_h$)
    \item \tcb{Number} of XBs and \tcb{availability of} OTs
    \item \tcb{Per-period overtime cost coefficient ($\costot$)}
    \item Reward \tcb{coefficients} for work pieces \tcb{($c_k$)}
\end{itemize}
We compare the performance of the approximate policy against several benchmarks: \tcb{a perfect-information} (PI) \tcb{benchmark}, a scenario-based rolling-horizon \tcb{(RH)} policy, a myopic policy, and a practice-oriented method.

\paragraph*{\tcb{Perfect-information (PI) benchmark}}
In the PI setting, all open work pieces during the workday are \tcb{known} at the beginning of the day, \tcb{allowing assignments to be optimized with complete knowledge of realized demand}. \tcb{The PI problem is a full-day version of the assignment IP} \eqref{eqn:obj}--\eqref{eqn:cons26}, \tcb{considering} all XBs $j\in\mathcal{J}$, all OTs $i\in\mathcal{I}$, and all \tcb{realized} open pieces $k\in\tcb{\bigcup_{t\in\mathcal{T}}\mathcal{W}_t}$. \tcb{Candidate XB sequences comprise non-overlapping eligible pieces throughout each XB's shift, while OT eligibility and the restriction to at most one assignment remain unchanged. Each piece may be covered at most once. The objective maximizes coverage rewards net of overtime costs, without future-value approximation terms. For each realization of open work, the PI optimal objective value provides an} upper bound \tcb{on the reward attainable by any online policy using the same roster and assignment restrictions}.

\paragraph*{\tcb{Rolling-horizon (RH) policy}}
\tcb{The scenario-based RH policy accounts for uncertain future work by solving a two-stage stochastic program in the XB-only setting. At a decision epoch $t\in\mathcal{T}$ with an available XB and an open piece starting at $t$, the first stage selects a plan for the currently available XBs $j\in\jtn$ from known open work in $\mathcal{D}_t$. Let $\mathcal{P}_t$ denote the set of feasible first-stage plans, represented by the sequence-selection and coverage variables $\ytjm$ and $\ztk$ defined above, collected in vectors $\bm{y}^t$ and $\bm{z}^t$. Each XB may be selected for at most one sequence in $\mtj$, and each piece may be covered at most once.}

\tcb{We generate scenarios $\Omega_t$ by independently sampling the openness outcomes of work sources starting in $[t+1,t+L]$, using their probabilities $p_h$. All pieces of each sampled open source enter the scenario, including pieces extending beyond the source-revelation window. Under each scenario, the same currently available XBs may be assigned additional eligible sequences of sampled work. Each XB may receive at most one such sequence, starting after its selected first-stage sequence ends, or from $t+1$ if no first-stage sequence is selected. Each sampled piece may be covered at most once within a scenario. The first-stage plan is common across scenarios, while these additional assignments are scenario dependent. The RH objective is
\begin{equation}
    \max_{(\bm{y}^t,\bm{z}^t)\in\mathcal{P}_t}
    \sum_{k\in\mathcal{D}_t} c_k\ztk
    + \frac{1}{|\Omega_t|} \sum_{\omega\in\Omega_t} Q_t(\bm{y}^t,\omega), \label{eqn:rh_obj}
\end{equation}
where $Q_t(\bm{y}^t,\omega)$ is the maximum additional reward under scenario $\omega$ subject to the second-stage assignment rules above. Thus, the model combines planned coverage rewards from known work with average scenario rewards, without a terminal value approximation.}

\tcb{Execution follows the same convention as the approximate policy: only selected first-stage sequences starting at $t$ are committed, in full. Future-starting selections remain provisional, and the other currently available XBs hold until the next period. The state is updated using the committed assignments and newly revealed work, and subsequent plans are recomputed using fresh scenarios. We test $L\in\{2,4,6\}$ periods, with 50 scenarios per solve. Each IP is solved with a 1\% relative optimality-gap target.}

\paragraph*{\tcb{Myopic policy}}
The myopic policy serves as a \tcb{rule-based benchmark that assigns open work without estimating future values}. \tcb{It} follows a greedy \tcb{priority} rule\tcb{: at each period} $t\in\mathcal{T}$, open pieces in $\dtn$\tcb{, which start at $t$, are processed sequentially in random order}. \tcb{When} the reserve \tcb{roster consists only} of XBs, each piece is assigned to an eligible \tcb{available} XB with the earliest shift end time\tcb{, if one exists}. When \tcb{OTs are also available}, we test two \tcb{alternatives}: ``XB-first'' and ``OT-first''. \tcb{Under} XB-first, \tcb{each piece is assigned to an eligible available XB whenever possible; otherwise, it is assigned to an available OT}. \tcb{Under} OT-first, \tcb{each piece is assigned to an available OT whenever possible; otherwise, an eligible available XB is selected}. \tcb{Both alternatives prioritize XBs} with the earliest shift end time \tcb{and OTs who have} been on standby the longest within their respective groups. \tcb{Each assignment covers a single piece starting at the current period. Pieces that cannot be assigned are dropped, while available operators who receive no assignment hold according to the state transitions described earlier.}

\paragraph*{\tcb{Practice-oriented method}}
The practice-oriented method \tcb{represents larger work assignments and possible departures from the XB-first priority rule}. \tcb{In this case study, representing individual round trips as work pieces can require several separate assignments within a regular operator run. To represent assignment in larger blocks, consecutive original pieces of each straight run are grouped into} two work pieces, \tcb{while each session of a split run} becomes a single work piece. \tcb{Extra-trip pieces remain unchanged. The resulting pieces are processed using the} XB-first myopic \tcb{rule described} above\tcb{, including random ordering of pieces starting in the same period. We additionally allow the initial operator choice to be overruled, representing, for example, an operator declining a particular assignment or a supervisor choosing a different operator}. In all \tcb{experiments with this benchmark}, an overrule \tcb{occurs with} probability 20\%. When \tcb{it} occurs, another \tcb{eligible} available XB or OT \tcb{is selected uniformly at random, excluding the initially selected operator. This replacement is assigned without another overrule draw; if no eligible alternative exists, the piece is dropped.}

\tcb{For performance evaluation, we generate} 1000 sample paths \tcb{for each openness-probability setting by independently sampling each work source's opening outcome according to its probability} ($p_h$). \tcb{Within each setting, the methods included in a performance comparison are evaluated on the same sample paths. Performance is measured by} the total \tcb{realized} reward \tcb{over the workday, with means and standard deviations calculated across the evaluation paths. The RH reward comparison uses the 1000 paths for the base-probability, XB-only setting, whereas its computational-time assessment under high openness probabilities and 18 XBs uses the first 20 paths. These evaluation paths are distinct from the 50 planning scenarios generated within each RH solve. We also examine computational times for the approximate and RH policies, with the measurement scope specified in their respective computational assessments.} \tcb{The reported computational-time results are obtained} on a PC \tcb{equipped with a} 2.80-GHz Intel Core i7-1165G7 \tcb{processor} and \tcb{16 GB of} RAM. \tcb{The optimization} models for the approximate policy\tcb{,} PI\tcb{, and RH} are solved \tcb{using} Gurobi version 12. \tcb{The additional component comparisons and complete synthetic PI calculations use Gurobi 12.0.3 with one solver thread on an ARM-based macOS system; their runtimes are not included in the reported timing comparisons.}

\subsubsection{Ordinary test setting: XBs only, simple reward coefficients}
The main set of numerical experiments \tcb{is} conducted with only XBs considered in decision-making and no OTs involved and with the reward for covering a work piece, $c_k$, set equal to its duration measured in hours. We closely examine the decisions involving only XBs because it is a useful setting if the garage is more concerned with the efficiency of XB utilization. \tcb{This focus is also relevant when sufficient OTs are available to cover any work left uncovered by XBs; such OT assignments are not modeled in this setting.} The weights of work pieces are chosen so that the total reward may be interpreted as the total number \tcb{of} hours of open work covered by XBs, or equivalently, the total productive hours \tcb{worked} by the XB roster. Recall that at the start of a given workday, the number of XBs and their report times are known from the output by upper decision levels. To account for the time loss due to administrative procedures at the bus garage, \tcb{we assume that} each XB works a straight shift of 7.5 hours (equivalent to 30 periods). We test various XB \tcb{roster} sizes $\noxb \in \{6,8,10,12,14,16,18\}$, with XBs' report times \tcb{held fixed across the compared policies and} generated using \tcb{a sequential coverage heuristic motivated by} \citet{kaysi1990scheduling}.

\subsubsection{\tcb{Threshold sensitivity and component assessment}}
\label{sec:component_assessment}
We examine the sensitivity of the approximate policy to the regression-fit threshold $\rthr$ in Algorithm~\ref{alg:loss}, which determines the range of positions used to estimate the XB loss parameters. For these fitting calculations, we use 1,000 simulated work-source opening scenarios for each XB/time pair requiring estimation and set $N^{\texttt{fit}}=9$. The resulting average values at positions $n=0,\ldots,9$ are reused across thresholds. Tables \ref{tab:r2_sensi_base} and \ref{tab:r2_sensi_high} report the average open work hours covered under the base and high openness probabilities, respectively.

\begin{table}[ht!]
    \centering
    \caption{\tcb{Open work hours covered by the approximate policy with varying $\rthr$ (XBs only; base $p_h$ levels)}}
    \label{tab:r2_sensi_base}
    \begin{tabular}{| c || c | c | c | c | c |}
    \hline
    $\noxb$ & \multicolumn{5}{c|}{$\rthr$ value} \\ \hline
    & 0.999 & 0.995 & 0.990 & 0.970 & 0.950 \\ \hline
    6 & 34.88 & 34.88 & 34.88 & 34.88 & 34.89 \\
    8 & 44.28 & 44.28 & 44.28 & 44.28 & 44.28 \\
    10 & 51.42 & 51.42 & 51.42 & 51.42 & 51.43 \\
    12 & 58.59 & 58.58 & 58.58 & 58.59 & 58.58 \\
    14 & 63.54 & 63.53 & 63.55 & 63.56 & 63.55 \\
    16 & 66.52 & 66.51 & 66.54 & 66.54 & 66.53 \\
    18 & 69.55 & 69.57 & 69.57 & 69.60 & 69.61 \\
    \hline
    \end{tabular}
\end{table}

\begin{table}[ht!]
    \centering
    \caption{\tcb{Open work hours covered by the approximate policy with varying $\rthr$ (XBs only; high $p_h$ levels)}}
    \label{tab:r2_sensi_high}
    \begin{tabular}{| c || c | c | c | c | c |}
    \hline
    $\noxb$ & \multicolumn{5}{c|}{$\rthr$ value} \\ \hline
    & 0.999 & 0.995 & 0.990 & 0.970 & 0.950 \\ \hline
    6 & 40.91 & 40.91 & 40.91 & 40.91 & 40.91 \\
    8 & 53.62 & 53.62 & 53.62 & 53.62 & 53.62 \\
    10 & 65.19 & 65.19 & 65.19 & 65.19 & 65.19 \\
    12 & 77.34 & 77.34 & 77.34 & 77.33 & 77.34 \\
    14 & 87.83 & 87.82 & 87.82 & 87.82 & 87.82 \\
    16 & 95.62 & 95.62 & 95.62 & 95.63 & 95.62 \\
    18 & 105.27 & 105.28 & 105.28 & 105.28 & 105.28 \\
    \hline
    \end{tabular}
\end{table}

Average open work covered is largely insensitive to $\rthr$ over the tested range, with a maximum difference of approximately 0.06 hours for a fixed roster size and openness-probability setting. The threshold affects the fitted holding-value adjustments, while assignment rewards and individual future-value estimates are computed independently of it. In this XB-only setting, these adjustments can reduce a provisional plan's value only when it leaves multiple available XBs with sufficient remaining shift time without selected sequences. The approximate-policy means in the subsequent ordinary XB-only base- and medium-probability tables correspond to $\rthr=0.970$; the high-probability table's approximate-policy means correspond to $\rthr=0.995$.

To assess the two value components separately, Table~\ref{tab:component_diagnostic} reports matched diagnostics under base openness probabilities for XB-only ordinary and weighted rewards and for the case with ten OTs at cost $0.1$. Within each setting, all three variants use the same 1,000 work-opening paths and feasible assignment decisions. The zero-continuation variant retains the online IP but removes its future-value terms; the individual-value variant restores those terms, and the full variant also includes the fitted competition adjustment. Each setting uses its matching $\rthr=0.970$ coefficient set.

Individual values increase mean reward in all six comparisons. In the two ordinary XB-only cases and the weighted case with six XBs, the adjustment changes mean reward only slightly. With 18 XBs under weighted rewards, it increases mean reward by 7.92 proxy passenger-hours. With ten OTs, it increases mean net reward by 0.13 and 0.69 hours for six and 14 XBs, respectively. To examine the XB-only pattern across intermediate roster sizes, Appendix~\ref{sec:app_component_diagnostics} extends the paired adjustment comparison to all seven tested rosters, holding the paths and fitted coefficients fixed. Ordinary-reward differences remain small throughout, whereas weighted-reward gains increase from 0.39 to 7.92 proxy passenger-hours across the rosters with eight to 18 XBs. The appendix reports all paired differences and their approximate 95\% intervals.

\begin{table}[ht!]
\centering
\caption{\tcb{Matched component diagnostics: mean reward on 1,000 common paths (base $p_h$)}}
\label{tab:component_diagnostic}
{\small\setlength{\tabcolsep}{4pt}
\begin{tabular}{@{}p{0.29\textwidth}rrrr@{}}
\hline
Setting & $\noxb$ & Zero continuation & Individual values & With adjustment \\
\hline
Ordinary, XBs only (h) & 6 & 34.09425 & 34.88225 & 34.88475 \\
& 18 & 67.05000 & 69.63475 & 69.59975 \\
Weighted, XBs only (proxy passenger-h) & 6 & 395.83833 & 405.40433 & 405.36400 \\
& 18 & 777.12183 & 795.80842 & 803.72783 \\
Ten OTs, cost $0.1$ (net h) & 6 & 44.68775 & 45.05975 & 45.18560 \\
& 14 & 66.46760 & 66.91190 & 67.60000 \\
\hline
\end{tabular}
}
\end{table}

\subsubsection{\tcb{Computational efficiency of the approximate policy}}
To \tcb{assess the computational requirements of} the approximate policy, we report \tcb{off-line parameter-computation times and online timing summaries for the} integer program (\ref{eqn:obj3})-(\ref{eqn:cons316}). These durations are reported for the high $p_h$ case\tcb{, which has the largest expected amount of open work among the tested probability settings}.

Off-line \tcb{parameter computation includes obtaining individual-value estimates} (Algorithms \ref{alg:bdp} and \ref{alg:bdpot} for XBs and OTs, respectively) and loss parameters (Algorithms \ref{alg:loss} and \ref{alg:loss_ot} for XBs and OTs, respectively). \tcb{The reported workload considers} XB report times from period 0 to 64 and OT start times from period 0 to 72. The \tcb{reported} run \tcb{times} for each of the four algorithms are \tcb{given} in Table \ref{tab:alg_runtime}.

\begin{table}[ht!]
    \centering
    \caption{\tcb{Reported off-line parameter-computation times} (high $p_h$ levels)}
    \label{tab:alg_runtime}
    \resizebox{0.375\textwidth}{!}{
    \begin{tabular}{cc}
    \hline
    Algorithm & Run Time (s) \\ \hline
    XB value functions (Alg. \ref{alg:bdp}) & 60.7 \\
    OT value functions (Alg. \ref{alg:bdpot}) & 1.3  \\
    XB loss parameters (Alg. \ref{alg:loss}) & 2172.8 \\
    OT loss parameters (Alg. \ref{alg:loss_ot}) & 285.5 \\ \hline
    \end{tabular}
    }
\end{table}

\tcb{The reported XB computations take longer than the OT computations; XBs have longer shifts, so their value estimates and loss parameters cover more potential future assignments.} \tcb{The loss-parameter algorithms are simulation-based, whereas the individual-value estimates are computed by backward recursion.} Altogether, these four algorithms take around 42 minutes to run. \tcb{Their parameters can be reused while the relevant inputs remain unchanged; the reported total suggests that daily recomputation is possible if needed.}

In addition, the integer program of the approximate policy has to be solved quickly enough to support making assignment decisions every 15-minute period. \tcb{Across the reported 1,000 sample paths, we examine online times for} $\noxb=18$. 
The median \tcb{reported online} time is 0: in more than half of the \tcb{sampled} decision epochs, the integer program does not need to be run because \tcb{no operator is available or no open piece starts at that epoch}. The maximum \tcb{reported online} time is less than 5 seconds\tcb{, and 99\% of the reported online times are} less than 0.1 seconds. These \tcb{reported times suggest that online optimization can fit within the 15-minute decision interval in this tested setting}.

\subsubsection{\tcb{Performance of the rolling-horizon method}}
We evaluate the reward and decision-time performance of the RH benchmark. The RH method does not require the off-line value and loss estimates used by the approximate policy. For its computational assessment, we record elapsed time for complete decision epochs, including scenario generation, model construction, optimization, and state updates, using the first 20 sample paths with $\noxb=18$ and high $p_h$ levels. For $L=6$, the average recorded epoch time is 3.41 seconds. The maximum epoch times are approximately 81, 112, and 217 seconds for $L=2,4,6$, respectively. These maxima are below the 15-minute interval between decisions, indicating computational feasibility for this RH implementation in the tested setting.

We also compare the performance of RH with that of the approximate and myopic policies in terms of average reward gained (i.e., hours of open work covered) across 1,000 common sample paths. The base $p_h$ level setting is used in this analysis. Table \ref{tab:xb_hours_w_rh} shows the results for all $\noxb$ tested. 

\begin{table}[ht!]
    \centering
    \caption{\tcb{Open work hours covered with RH (XBs only; base $p_h$ levels)}}
    \label{tab:xb_hours_w_rh}
    \resizebox{0.75\textwidth}{!}{
    \begin{tabular}{| c | c | c | c | c | c |}
    \hline
    $\noxb$ & Approx. policy & Myopic policy & RH ($L = 2$) & RH ($L = 4$) & RH ($L = 6$) \\ \hline
    6 & 34.88 & 34.25 & 34.23 & 33.93 & 34.29 \\
    8 & 44.28 & 43.50 & 43.21 & 42.84 & 43.24 \\
    10 & 51.42 & 50.45 & 50.32 & 49.85 & 50.30 \\
    12 & 58.59 & 57.41 & 57.34 & 56.91 & 57.41 \\
    14 & 63.56 & 62.21 & 62.29 & 61.87 & 62.38 \\
    16 & 66.54 & 65.35 & 65.17 & 64.86 & 65.30 \\
    18 & 69.60 & 68.98 & 68.42 & 68.20 & 68.57 \\
    \hline
    \end{tabular}
    }
\end{table}

In these tests, larger $L$ does not necessarily lead to better performance. The $L=4$ results consistently perform worse than the other two RH candidates. $L=6$ only marginally outperforms $L=2$ in most cases. At $\noxb=10$, $L=2$ has the highest reported mean among RH methods. Overall, the highest RH mean at each roster size is close to the myopic mean and below the approximate-policy mean. Together, the reward and timing assessments show the approximate policy's value against a computationally feasible stochastic-lookahead benchmark in the tested settings. The subsequent comparisons focus on the other methods across additional probability, operator, and reward settings; Table~\ref{tab:xb_hours_w_rh} provides the RH comparison for the tested XB-only base-probability setting.

\subsubsection{Performance comparisons for the ordinary test setting}
\tcb{We next compare} PI, the approximate policy, the myopic policy, and the practice-oriented method under the ordinary test setting (XBs only, simple reward coefficients) at base $p_h$ levels, where the average total amount of open work is 76.86 hours \tcb{across the 1,000 evaluation paths}. Table \ref{tab:xb_hours} summarizes the average total reward and the corresponding SD for each method at each $\noxb$ level. For the approximate policy, the myopic policy, and the practice-oriented method, we also report the gap between their respective total reward and that of the PI solution. In terms of averages, the performance of the practice-oriented method lags \tcb{substantially} behind even the myopic policy for all $\noxb$ values examined. \tcb{This benchmark combines larger assignment blocks with possible overrides of the XB-first priority rule.} Among the three methods outside the practice-oriented method, where each work piece corresponds to a single round trip, the approximate policy consistently \tcb{has a higher mean reward than} the myopic policy while being behind the PI solution. The gaps between the myopic policy \tcb{and} both PI and approximate policy appear modest, especially at larger $\noxb$. \tcb{One possible contributor is greater aggregate XB capacity relative to expected open work} (e.g., at $\noxb=18$, 135 hours of XB shift \tcb{capacity} vs. 76.86 hours of average open work)\tcb{, although work timing and eligibility still restrict coverage}. When there is an abundant supply of XBs, there is more tolerance for making \tcb{less productive} assignment decisions. Other potential factors for the modest gaps are \tcb{the lower expected amount of open work} at base $p_h$ levels, homogeneity of the reserve operator roster (with XBs only), and all work pieces having the same reward per time period. In later test settings, we explore how changing these parameters can affect model performances.

Standard deviations of \tcb{covered} open work hours increase along with the mean as $\noxb$ increases. However, the coefficient of variation (CV), defined as the ratio between SD and mean, also increases with more XBs added. \tcb{Thus, across the evaluated paths, relative day-to-day variation in covered work is higher for the larger rosters tested.} \tcb{At a given $\noxb$, the methods' CVs are generally similar, with somewhat larger differences for smaller rosters.}

\setlength{\tabcolsep}{4pt} 
\renewcommand{\arraystretch}{1.15} 
\begin{table}[ht!]
    \centering
    \caption{Open work hours covered (XBs only; base $p_h$ levels)}
    \label{tab:xb_hours}
    \resizebox{0.75\textwidth}{!}{
    \begin{tabular}{| c | c | c | c | c |}
    \hline
    $N_\texttt{XB}$ &Perfect information &Approx. policy &Myopic policy & Pract. method \\
    \hline
    6  & 36.33 (SD = 3.97) & 34.88 (SD = 4.16)  & 34.25 (SD = 3.98)  & 28.83 (SD = 4.10) \\
    && {[}gap = 3.98\%{]}& {[}gap = 5.73\%{]}& {[}gap = 20.63\%{]} \\ \hline
    8  & 46.32 (SD = 5.76) & 44.28 (SD = 5.79)  & 43.50 (SD = 5.60)  & 37.41 (SD = 5.07) \\
    && {[}gap = 4.42\%{]}& {[}gap = 6.09\%{]}& {[}gap = 19.25\%{]} \\ \hline
    10 & 53.60 (SD = 7.96) & 51.42 (SD = 7.82) & 50.45 (SD = \tcb{7.56}) & 43.97 (SD = 6.53) \\
    && {[}gap = 4.08\%{]}& {[}gap = 5.88\%{]}& {[}gap = 17.96\%{]} \\ \hline
    12 & 60.61 (SD = 10.10) & 58.59 (SD = 9.71) & 57.41 (SD = 9.36) & 50.90 (SD = 7.97) \\
    && {[}gap = 3.34\%{]}& {[}gap = 5.28\%{]}& {[}gap = 16.03\%{]} \\ \hline
    14 & 65.41 (SD = 12.05) & 63.56 (SD = 11.58) & 62.21 (SD = 11.05) & 55.97 (SD = 9.41) \\
    && {[}gap = 2.82\%{]}& {[}gap = 4.89\%{]}& {[}gap = 14.43\%{]} \\ \hline
    16 & 68.26 (SD = 13.36) & 66.54 (SD = 12.85) & 65.35 (SD = 12.30) & 59.54 (SD = 10.85) \\
    && {[}gap = 2.53\%{]}& {[}gap = 4.26\%{]}& {[}gap = 12.77\%{]} \\ \hline
    18 & 71.14 (SD = 14.99) & 69.60 (SD = 14.40) & 68.98 (SD = 13.94) & 63.31 (SD = 12.25) \\
    && {[}gap = 2.16\%{]}& {[}gap = 3.03\%{]}& {[}gap = 11.01\%{]} \\
    \hline
    \end{tabular}
    }
\end{table}

In Table \ref{tab:uow_hours}, we report the amount of uncovered open work (defined as total open work minus covered open work) as well as its percentage out of total open work. Recall that there is an average of 76.86 hours of open work per \tcb{simulated} scenario. \tcb{Dividing this average by the 7.5-hour XB shift gives approximately 10.25 XB shifts, but this aggregate calculation ignores when work occurs and which XBs are eligible.} \tcb{With a roster of 10 XBs,} there is still \tcb{a substantial} amount of uncovered open work left (33\% of total open work for \tcb{the} approximate policy; 43\% for \tcb{the} practice-oriented method) to be either assigned to OTs or \tcb{left uncovered}.

\begin{table}[ht!]
    \centering
    \caption{Uncovered open work hours and percentage out of total open work (XBs only; base $p_h$ levels)}
    \label{tab:uow_hours}
    \resizebox{0.705\textwidth}{!}{
    \begin{tabular}{| c | c | c | c | c |}
    \hline
    $N_\texttt{XB}$ &Perfect information &Approx. policy &Myopic policy& Pract. method  \\
    \hline
    6 & 40.53 (52.73\%) & 41.98 (54.61\%) & 42.61 (55.44\%) & 48.03 (62.48\%) \\
    8 & 30.54 (39.73\%) & 32.58 (42.39\%) & 33.36 (43.40\%) & 39.45 (51.33\%) \\
    10 & 23.26 (30.26\%) & 25.44 (33.10\%) & 26.41 (34.36\%) & 32.89 (42.79\%) \\
    12 & 16.25 (21.14\%) & 18.27 (23.78\%) & 19.45 (25.30\%) & 25.96 (33.78\%) \\
    14 & 11.45 (14.90\%) & 13.30 (17.30\%) & 14.65 (19.06\%) & 20.89 (27.18\%) \\
    16 & 8.60 (11.18\%) & 10.32 (13.43\%) & 11.51 (14.97\%) & 17.32 (22.53\%) \\
    18 & 5.72 (7.45\%) & 7.26 (9.45\%) & 7.88 (10.25\%) & 13.55 (17.64\%) \\
    \hline
    \end{tabular}
    }
\end{table}

We examine the utilization rate of XBs. The utilization rate is defined as the percentage of total productive hours (i.e., total reward) out of the total shift duration of the XB roster, calculated as $7.5\times \noxb$. Table \ref{tab:xb_utilize} displays the average utilization rates under the four methods. \tcb{Average utilization decreases as the roster grows while the open-work scenarios are held fixed.} \tcb{Even for smaller rosters, some XB shift time remains idle while open work is uncovered (Table~\ref{tab:uow_hours}); aggregate capacity alone does not ensure that an eligible XB is available when each piece starts.}

\begin{table}[ht!]
    \centering
    \caption{Average XB utilization rate (XBs only; base $p_h$ levels)}
    \label{tab:xb_utilize}
    \resizebox{0.65\textwidth}{!}{
    \begin{tabular}{| c | c | c | c | c |}
    \hline
    $N_\texttt{XB}$ &Perfect information &Approx. policy &Myopic policy& Pract. method  \\
    \hline
    6 & 80.73\% & 77.52\% & 76.10\% & 64.08\% \\
    8 & 77.21\% & 73.80\% & 72.50\% & 62.35\% \\
    10 & 71.47\% & 68.56\% & 67.27\% & 58.63\% \\
    12 & 67.35\% & 65.10\% & 63.79\% & 56.55\% \\
    14 & 62.29\% & 60.53\% & 59.25\% & 53.31\% \\
    16 & 56.89\% & 55.45\% & 54.46\% & 49.62\% \\
    18 & 52.69\% & 51.56\% & 51.10\% & 46.89\% \\
    \hline
    \end{tabular}
    }
\end{table}

Holding other parameters constant, we examine two increased $p_h$ levels: ``medium'' $p_h$ level is 1.5 times the base $p_h$, and ``high'' level is 2 times the base $p_h$. Tables \ref{tab:xb_hours_med} and \ref{tab:xb_hours_high} report the model performances at medium and high $p_h$ levels. In terms of average performance, similar patterns of relative gaps are observed among the four methods when compared against the base $p_h$ case from Table \ref{tab:xb_hours}. For the practice-oriented method, the gap from PI becomes \tcb{less variable} across different $\noxb$ values, hovering around 16\% for most instances. The approximate policy again \tcb{has a higher mean reward than} the myopic policy for all $\noxb$ values. The performance gap between the approximate policy solution and the PI solution is generally in the 3\% -- 4\% range, similar to the gaps observed in the base $p_h$. On the other hand, solutions for the myopic policy are observed to have larger gaps. \tcb{Thus, the approximate policy's average advantage over the myopic policy is larger in these higher-openness settings, which also have higher expected open-work demand.}

\begin{table}[ht!]
    \centering
    \caption{Open work hours covered (XBs only; medium $p_h$ levels)}
    \label{tab:xb_hours_med}
    \resizebox{0.75\textwidth}{!}{
    \begin{tabular}{| c | c | c | c | c |}
    \hline
    $N_\texttt{XB}$ &Perfect information &Approx. policy &Myopic policy &Pract. method \\
    \hline
    6  & 40.16 (SD = 2.51) & 38.86 (SD = 2.87)  & 37.78 (SD = 2.73)  & 32.73 (SD = 3.28) \\
    && {[}gap = 3.24\%{]}& {[}gap = 5.93\%{]}& {[}gap = 18.50\%{]} \\ \hline
    8  & 52.52 (SD = 3.51) & 50.55 (SD = 3.84) & 49.15 (SD = 3.63) & 42.86 (SD = 3.97) \\
    && {[}gap = 3.76\%{]}& {[}gap = 6.43\%{]}& {[}gap = 18.39\%{]} \\ \hline
    10 & 62.87 (SD = 5.44) & 60.54 (SD = 5.58) & 58.88 (SD = 5.24) & 51.87 (SD = 5.16) \\
    && {[}gap = 3.71\%{]}& {[}gap = 6.34\%{]}& {[}gap = 17.50\%{]} \\ \hline
    12 & 73.61 (SD = 6.96) & 70.93 (SD = 6.99)  & 69.04 (SD = 6.61) & 61.31 (SD = 6.12) \\
    && {[}gap = 3.65\%{]}& {[}gap = 6.21\%{]}& {[}gap = 16.71\%{]} \\ \hline
    14 & 82.11 (SD = 8.93) & 79.17 (SD = 8.90) & 76.96 (SD = 8.37) & 68.83 (SD = 7.76) \\
    && {[}gap = 3.59\%{]}& {[}gap = 6.27\%{]}& {[}gap = 16.17\%{]} \\ \hline
    16 & 87.93 (SD = 10.69) & 84.93 (SD = 10.52) & 82.93 (SD = 10.02) & 74.92 (SD = 9.00) \\
    && {[}gap = 3.41\%{]}& {[}gap = 5.69\%{]}& {[}gap = 14.80\%{]} \\ \hline
    18 & 94.96 (SD = 12.55) & 92.03 (SD = 12.11) & 89.97 (SD = 11.57) & 81.90 (SD = 10.18) \\
    && {[}gap = 3.09\%{]}& {[}gap = 5.25\%{]}& {[}gap = 13.75\%{]} \\
    \hline
    \end{tabular}
    }
\end{table}

\begin{table}[ht!]
    \centering
    \caption{Open work hours covered (XBs only; high $p_h$ levels)}
    \label{tab:xb_hours_high}
    \resizebox{0.75\textwidth}{!}{
    \begin{tabular}{| c | c | c | c | c |}
    \hline
    $N_\texttt{XB}$ &Perfect information &Approx. policy &Myopic policy &Pract. method  \\
    \hline
    6  & 41.95 (SD = 1.57)  & 40.91 (SD = 2.00)  & 39.43 (SD = 1.98)  & 34.79 (SD = 2.86)  \\
    && {[}gap = 2.47\%{]}& {[}gap = 5.99\%{]}& {[}gap = 17.05\%{]} \\ \hline
    8  & 55.30 (SD = 2.13)  & 53.62 (SD = 2.65)  & 51.71 (SD = 2.51)  & 46.04 (SD = 3.44)  \\
    && {[}gap = 3.03\%{]}& {[}gap = 6.49\%{]}& {[}gap = 16.74\%{]} \\ \hline
    10 & 67.43 (SD = 3.73)  & 65.19 (SD = 4.19)  & 63.11 (SD = 3.88)  & 56.21 (SD = 4.50) \\
    && {[}gap = 3.33\%{]}& {[}gap = 6.42\%{]}& {[}gap = 16.65\%{]} \\ \hline
    12 & 80.00 (SD = 4.68)  & 77.34 (SD = 5.04)  & 74.72 (SD = 4.80)  & 67.06 (SD = 5.05) \\
    && {[}gap = 3.33\%{]}& {[}gap = 6.59\%{]}& {[}gap = 16.17\%{]} \\ \hline
    14 & 90.86 (SD = 6.29)  & 87.82 (SD = 6.53)  & 84.91 (SD = 6.12)  & 76.51 (SD = 6.20) \\
    && {[}gap = 3.34\%{]}& {[}gap = 6.54\%{]}& {[}gap = 15.79\%{]} \\ \hline
    16 & 98.98 (SD = 8.23)  & 95.63 (SD = 8.33)  & 92.81 (SD = 7.82)  & 84.18 (SD = 7.48) \\
    && {[}gap = 3.38\%{]}& {[}gap = 6.22\%{]}& {[}gap = 14.95\%{]} \\ \hline
    18 & 108.92 (SD = 9.50) & 105.28 (SD = 9.35) & 102.49 (SD = 8.84) & 93.31 (SD = 8.22) \\
    && {[}gap = 3.34\%{]}& {[}gap = 5.91\%{]}& {[}gap = 14.34\%{]} \\
    \hline
    \end{tabular}
    }
\end{table}

\tcb{At a fixed $\noxb$, CVs generally decrease as $p_h$ increases in these experiments.} \tcb{Differences in CV among the four methods are generally smaller at larger roster sizes.}

\tcb{The higher $p_h$ levels also yield more covered work and higher XB utilization at a fixed roster size, as more work is expected to open. To assess how much of that work remains uncovered,} \tcb{we report} \tcb{the corresponding results} in Tables \ref{tab:uow_hours_med} and \ref{tab:uow_hours_high} for medium and high $p_h$, respectively.

\begin{table}[ht!]
    \centering
    \caption{Uncovered open work hours and percentage out of total open work (XBs only; medium $p_h$ levels)}
    \label{tab:uow_hours_med}
    \resizebox{0.705\textwidth}{!}{
    \begin{tabular}{| c | c | c | c | c |}
    \hline
    $N_\texttt{XB}$ &Perfect information &Approx. policy &Myopic policy &Pract. method \\
    \hline
    6 & 75.53 (65.28\%) & 76.83 (66.41\%) & 77.91 (67.34\%) & 82.96 (71.71\%) \\
    8 & 63.17 (54.60\%) & 65.14 (56.31\%) & 66.55 (57.52\%) & 72.83 (62.95\%) \\
    10 & 52.82 (45.66\%) & 55.15 (47.67\%) & 56.81 (49.10\%) & 63.82 (55.17\%) \\
    12 & 42.08 (36.37\%) & 44.76 (38.69\%) & 46.65 (40.32\%) & 54.38 (47.00\%) \\
    14 & 33.58 (29.02\%) & 36.52 (31.56\%) & 38.73 (33.48\%) & 46.86 (40.50\%) \\
    16 & 27.76 (23.99\%) & 30.76 (26.59\%) & 32.76 (28.32\%) & 40.78 (35.25\%) \\
    18 & 20.73 (17.92\%) & 23.66 (20.45\%) & 25.72 (22.23\%) & 33.79 (29.21\%) \\
    \hline
    \end{tabular}
    }
\end{table}

\begin{table}[ht!]
    \centering
    \caption{Uncovered open work hours and percentage out of total open work (XBs only; high $p_h$ levels)}
    \label{tab:uow_hours_high}
    \resizebox{0.705\textwidth}{!}{
    \begin{tabular}{| c | c | c | c | c |}
    \hline
    $N_\texttt{XB}$ &Perfect information &Approx. policy &Myopic policy &Pract. method  \\
    \hline
    6  & 111.51 (72.67\%) & 112.55 (73.34\%) & 114.03 (74.30\%) & 118.66 (77.33\%) \\
    8  & 98.16 (63.96\%)  & 99.84 (65.06\%)  & 101.75 (66.30\%) & 107.42 (70.00\%) \\
    10 & 86.02 (56.06\%)  & 88.27 (57.52\%)  & 90.35 (58.88\%)  & 97.25 (63.37\%)  \\
    12 & 73.46 (47.87\%)  & 76.12 (49.60\%)  & 78.73 (51.31\%)  & 86.40 (56.30\%)  \\
    14 & 62.60 (40.79\%)  & 65.64 (42.77\%)  & 68.55 (44.67\%)  & 76.95 (50.14\%)  \\
    16 & 54.48 (35.50\%)  & 57.83 (37.69\%)  & 60.64 (39.52\%)  & 69.28 (45.14\%)  \\
    18 & 44.54 (29.02\%)  & 48.17 (31.39\%)  & 50.97 (33.22\%)  & 60.15 (39.20\%) \\
    \hline
    \end{tabular}
    }
\end{table}

At any given $\noxb$, \tcb{uncovered open work increases more than proportionally} with elevated $p_h$. For example, there are about 25 hours of uncovered open work for base $p_h$ with $\noxb=10$ \tcb{under} the approximate policy. This value more than doubles (55) for medium $p_h$ and more than triples (88) for high $p_h$, even though these $p_h$ levels are only 1.5 and 2 times the base $p_h$. \tcb{For the fixed rosters tested, the additional work opening at higher probability levels exceeds the additional work that XBs cover.}

\subsubsection{\tcb{Example visualization of extraboard assignment decisions}}
\begin{figure}[!ht]
    \centering
    \includegraphics[width=0.75\textwidth]{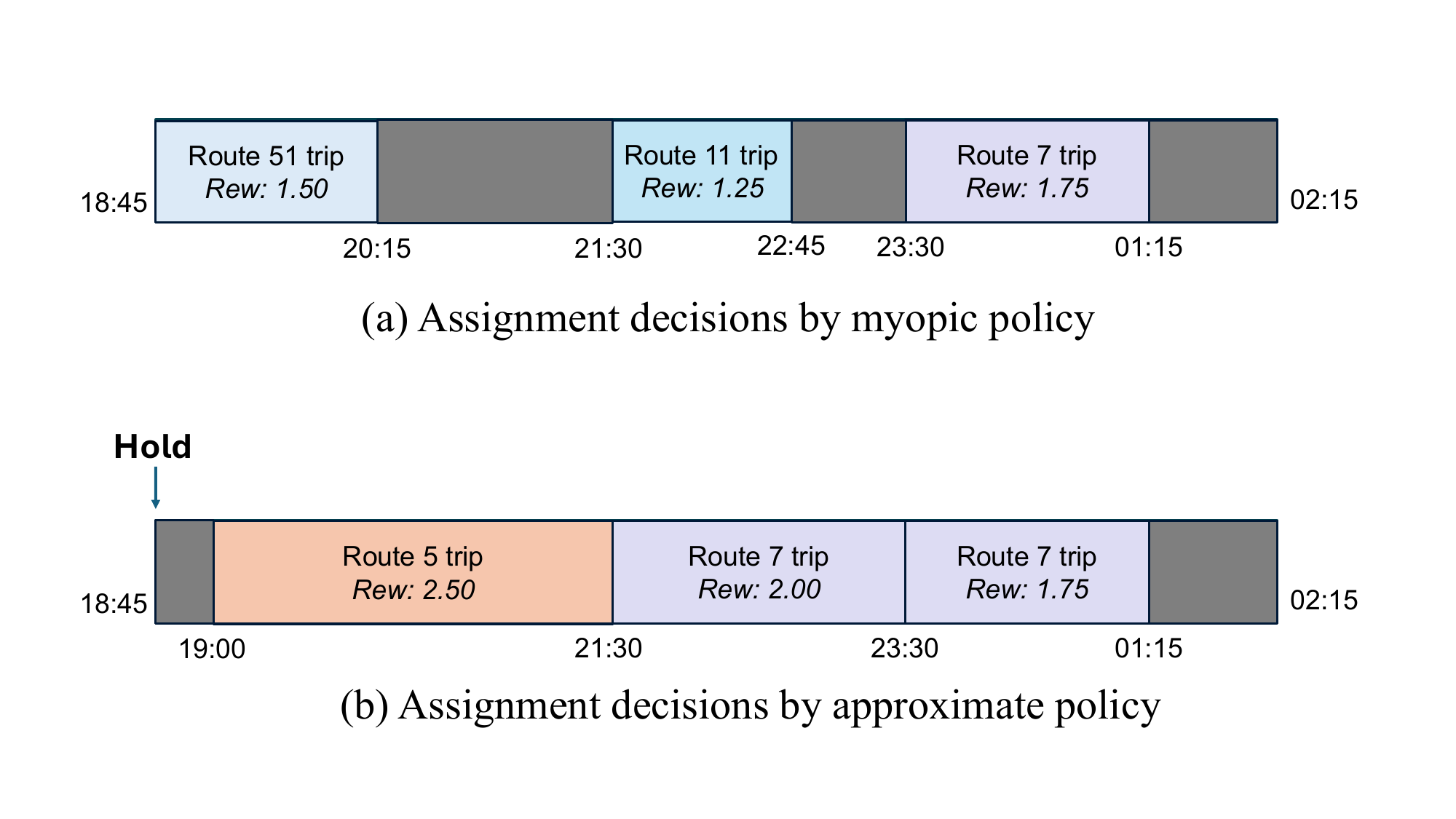}
    \caption{Example decision timelines by myopic and approximate policies}
    \label{fig:visual_assign}
\end{figure}

The numerical results presented thus far use aggregate statistics across all sample paths as primary performance metrics. To provide a detailed perspective on XB assignment decisions, Figure \ref{fig:visual_assign} compares the decisions made by the myopic and approximate policies for a single XB (hereafter referred to as the focal XB) along one illustrative sample path.

In this example, the focal XB reports at 18:45 and finishes their shift at 2:15 on the next day. The focal XB's report time is the latest and most isolated among those of the six XBs in the roster. The reward for covering a work piece, denoted by ``Rew'' in the figure, is measured by the duration of the work piece in hours. Two key differences between the assignment decisions of the myopic and approximate policies are highlighted.

The first difference occurs at 18:45, immediately after the focal XB reports for duty. The myopic policy assigns the operator to a newly opened round trip on Route 51. In contrast, the approximate policy chooses to hold the operator, forgoing the immediate reward from this relatively short work piece while its decision model also accounts for estimated future opportunities. On this realized path, the focal XB is assigned at the very next period to a newly opened work piece on Route 5 that offers a higher reward. The later opening illustrates the potential benefit of holding, although it was unknown at 18:45.

After the focal XB finishes the respective assignments under both policies, the next illustrated decision occurs at 21:30. The approximate policy selects a work sequence consisting of two round trips on Route 7, which are both known open at 21:30 and together offer more reward to this XB than the Route 11 trip. Under the myopic policy, the Route 11 trip is assigned to the focal XB. At 23:30, after completing this trip at 22:45, the focal XB is assigned to the second round trip on Route 7, which starts at 23:30. Across the illustrated workday, the focal XB covers 6.25 hours under the approximate policy and 4.50 hours under the myopic policy, a difference of 1.75 hours.

\subsubsection{Decision-making involving OTs}
\label{sec:ot_experiments}
In previous numerical experiments, assignment decisions are made for reserve driver rosters consisting of XBs only. Model performances are now examined with OTs added to the decision space. The size of the OT workforce and their report times are additional inputs to the problem. \tcb{For this experiment, we specify ten OT report times: 06:00, 07:00, 08:00, 15:00, 16:00, 16:30, 17:00, 17:30, 18:00, and 19:00. These place three OTs in the morning and seven in the afternoon and evening, with more closely spaced reports during the PM peak.} \tcb{The minimum-pay duration is one hour ($\delot=4$ periods), equal to the shortest work-piece duration ($\delmin=4$). Pieces within each source do not overlap. Consequently, the general OT recursion in \eqref{eqn:ot_general} reduces to the scalar form in Algorithm~\ref{alg:bdpot}; the source scores in Algorithm~\ref{alg:loss_ot} likewise reduce to first-piece rewards for sources not yet revealed.} The per-period cost coefficient for OTs, $\costot$, \tcb{represents the incremental cost of an OT assignment beyond the guaranteed minimum-pay duration in the model objective}. We test three different values of $\costot$: $\{0,0.1,0.2\}$. Because each period is 0.25 hours long, these $\costot$ values correspond to \tcb{incremental penalties equal to} 0\%, 40\%, and 80\% of the per-period reward for covering an open work piece, respectively.

Numerical experiments are conducted for varying $\noxb$ values at base $p_h$ levels. The average performance in terms of total reward (open work hours covered minus \tcb{incremental} overtime costs) is reported for the PI, the approximate policy, both XB-first and OT-first myopic policies (denoted by ``Myo-XB'' and ``Myo-OT'', respectively), and the practice-oriented method. \tcb{Table~\ref{tab:xb_ot_hours} reports the SD of total reward for PI alongside the mean for each method.}

\begin{table}[ht!]
    \centering
    \caption{\tcb{Net reward: open work hours covered minus incremental OT costs} (XBs + 10 OTs; base $p_h$ levels)}
    \label{tab:xb_ot_hours}
    \resizebox{0.75\textwidth}{!}{
    \begin{tabular}{| c | c | c | c | c | c | c |}
    \hline
    $N_\texttt{XB}$ & $\costot$ &PI (SD) &Approx. [\%gap] &Myo-XB {[}\%gap{]}&Myo-OT {[}\%gap{]} &Pract. {[}\%gap{]} \\
    \hline
    6  & 0   & 52.22 (8.13) & 49.91 {[}4.42\%{]} & 48.24 {[}7.63\%{]} & 47.48 {[}9.08\%{]}  & \tcb{44.39} {[}\tcb{14.99}\%{]} \\
    6  & 0.1 & 47.61 (6.\tcb{7}3) & 45.17 {[}5.11\%{]} & 44.34 {[}6.86\%{]} & 43.10 {[}9.48\%{]}  & \tcb{39.47} {[}\tcb{17.10}\%{]} \\
    6  & 0.2 & 43.47 (5.64) & 41.48 {[}4.58\%{]} & 40.50 {[}6.83\%{]} & 38.56 {[}11.29\%{]} & \tcb{34.55} {[}\tcb{20.53}\%{]} \\
    \hline
    8  & 0   & 59.75 (10.05) & 57.50 {[}3.78\%{]} & 55.77 {[}6.67\%{]} & 55.06 {[}7.85\%{]}  & \tcb{51.74} {[}\tcb{13.41}\%{]} \\
    8  & 0.1 & 55.84 (8.64) & 53.10 {[}4.91\%{]} & 52.25 {[}6.43\%{]} & 50.65 {[}9.28\%{]}  & \tcb{47.08} {[}\tcb{15.68}\%{]} \\
    8  & 0.2 & 52.35 (7.53) & 49.82 {[}4.84\%{]} & 48.70 {[}6.98\%{]} & 46.17 {[}11.80\%{]} & \tcb{42.43} {[}\tcb{18.95}\%{]} \\
    \hline
    10 & 0   & 65.41 (12.09) & 63.31 {[}3.21\%{]} & 61.66 {[}5.73\%{]} & 61.16 {[}6.50\%{]}  & \tcb{57.33} {[}\tcb{12.35}\%{]} \\
    10 & 0.1 & 62.09 (10.73) & 59.19 {[}4.67\%{]} & 58.44 {[}5.87\%{]} & 56.78 {[}8.56\%{]}  & \tcb{52.99} {[}\tcb{14.66}\%{]} \\
    10 & 0.2 & 59.09 (9.65) & 56.34 {[}4.65\%{]} & 55.19 {[}6.61\%{]} & 52.31 {[}11.48\%{]} & \tcb{48.65} {[}\tcb{17.68}\%{]} \\
    \hline
    12 & 0   & 69.13 (13.94) & 67.56 {[}2.27\%{]} & 66.08 {[}4.41\%{]} & 65.74 {[}4.90\%{]}  & \tcb{62.26} {[}\tcb{9.94}\%{]} \\
    12 & 0.1 & 66.75 (12.72) & 64.09 {[}3.99\%{]} & 63.47 {[}4.92\%{]} & 61.33 {[}8.13\%{]}  & \tcb{58.38} {[}\tcb{12.54}\%{]} \\
    12 & 0.2 & 64.61 (11.72) & 62.08 {[}3.91\%{]} & 60.88 {[}5.78\%{]} & 56.90 {[}11.93\%{]} & \tcb{54.51} {[}\tcb{15.63}\%{]} \\
    \hline
    14 & 0   & 71.73 (15.30) & 70.60 {[}1.57\%{]} & 69.28 {[}3.41\%{]} & 69.06 {[}3.72\%{]}  & \tcb{65.71} {[}\tcb{8.39}\%{]}  \\
    14 & 0.1 & 70.02 (14.30) & 67.63 {[}3.42\%{]} & 67.12 {[}4.15\%{]} & 64.58 {[}7.77\%{]}  & \tcb{62.26} {[}\tcb{11.08}\%{]} \\
    14 & 0.2 & 68.47 (13.47) & 66.16 {[}3.37\%{]} & 64.92 {[}5.18\%{]} & 60.17 {[}12.12\%{]} & \tcb{58.82} {[}\tcb{14.09}\%{]} \\
    \hline
    \end{tabular}
    }
\end{table}

In every instance examined, the tested methods show the same ordering in terms of average performance: PI > Approximate policy > Myo-XB > Myo-OT > Practice-oriented method. \tcb{Net rewards for the OT-first myopic policy and the practice-oriented method decrease as} $\costot$ \tcb{increases, reflecting the higher penalty assigned to overtime use in the objective}. \tcb{The largest absolute mean-reward advantages} of the approximate policy over Myo-XB are observed for instances with $\costot=0$, which also coincide with the smallest gaps from the PI solution. \tcb{At $\costot=0$, the PI CV is higher than in the corresponding XB-only PI setting for each tested $\noxb$.}

\subsubsection{\tcb{Passenger-delay-weighted coverage rewards}}
\label{sec:passenger_weights}

In previous test instances, we used the duration of a work piece as its coverage reward ($c_k$). That \tcb{made the total reward interpretable as hours of open work covered}. \tcb{To examine a setting in which the relative importance of trips varies, we also use passenger-delay-weighted coverage rewards. These rewards provide a proxy for the passenger impact of leaving work uncovered.}

\tcb{As a weighting proxy, we assume that leaving a trip uncovered delays its riders by one scheduled headway, without accounting for alternative service, switching modes, or crowding.} \tcb{For work piece $k$, we therefore set} $c_k=\text{Round trip ridership} \times \text{Headway}$\tcb{, with headway measured in hours.}

The \tcb{specified} headway varies throughout the day. For each of the analyzed bus routes at MiWay, a constant headway is \tcb{used in the case-study input} for each operation period (e.g., AM peak, midday, late evening, etc.). \tcb{These route- and time-dependent headways reflect the service patterns used to construct the instance.} For example, the pattern of operation periods and corresponding headways for \tcb{Route 11} is summarized in Table \ref{tab:hways}.

\begin{table}[ht!]
    \centering
    \caption{Example operation periods and headways \tcb{for Route 11}}
    \label{tab:hways}
    \resizebox{0.505\textwidth}{!}{
    \begin{tabular}{| c | c | c |}
    \hline
    Operation period &Start and end time &Headway (min)  \\
    \hline
    Early morning & 4:00 AM -- 5:00 AM & 30 \\
    AM peak & 5:00 AM -- 9:00 AM & 15 \\
    Midday & 9:00 AM -- 3:00 PM & 20 \\
    PM peak & 3:00 PM -- 7:00 PM & 15 \\
    Evening & 7:00 PM -- 9:00 PM & 20 \\
    Late evening & 9:00 PM -- 3:30 AM & 30 \\
    \hline
    \end{tabular}
    }
\end{table}

Round trip ridership values are synthesized for each work piece by route and by operation period\tcb{; trip-level ridership observations were not available for the case study}. \tcb{A historical ranking of 2008 daily boardings \citep{mississauga2009ridership} provides context for the relative importance of the analyzed routes, while the numerical input specifies ridership and headway lookups by route and time of day. Each work piece receives the lookup values for its route at its start period. Its reward is the assigned ridership multiplied by the assigned headway in minutes and divided by 60, giving a passenger-hour weighting.}
With the newly defined weights for all work pieces, the total reward gained can be interpreted as \tcb{passenger-delay-weighted work coverage}. \tcb{We compare the four methods across} various $\noxb$ levels. The results for base and medium $p_h$ levels are presented in Tables \ref{tab:xb_waitt} and \ref{tab:xb_waitt_med}.

\begin{table}[ht!]
    \centering
    \caption{\tcb{Passenger-delay-weighted coverage reward, in passenger-hours} (XBs only; base $p_h$ levels)}
    \label{tab:xb_waitt}
    \resizebox{0.75\textwidth}{!}{
    \begin{tabular}{| c | c | c | c | c |}
    \hline
    $N_\texttt{XB}$ &Perfect info. (SD) &Approx. policy [\%gap] &Myopic policy {[}\%gap{]} &Pract. method {[}\%gap{]} \\
    \hline
    6  & 425.11 (51.20) & 405.37 {[}4.64\%{]} & 393.22 {[}7.50\%{]} & 333.32 {[}21.59\%{]} \\
    8  & 541.64 (71.39) & 513.52 {[}5.19\%{]} & 499.50 {[}7.78\%{]} & 432.67 {[}20.12\%{]} \\
    10 & 628.52 (100.56) & 597.96 {[}4.86\%{]} & 583.04 {[}7.24\%{]} & 511.79 {[}18.57\%{]} \\
    12 & 707.19 (123.97) & 678.32 {[}4.08\%{]} & 662.27 {[}6.35\%{]} & 590.60 {[}16.49\%{]} \\
    14 & 759.40 (145.66) & 732.42 {[}3.55\%{]} & 716.70 {[}5.62\%{]} & 648.47 {[}14.61\%{]} \\
    16 & 792.30 (160.87) & 767.96 {[}3.07\%{]} & 755.65 {[}4.63\%{]} & 692.91 {[}12.54\%{]} \\
    18 & 823.53 (178.35) & 803.50 {[}2.43\%{]} & 795.35 {[}3.42\%{]} & 734.57 {[}10.80\%{]} \\
    \hline
    \end{tabular}
    }
\end{table}

\begin{table}[ht!]
    \centering
    \caption{\tcb{Passenger-delay-weighted coverage reward, in passenger-hours} (XBs only; medium $p_h$ levels)}
    \label{tab:xb_waitt_med}
    \resizebox{0.75\textwidth}{!}{
    \begin{tabular}{| c | c | c | c | c |}
    \hline
    $N_\texttt{XB}$ &Perfect info. (SD) &Approx. policy [\%gap] &Myopic policy {[}\%gap{]} &Pract. method {[}\%gap{]}  \\
    \hline
    6  & 477.16 (36.62)  & 457.35 {[}4.15\%{]}  & 434.63 {[}8.91\%{]}  & 378.27 {[}20.73\%{]} \\
    8  & 621.96 (47.80)  & 592.97 {[}4.66\%{]}  & 565.03 {[}9.15\%{]}  & 495.31 {[}20.36\%{]} \\
    10 & 746.49 (75.30) & 711.54 {[}4.68\%{]}  & 681.64 {[}8.69\%{]}  & 604.13 {[}19.07\%{]} \\
    12 & 869.52 (91.96)  & 829.06 {[}4.65\%{]}  & 796.95 {[}8.35\%{]}  & 713.53 {[}17.94\%{]} \\
    14 & 963.69 (114.66)  & 919.89 {[}4.54\%{]}  & 890.53 {[}7.59\%{]}  & 802.47 {[}16.73\%{]} \\
    16 & 1030.86 (135.18) & 986.44 {[}4.31\%{]}  & 959.61 {[}6.91\%{]}  & 873.84 {[}15.23\%{]} \\
    18 & 1108.76 (155.78) & 1067.26 {[}3.74\%{]} & 1040.27 {[}6.18\%{]} & 953.77 {[}13.98\%{]} \\
    \hline
    \end{tabular}
    }
\end{table}

The performance ranking among the four methods remains the same from previous test settings. Gaps from PI are generally smaller as $\noxb$ increases \tcb{in these tests}. Compared with tests involving ordinary $c_k$ weights, the gap values are more sensitive to $\noxb$, especially for the practice-oriented method. \tcb{At matching openness levels, the approximate policy has a larger percentage advantage over the myopic policy under passenger-delay-weighted rewards than under ordinary rewards for most tested roster sizes.} \tcb{The myopic rule does not prioritize pieces by these weights, while the approximate policy scores known pieces using their specified rewards.} Greater gaps to the PI solution performance are also observed \tcb{in most tested cases} when compared with the results from ordinary $c_k$.

\subsection{Synthetic Case Study}
\label{sec:SyntheticCaseStudy}
To assess \tcb{performance on a larger constructed instance} beyond the MiWay instance, we also tested the proposed method on a synthetic schedule with 13 routes, \tcb{300 regular-duty sessions, and 2,130 productive operator-hours}. The full case-study description and detailed numerical results are provided in \tcb{Appendix}~\ref{sec:app_synthetic_results}. The synthetic results \tcb{show the reported approximate-policy means above the myopic and practice-oriented means in the tested settings. The PI benchmark optimizes over the complete eligible assignment set and provides an upper bound on the reward attainable by the online policies on each evaluated realization.}

Table~\ref{tab:synthetic_summary} summarizes the \tcb{reported comparisons}. Under ordinary reward coefficients, the approximate policy \tcb{has relatively small gaps to the complete PI benchmark} for most tested \tcb{roster} sizes, whereas the myopic policy often exhibits \tcb{larger gaps} in capacity-constrained settings. Under \tcb{passenger-delay-weighted} reward coefficients, \tcb{the reported difference between the approximate and myopic policies is larger for smaller rosters}. Because route-level rewards are highly heterogeneous in this setting, the myopic policy may assign scarce XBs to low-value work pieces, while the approximate policy \tcb{scores both known work and estimated future opportunities}.

\begin{table}[ht!]
\small
\centering
\caption{\tcb{Summary of synthetic case study results}}
\label{tab:synthetic_summary}
{\setlength{\tabcolsep}{4pt}
\renewcommand{\arraystretch}{1.2}
\begin{tabular}{@{}L{0.225\textwidth} *{3}{>{\centering\arraybackslash}p{0.15\textwidth}} L{0.244\textwidth}@{}}
\hline
Setting & \makecell[c]{Approx. gap\\to complete PI} & \makecell[c]{Myopic gap\\to complete PI} & \makecell[c]{Practice gap\\to complete PI} & \makecell[l]{Main\\observation} \\
\hline
\makecell[l]{Ordinary rewards,\\base $p_h$}
& $\approx$ 1.5--3.7\% 
& $\approx$ 3--8\% 
& $\approx$ 11--17\% 
& Small gap to complete PI \\

\makecell[l]{Ordinary rewards,\\medium $p_h$}
  & $\approx$ 2.4--3.4\% 
& $\approx$ 4.5--8.1\% 
& $\approx$ 12--15\% 
& Same ranking under higher openness \\

Passenger-delay-weighted rewards
& $\approx$ 4--11\% 
& $\approx$ 5--38\% 
& $\approx$ 9--39\% 
& Larger policy differences at smaller rosters \\
\hline
\end{tabular}
}
\end{table}

\section{Conclusion}
\label{sec:conclusion}
This study investigates the real-time assignment problem for extraboard transit operators. The problem is formulated as a Markov decision process (MDP) which is designed to capture its sequential and stochastic nature. The proposed MDP features a very large state space due to the heterogeneity in driver types and shifts, making it difficult to find an optimal policy for the MDP. To tackle this, an approximate \tcb{dynamic programming (ADP)} policy is proposed in the form of an \tcb{integer program}. \tcb{The integer program scores known open work together with estimated future rewards, then implements selected XB sequences beginning at the current decision epoch and selected OT pieces starting then.} The \tcb{individual future-value estimates are computed offline by} customized backward dynamic programming \tcb{recursions using specified work-source opening probabilities}. Then, a separate, simulation-based algorithm \tcb{estimates coefficients for an assignment-dependent competition adjustment used in the online integer program}.

Numerical analyses are conducted using two case studies: one based on \tcb{MiWay's published service schedules and constructed operator duties}; the other \tcb{is a larger synthetic instance}. The approximate policy is compared against \tcb{a myopic policy, a practice-oriented rule, and a scenario-based rolling-horizon method in selected MiWay tests}. The approximate policy \tcb{achieves higher mean reward than the tested myopic, practice-oriented, and scenario-based rolling-horizon benchmarks}. \tcb{Both case studies include perfect-information benchmarks for evaluating the approximate policy.} 

Numerical experiments \tcb{also include OTs in the MiWay rosters}. \tcb{In the tested settings, the approximate policy retains the highest mean net reward among the online methods, and the myopic rule prioritizing XBs over OTs performs better than its OT-first counterpart.} Lastly, additional computational analyses \tcb{use passenger-delay-weighted coverage rewards as a proxy for the service impact of uncovered work}. \tcb{In the weighted-reward experiments, the approximate policy's percentage advantage over the myopic rule decreases as XB roster size increases, indicating that assignment choices are particularly consequential when capacity is scarce and work pieces differ in passenger impact.} \tcb{Matched component comparisons show gains from individual future values in all six tested settings. The competition adjustment adds gains in the tested OT settings and larger weighted-reward XB rosters, with little change in the ordinary-reward XB-only settings.}

We recognize some limitations of this study, which suggest possible directions for future research. \tcb{When an XB is assigned a sequence beginning at the current epoch, any later pieces in that sequence become binding. This limits the ability to reassign the XB if additional work opens before the sequence ends; allowing such commitments to be revisited could provide greater flexibility.} Another limitation is that \tcb{coverage of each open work piece must begin at its specified start time}. For certain routes, it may be possible to cover a trip late by employing short-turn or skip-stop operations, depending on operational policies of the agency. \tcb{Incorporating these options would enlarge the decision and state spaces and may call for different value-approximation methods.} Lastly, this study addresses the real-time decisions for extraboard operators in isolation from other decision levels in extraboard scheduling. Upper level decisions are used as input to the real-time problem. In particular, extraboard report times from operation\tcb{al}-level decisions directly affect the performance of real-time assignment of extraboards. It would be interesting to integrate these two decision levels and possibly improve upon existing report time scheduling methods with the knowledge from this study.

\section*{Data availability}
The data will be made available on a reasonable request.

\section*{Declaration of competing interest}
The authors declare that they have no known competing financial interests or personal relationships that could have appeared to influence the work reported in this paper. 

\small
\onehalfspacing
\bibliographystyle{plainnat}

\appendix
\renewcommand{\thefigure}{A\arabic{figure}}
\setcounter{figure}{0}
\renewcommand{\thetable}{A\arabic{table}}
\setcounter{table}{0}



\section{Detailed Results with the Synthetic Case Study}
\label{sec:app_synthetic_results}
We used a synthetic workday schedule to analyze model performances at a larger scale. The synthetic case contains 13 regular bus routes with daily service patterns inspired by those of MiWay bus routes, with a workday still defined as 4:00 AM to 3:30 AM on the next day. The synthesized bus routes are designed to be heterogeneous in many aspects. Round trip duration ranges from \tcb{1.5} to 3.75 hours depending on route and time-of-the-day. Some bus routes have a flat temporal profile with a uniform service frequency throughout the workday; while some have increased service during AM and PM peak hours. Two of the 13 routes operate during peak hours only. 

\begin{figure}[!ht]
    \centering
    \includegraphics[width=0.65\textwidth]{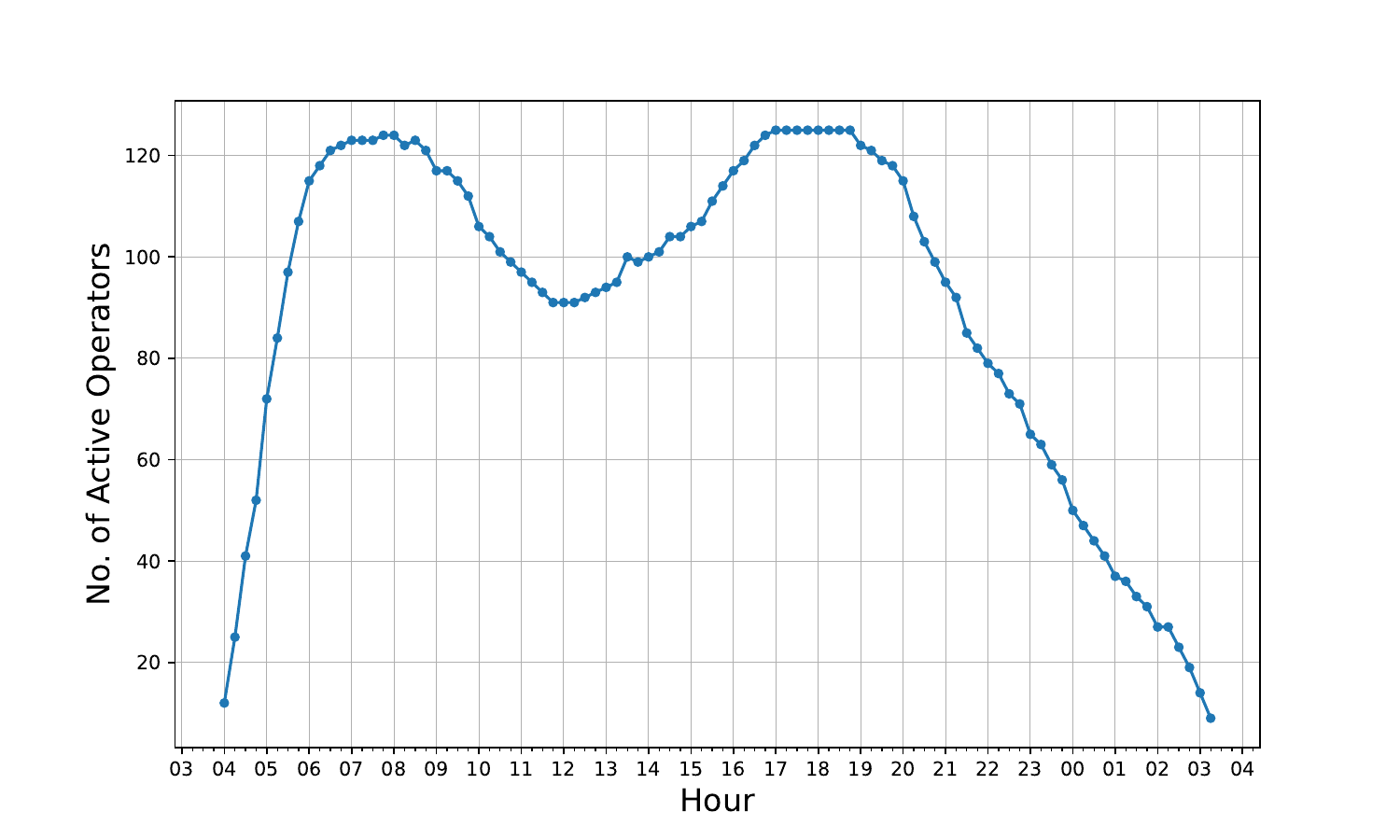}
    \caption{Daily Temporal Profile of Active Operators for the Synthesized Schedule}
    \label{fig:temp_profile_syn}
\end{figure}

There are \tcb{300 regular-duty work sessions} with various shift patterns. \tcb{Straight sessions of 7.5 hours are the most common, and 9.5-hour sessions are also included.} Split \tcb{sessions are included} for routes with peaked service patterns, similar to those at MiWay. Lastly, routes that only serve peak hours \tcb{have shorter straight sessions of 5.5--6 hours}. In total, there are \tcb{2,130 productive operator-hours} scheduled for a workday. The daily temporal profile for the synthesized schedule is shown in Figure \ref{fig:temp_profile_syn}.

As with the MiWay case study, scheduled \tcb{regular-duty sessions} make up type (1) work sources. The openness probability for these work sources ($p_h$) \tcb{represents a specified absence probability}. The base $p_h$ value is \tcb{7\% for 250 sessions and 10\% for 50 sessions}. Type (2) work sources, or extra trips dispatched against unexpected events, are defined for four bus routes that have peaked service patterns. $p_h$ values for type (2) sources are \tcb{specified by route and start time, ranging from 2\% to 18\% at the base level}.

Numerical experiments are conducted to compare performances among the four aforementioned methods: PI, the approximate policy, the myopic policy, and the practice-oriented method. \tcb{For the synthetic PI benchmark, we solve a binary assignment model with an assignment variable for each eligible XB--piece pair, at-most-once coverage for each piece, and per-period nonoverlap for each XB. This model includes all feasible nonoverlapping assignments. Across the three synthetic result settings and their 24 roster configurations, all 24,000 path-level PI models were solved to optimality (maximum recorded relative gap below $10^{-10}$).} The tests involve reserve operator rosters made of XBs only, with roster sizes $\noxb\tcb{\in}\{5,10,15,\hdots,40\}$ for base $p_h$ and $\noxb\tcb{\in}\{15,20,25,\hdots,50\}$ for medium $p_h$. XB report times are again determined using the method proposed by \citep{kaysi1990scheduling}. \tcb{The reported online-policy experiments were run} on the same PC as before. Slight modifications are made to the Gurobi model parameters for the \tcb{online} integer program. The model is set to ``MIPFocus'' such that finding a feasible integer solution is prioritized. \tcb{The PI models used one solver thread and zero requested MIP gap.}

We first examine the ordinary reward setting: the reward for covering an open work piece is equal to the piece's duration in hours. Tables \ref{tab:xb_hours_syn} and \ref{tab:uow_hours_syn} summarize average open work hours covered by XBs and the amount of uncovered open work left, respectively, at base $p_h$ levels. \tcb{The means have the same ranking as in the MiWay comparison (complete PI $>$ approximate policy $>$ myopic policy $>$ practice-oriented method). The approximate policy has smaller gaps to complete PI than the myopic policy across the tested rosters.} Another interesting observation is that, starting at $\noxb=30$, performance gaps shrink \tcb{for the approximate and myopic policies}. These high $\noxb$ values correspond to relatively abundant supply of XBs such that there is \tcb{less} uncovered open work left (from Table \ref{tab:uow_hours_syn}).

\begin{table}[ht!]
    \centering
    \caption{Open work hours covered (synthetic case study; XBs only; base $p_h$ levels)}
    \label{tab:xb_hours_syn}
    \resizebox{0.75\textwidth}{!}{
    \begin{tabular}{| c | c | c | c | c |}
    \hline
    $N_\texttt{XB}$ &\tcb{Complete PI} &Approx. policy &Myopic policy& Pract. method  \\
    \hline
    5  & \tcb{35.91} (SD = \tcb{0.87})  & 34.86 (SD = 1.19)  & 33.13 (SD = 1.47)  & 30.19 (SD = 2.23)  \\
    && {[}gap = \tcb{2.91}\%{]}& {[}gap = \tcb{7.73}\%{]}& {[}gap = \tcb{15.92}\%{]} \\ \hline
    10 & \tcb{71.06} (SD = \tcb{1.73})  & 68.64 (SD = 2.26)  & 65.48 (SD = 2.41)  & 59.17 (SD = 3.51)  \\
    && {[}gap = \tcb{3.41}\%{]}& {[}gap = \tcb{7.86}\%{]}& {[}gap = \tcb{16.74}\%{]} \\ \hline
    15 & \tcb{104.58} (SD = \tcb{3.50}) & 100.82 (SD = 4.07) & 96.20 (SD = 3.92)  & 86.66 (SD = 5.11)  \\
    && {[}gap = \tcb{3.59}\%{]}& {[}gap = \tcb{8.01}\%{]}& {[}gap = \tcb{17.13}\%{]} \\ \hline
    20 & \tcb{135.67} (SD = \tcb{6.79}) & 130.70 (SD = 7.06) & 125.10 (SD = 6.34) & 112.57 (SD = 7.08)  \\
    && {[}gap = \tcb{3.66}\%{]}& {[}gap = \tcb{7.79}\%{]}& {[}gap = \tcb{17.02}\%{]} \\ \hline
    25 & \tcb{163.05} (SD = \tcb{11.66}) & 157.34 (SD = 11.07) & 151.23 (SD = 9.98) & 136.24 (SD = 9.64) \\
    && {[}gap = \tcb{3.50}\%{]}& {[}gap = \tcb{7.25}\%{]}& {[}gap = \tcb{16.44}\%{]} \\ \hline
    30 & \tcb{183.34} (SD = \tcb{17.43}) & 177.87 (SD = 16.12) & 172.14 (SD = 14.49) & 155.93 (SD = 12.88)  \\
    && {[}gap = \tcb{2.98}\%{]}& {[}gap = \tcb{6.11}\%{]}& {[}gap = \tcb{14.95}\%{]} \\ \hline
    35 & \tcb{196.00} (SD = \tcb{22.75}) & 191.40 (SD = 21.04) & 187.11 (SD = 19.29) & 170.69 (SD = 16.59)  \\
    && {[}gap = \tcb{2.35}\%{]}& {[}gap = \tcb{4.54}\%{]}& {[}gap = \tcb{12.92}\%{]} \\ \hline
    40 & \tcb{204.77} (SD = \tcb{27.48}) & 201.67 (SD = 25.86) & 198.38 (SD = 24.11) & 182.64 (SD = 19.86)  \\
    && {[}gap = 1.51\%{]}& {[}gap = \tcb{3.12}\%{]}& {[}gap = \tcb{10.81}\%{]} \\
    \hline
    \end{tabular}
    }
\end{table}

\begin{table}[ht!]
    \centering
    \caption{Uncovered open work hours and percentage out of total open work (synthetic case study; XBs only; base $p_h$ levels)}
    \label{tab:uow_hours_syn}
    \resizebox{0.705\textwidth}{!}{
    \begin{tabular}{| c | c | c | c | c |}
    \hline
    $N_\texttt{XB}$ &\tcb{Complete PI} &Approx. policy &Myopic policy& Pract. method  \\
    \hline
    5  & \tcb{176.78} (\tcb{83.12}\%) & 177.82 (83.61\%) & 179.55 (84.42\%) & 182.49 (85.80\%) \\
    10 & \tcb{141.62} (\tcb{66.59}\%) & 144.05 (67.73\%) & 147.20 (69.21\%) & 153.51 (72.18\%) \\
    15 & \tcb{108.10} (\tcb{50.83}\%) & 111.87 (52.60\%) & 116.48 (54.77\%) & 126.02 (59.25\%) \\
    20 & \tcb{77.02} (\tcb{36.21}\%)  & 81.99 (38.55\%)  & 87.58 (41.18\%)  & 100.11 (47.07\%) \\
    25 & \tcb{49.64} (\tcb{23.34}\%)  & 55.34 (26.02\%)  & 61.45 (28.89\%)  & 76.44 (35.94\%)  \\
    30 & \tcb{29.34} (\tcb{13.80}\%)  & 34.81 (16.37\%)  & 40.54 (19.06\%)  & 56.75 (26.68\%)  \\
    35 & \tcb{16.68} (\tcb{7.84}\%)   & 21.28 (10.01\%)  & 25.57 (12.02\%)  & 41.99 (19.74\%)  \\
    40 & \tcb{7.92} (\tcb{3.72}\%)    & 11.02 (5.18\%)   & 14.30 (6.72\%)   & 30.05 (14.13\%)  \\
    \hline
    \end{tabular}
    }
\end{table}

The average XB utilization rates are reported in Table \ref{tab:xb_utilize_syn}. As explained before, the average utilization rate decreases as $\noxb$ \tcb{increases}. However, these rates decrease slowly for small $\noxb$ numbers and remain at a relatively high level (e.g., $> 90\%$ for \tcb{complete PI} when $\noxb<20$). \tcb{At larger roster sizes, utilization falls more quickly, consistent with XB supply becoming less scarce relative to the sampled open work.}

\begin{table}[ht!]
    \centering
    \caption{Average XB utilization rate (synthetic case study; XBs only; base $p_h$ levels)}
    \label{tab:xb_utilize_syn}
    \resizebox{0.65\textwidth}{!}{
    \begin{tabular}{| c | c | c | c | c |}
    \hline
    $N_\texttt{XB}$ &\tcb{Complete PI} &Approx. policy &Myopic policy& Pract. method  \\
    \hline
    5  & \tcb{95.75}\% & 92.96\% & 88.35\% & 80.52\% \\
    10 & \tcb{94.75}\% & 91.51\% & 87.31\% & 78.90\% \\
    15 & \tcb{92.96}\% & 89.61\% & 85.51\% & 77.03\% \\
    20 & \tcb{90.44}\% & 87.13\% & 83.40\% & 75.05\% \\
    25 & \tcb{86.96}\% & 83.91\% & 80.66\% & 72.66\% \\
    30 & \tcb{81.48}\% & 79.05\% & 76.51\% & 69.30\% \\
    35 & \tcb{74.67}\% & 72.91\% & 71.28\% & 65.03\% \\
    40 & \tcb{68.26}\% & 67.22\% & 66.13\% & 60.88\% \\
    \hline
    \end{tabular}
    }
\end{table}

The same numerical experiments are run for the medium $p_h$ levels, which are 1.5 times the base values. Performance comparisons among the four methods are reported in Tables \ref{tab:xb_hours_med_syn} and \ref{tab:uow_hours_med_syn}. The observed gaps exhibit similar patterns as those of the base $p_h$ case. \tcb{The approximate policy has smaller gaps to complete PI than the other online methods.} \tcb{At larger rosters, its gap and the myopic policy's gap narrow as the amount of uncovered work falls.}

\begin{table}[ht!]
    \centering
    \caption{Open work hours covered (synthetic case study; XBs only; medium $p_h$ levels)}
    \label{tab:xb_hours_med_syn}
    \resizebox{0.75\textwidth}{!}{
    \begin{tabular}{| c | c | c | c | c |}
    \hline
    $N_\texttt{XB}$ &\tcb{Complete PI} &Approx. policy &Myopic policy& Pract. method \\
    \hline
    15 & \tcb{109.37} (SD = \tcb{1.59}) & 106.67 (SD = 2.29) & 100.61 (SD = 2.53) & 93.80 (SD = 3.74)  \\
    && {[}gap = \tcb{2.47}\%{]}& {[}gap = \tcb{8.01}\%{]}& {[}gap = \tcb{14.23}\%{]} \\ \hline
    20 & \tcb{144.91} (SD = \tcb{2.64}) & 141.09 (SD = 3.52) & 133.31 (SD = 3.48) & 123.65 (SD = 4.96) \\
    && {[}gap = \tcb{2.64}\%{]}& {[}gap = \tcb{8.01}\%{]}& {[}gap = \tcb{14.67}\%{]} \\ \hline
    25 & \tcb{179.65} (SD = \tcb{4.14}) & 174.49 (SD = 5.03) & 165.12 (SD = 4.67) & 152.69 (SD = 6.14) \\
    && {[}gap = \tcb{2.87}\%{]}& {[}gap = \tcb{8.09}\%{]}& {[}gap = \tcb{15.01}\%{]} \\ \hline
    30 & \tcb{211.54} (SD = \tcb{6.93}) & 204.87 (SD = 7.57) & 194.95 (SD = 6.55) & 179.52 (SD = 8.19) \\
    && {[}gap = \tcb{3.16}\%{]}& {[}gap = \tcb{7.84}\%{]}& {[}gap = \tcb{15.14}\%{]} \\ \hline
    35 & \tcb{238.45} (SD = \tcb{11.03}) & 230.49 (SD = 11.14) & 221.09 (SD = 9.66) & 202.83 (SD = 10.42) \\
    && {[}gap = \tcb{3.34}\%{]}& {[}gap = \tcb{7.28}\%{]}& {[}gap = \tcb{14.94}\%{]} \\ \hline
    40 & \tcb{263.38} (SD = \tcb{15.58}) & 254.96 (SD = 14.94) & 245.91 (SD = 13.07) & 224.70 (SD = 12.72) \\
    && {[}gap = \tcb{3.20}\%{]}& {[}gap = \tcb{6.63}\%{]}& {[}gap = \tcb{14.68}\%{]} \\ \hline
    45 & \tcb{283.51} (SD = \tcb{21.00}) & 275.33 (SD = 19.27) & 267.54 (SD = 17.34) & 244.75 (SD = 15.43) \\
    && {[}gap = \tcb{2.88}\%{]}& {[}gap = \tcb{5.63}\%{]}& {[}gap = \tcb{13.67}\%{]} \\ \hline
    50 & \tcb{298.16} (SD = \tcb{26.47}) & 290.99 (SD = 24.24) & 284.55 (SD = 21.97) & 261.94 (SD = 18.68) \\
    && {[}gap = \tcb{2.41}\%{]}& {[}gap = \tcb{4.57}\%{]}& {[}gap = \tcb{12.15}\%{]} \\
    \hline
    \end{tabular}
    }
\end{table}

\begin{table}[ht!]
    \centering
    \caption{Uncovered open work hours and percentage out of total open work (synthetic case study; XBs only; medium $p_h$ levels)}
    \label{tab:uow_hours_med_syn}
    \resizebox{0.75\textwidth}{!}{
    \begin{tabular}{| c | c | c | c | c |}
    \hline
    $N_\texttt{XB}$ &\tcb{Complete PI} &Approx. policy &Myopic policy& Pract. method  \\
    \hline
    15 & \tcb{211.71} (\tcb{65.94}\%) & 214.40 (66.78\%) & 220.46 (68.66\%) & 227.27 (70.78\%) \\
    20 & \tcb{176.16} (\tcb{54.87}\%) & 179.98 (56.06\%) & 187.77 (58.48\%) & 197.43 (61.49\%) \\
    25 & \tcb{141.43} (\tcb{44.05}\%) & 146.59 (45.66\%) & 155.96 (48.57\%) & 168.38 (52.44\%) \\
    30 & \tcb{109.53} (\tcb{34.11}\%) & 116.21 (36.19\%) & 126.12 (39.28\%) & 141.55 (44.09\%) \\
    35 & \tcb{82.62} (\tcb{25.73}\%)  & 90.58 (28.21\%)  & 99.99 (31.14\%)  & 118.25 (36.83\%) \\
    40 & \tcb{57.70} (\tcb{17.97}\%)  & 66.11 (20.59\%)  & 75.16 (23.41\%)  & 96.38 (30.02\%)  \\
    45 & \tcb{37.57} (\tcb{11.70}\%)  & 45.75 (14.25\%)  & 53.53 (16.67\%)  & 76.33 (23.77\%)  \\
    50 & \tcb{22.91} (\tcb{7.14}\%)   & 30.08 (9.37\%)   & 36.53 (11.38\%)  & 59.14 (18.42\%) \\
    \hline
    \end{tabular}
    }
\end{table}

Computational \tcb{performance is assessed through reported offline parameter-computation times and online integer-program solve times}, similar to the analyses in the MiWay case study. For medium $p_h$ levels, the \tcb{reported total offline time} for the four algorithms is around 52 minutes. \tcb{These estimates can be reused when their input settings do not change.} Examining the \tcb{reported} solution times for the integer program at $\noxb=50$ \tcb{shows that most recorded solves take less than 0.2 seconds. The reported maximum is slightly above 90 seconds, and approximately 0.0079\% of recorded solves reach the solver time limit.}


Next, we conduct numerical experiments \tcb{using passenger-delay-weighted coverage rewards as a proxy for the impact of uncovered work}. The approach to generating individual weights of work pieces is the same as the one described in the MiWay case study. For this synthetic schedule, \tcb{the specified ridership weights vary across routes; the highest route weight is three times the lowest}. This \tcb{constructed} variability is designed to examine how the approximate policy fares \tcb{when route rewards are heterogeneous}. 

The results for base $p_h$ levels are presented in Table \ref{tab:xb_waitt_syn}. \tcb{The policy gaps to complete PI are larger than in the ordinary-reward setting, especially for the myopic and practice-oriented methods at smaller rosters. One possible explanation is that these rules can assign scarce XBs to lower-reward pieces, leaving higher-reward work uncovered.} Under the approximate policy, however, XBs may be ordered to hold even when a low-value piece is open. \tcb{Its assignments are selected using a score for known open work and estimated future rewards.} As $\noxb$ increases, \tcb{the difference between the approximate and myopic policies narrows in these tests. With more XB capacity relative to open work, immediate assignment to a lower-value piece may be less likely to prevent coverage of later higher-value work.}

\setlength{\tabcolsep}{4pt} 
\renewcommand{\arraystretch}{1.15} 
\begin{table}[ht!]
    \centering
    \caption{\tcb{Passenger-delay-weighted coverage reward, in passenger-hours} (synthetic case study; XBs only; base $p_h$ levels)}
    \label{tab:xb_waitt_syn}
    \resizebox{0.75\textwidth}{!}{
    \begin{tabular}{| c | c | c | c | c |}
    \hline
    $N_\texttt{XB}$ &\tcb{Complete PI} (SD) &Approx. policy [\%gap] &Myopic policy {[}\%gap{]} &Pract. method {[}\%gap{]} \\
    \hline
    5  & \tcb{279.25} (\tcb{21.69})  & 255.66 {[}\tcb{8.45}\%{]}  & 172.86 {[}\tcb{38.10}\%{]} & 170.94 {[}\tcb{38.79}\%{]} \\
    10 & \tcb{525.43} (\tcb{42.96})  & 473.60 {[}\tcb{9.86}\%{]}  & 340.49 {[}\tcb{35.20}\%{]} & 339.03 {[}\tcb{35.48}\%{]} \\
    15 & \tcb{723.11} (68.15)  & 645.28 {[}\tcb{10.76}\%{]} & 503.18 {[}\tcb{30.41}\%{]} & 491.63 {[}\tcb{32.01}\%{]} \\
    20 & \tcb{872.50} (\tcb{93.40})  & 777.47 {[}\tcb{10.89}\%{]} & 657.27 {[}\tcb{24.67}\%{]} & 635.50 {[}\tcb{27.16}\%{]} \\
    25 & \tcb{983.43} (\tcb{117.33})  & 879.53 {[}\tcb{10.56}\%{]} & 797.35 {[}\tcb{18.92}\%{]} & 765.79 {[}\tcb{22.13}\%{]} \\
    30 & \tcb{1051.91} (\tcb{137.29}) & 976.43 {[}\tcb{7.18}\%{]}  & 913.58 {[}\tcb{13.15}\%{]} & 874.27 {[}\tcb{16.89}\%{]} \\
    35 & \tcb{1091.10} (\tcb{152.32}) & 1023.81 {[}\tcb{6.17}\%{]} & 999.10 {[}8.43\%{]}  & 950.64 {[}12.87\%{]} \\
    40 & \tcb{1116.40} (164.79) & 1068.07 {[}4.33\%{]} & 1060.49 {[}5.01\%{]} & 1012.82 {[}9.28\%{]} \\
    \hline
    \end{tabular}
    }
\end{table}

\section{\tcb{Additional component diagnostics}}
\label{sec:app_component_diagnostics}
\renewcommand{\thetable}{\thesection\arabic{table}}
\setcounter{table}{0}

To examine how the competition adjustment's contribution varies across XB rosters, we extend the XB-only comparisons in Section~\ref{sec:component_assessment} to $\noxb\in\{6,8,10,12,14,16,18\}$. For each roster, we compare individual values alone with individual values plus the adjustment under ordinary and passenger-delay-weighted rewards. Both variants use the same 1,000 base-probability opening paths, report times, feasible assignment decisions, and individual-value estimates. The reward-specific coefficient sets fitted at $\rthr=0.970$ are fixed throughout the comparison. The six- and 18-XB results are the same as those in Table~\ref{tab:component_diagnostic}.

Table~\ref{tab:adjustment_roster_sweep} reports the mean reward with the adjustment minus that without it. Each approximate paired 95\% interval is the mean difference plus or minus 1.96 times the sample standard deviation of the 1,000 pathwise differences divided by $\sqrt{1000}$. These descriptive, pointwise intervals quantify evaluation-path variability conditional on the fixed fitted coefficients.

\begin{table}[ht!]
\centering
\caption{\tcb{Incremental reward from the competition adjustment across XB rosters (base $p_h$; no OTs)}}
\label{tab:adjustment_roster_sweep}
{\small\setlength{\tabcolsep}{6pt}\renewcommand{\arraystretch}{1.15}
\begin{tabular}{@{}crrrr@{}}
\toprule
& \multicolumn{2}{c}{Ordinary reward (h)} & \multicolumn{2}{c}{Weighted reward (proxy passenger-h)} \\
\cmidrule(lr){2-3}\cmidrule(l){4-5}
$\noxb$ & Mean difference & Paired 95\% interval & Mean difference & Paired 95\% interval \\
\midrule
6  & $+0.0025$ & $[-0.0100,\ 0.0150]$ & $-0.040$ & $[-0.240,\ 0.160]$ \\
8  & $-0.0178$ & $[-0.0360,\ 0.0005]$ & $+0.393$ & $[0.075,\ 0.712]$ \\
10 & $-0.0085$ & $[-0.0305,\ 0.0135]$ & $+0.878$ & $[0.496,\ 1.259]$ \\
12 & $+0.0108$ & $[-0.0170,\ 0.0385]$ & $+2.802$ & $[2.196,\ 3.408]$ \\
14 & $+0.0293$ & $[-0.0042,\ 0.0627]$ & $+4.671$ & $[3.986,\ 5.357]$ \\
16 & $-0.0088$ & $[-0.0454,\ 0.0279]$ & $+5.454$ & $[4.728,\ 6.180]$ \\
18 & $-0.0350$ & $[-0.0817,\ 0.0117]$ & $+7.919$ & $[7.042,\ 8.797]$ \\
\bottomrule
\end{tabular}
}
\end{table}

Under ordinary rewards, mean differences range from $-0.035$ to $+0.029$ hours and all seven intervals include zero. Under weighted rewards, the difference is close to zero with six XBs, while gains increase across the tested rosters from eight to 18 XBs, with each corresponding interval above zero. This pattern is consistent with the adjustment becoming more consequential when more XBs compete for work with heterogeneous rewards. Because report-time profiles also change with roster size, this trend reflects the tested combinations of roster size and shift timing.

\end{document}